\documentclass{lmcs} 
\pdfoutput=1

\usepackage{hyperref}

\usepackage{tikz}
\usetikzlibrary{arrows,decorations,automata,backgrounds,positioning,calc,cd}
\usetikzlibrary{shapes.geometric,shapes.misc,matrix,arrows.meta,decorations.markings}
\usetikzlibrary{circuits.logic.US,patterns,shapes}
\tikzcdset{scale cd/.style={every label/.append style={scale=#1},
    cells={nodes={scale=#1}}}}
\usepackage{stmaryrd}
\usepackage{microtype}
\usepackage{tikz-cd}
\usepackage{quiver}
\usepackage{xfrac}
\usepackage{subcaption}
\usepackage{amsmath}

\usepackage{mathcommand,mathtools}
\usepackage[scope,quotation,notion,paper]{knowledge}
\knowledgeconfigure{protect quotation={tikzcd}}
\knowledgestyle{intro notion}{}
\knowledge{notion}
| prop
| props
| Prop
| Props
| coloured prop
| coloured props

\knowledge{notion}
| graded prop
| graded props

\knowledge{notion}
| equal up to regrading
| up to regrading

\knowledge{notion}
| SMC
| SMCs

\knowledge{notion}
| underlying SMC
| underlying SMCs
| underlying
| underlying prop
| underlying props

\knowledge{notion}
| -SMC
| graded SMC
| graded SMCs
| graded category
| graded categories
| graded symmetric monoidal category
| graded symmetric monoidal categories
| graded symmetric strict monoidal category
| graded symmetric strict monoidal categories

\knowledge{notion}
| para
| para construction

\knowledge{notion}
| graded functor
| graded symmetric monoidal functor
| graded symmetric strict monoidal functor
| graded functors
| graded symmetric monoidal functors

\knowledge{notion}
| regrade
| regrading
| regradings
| regraded
| Regrading

\knowledge{notion}
| actegory
| actegories
| Actegories
| action
| symmetric monoidal action

\knowledge{notion}
| theory@MON
| monoidal theory
| monoidal theories

\knowledge{notion}
| graded theory
| graded theories
| theory@GRD
| theories@GRD
| graded monoidal theory
| graded monoidal theories

\knowledge{notion}
| regrading collapse
| collapse

\knowledge{notion}
| flattened theory
| flattened theories
| flattening
| flattened generator
| flattened generators

\knowledge{notion}
| monoidal signature
| monoidal signatures
| monoidal term
| monoidal terms

\knowledge{notion}
| syntactic category@MON
| syntactic categories@MON

\knowledge{notion}
| syntactic category@GRD
| syntactic categories@GRD

\knowledge{notion}
| graded signature
| graded signatures
| graded term
| graded terms
| Graded terms
| term@GRD
| terms@GRD

\usepackage[capitalise,nameinlink]{cleveref}
\crefname{equation}{}{}
\crefname{equations}{}{}
\crefname{figure}{Figure}{Figures}
\crefname{defi}{Definition}{Definitions}
\crefname{defiC}{Definition}{Definitions}
\crefname{rem}{Remark}{Remarks}
\crefname{lem}{Lemma}{Lemmas}
\crefname{section}{Section}{Sections}
\crefname{app}{Appendix}{Appendices}
\crefname{thm}{Theorem}{Theorems}
\crefname{thmC}{Theorem}{Theorems}
\crefname{exa}{Example}{Examples}
\crefname{cor}{Corollary}{Corollaries}
\crefname{prop}{Proposition}{Propositions}
\crefname{propC}{Proposition}{Propositions}

\usepackage{bussproofs}
\newenvironment{bprooftree}
{\leavevmode\hbox\bgroup}
{\DisplayProof\egroup}
\usepackage{arydshln}

\usepackage{tikzit}
\ifdef{\macrotikzfig}{}{\newcommand{\macrotikzfig}[2]{\input{#1.tikz}}}

\tikzstyle{black}=[fill=black, draw=black, shape=circle, minimum size=4pt, inner sep=0pt, outer sep=0pt]
\tikzstyle{delgrade}=[fill={rgb,255: red,204; green,112; blue,0}, draw={rgb,255: red,204; green,112; blue,0}, shape=circle, minimum size=4pt, inner sep=0pt, outer sep=0pt]
\tikzstyle{conv}=[fill=none, draw=black, shape=rectangle, inner sep=3pt]
\tikzstyle{letter}=[fill=white, line width=0.5pt, inner sep=1pt, draw=black, rectangle, rounded corners, minimum height=1.75em, minimum width=1.75em]
\tikzstyle{gletter}=[fill=white, line width=0.5pt, inner sep=1pt, draw={rgb,255: red,204; green,112; blue,0}, text={rgb,255: red,204; green,112; blue,0}, regular polygon, regular polygon sides=4, rounded corners,minimum height=1.8em, minimum width={width("$f$")+3pt}]
\tikzstyle{dist}=[fill=none, draw=black, isosceles triangle, isosceles triangle apex angle=60, inner sep=3pt, rotate=180, anchor=center, shape border uses incircle]
\tikzstyle{white}=[fill=white, draw=black, shape=circle, minimum size=4pt, inner sep=0pt, outer sep=0pt, tikzit draw=black, tikzit fill=white]
\tikzstyle{basic box}=[draw, fill=white, rectangle, minimum height=1.2em, minimum width=1em]
\tikzstyle{ancilla}=[fill=black, draw=black, shape=rectangle, minimum width=0.01cm, minimum height=0.2cm, inner sep=0.01em]
\tikzstyle{and}=[fill=white, draw=black, and gate, anchor=center]
\tikzstyle{or}=[fill=white, draw=black, or gate, anchor=center]
\tikzstyle{not}=[fill=white, draw=black, not gate, anchor=center]
\tikzstyle{xor}=[fill=white, draw=black, xor gate, scale=.6]
\tikzstyle{flip}=[draw={rgb,255: red,86; green,86; blue,86}, fill={rgb,255: red,86; green,86; blue,86}, 
rounded rectangle, rounded rectangle right arc=none, text=white, node font={\small}, inner sep =0.2em]
\tikzstyle{flipop}=[draw={rgb,255: red,86; green,86; blue,86}, fill={rgb,255: red,86; green,86; blue,86}, 
rounded rectangle, rounded rectangle left arc=none, text=white, node font={\small}, inner sep =0.2em]
\tikzstyle{flipup}=[draw={rgb,255: red,204; green,112; blue,0}, fill={rgb,255: red,204; green,112; blue,0}, 
rounded rectangle, rounded rectangle right arc=none, rotate=90, text=black, node font={\small}, inner sep =0.2em]
\tikzstyle{scalar}=[draw=black,rounded rectangle, rounded rectangle left arc=none, node font={\small}, inner sep =0.2em]
\tikzstyle{negone}=[draw=black,rounded rectangle, rounded rectangle left arc=none, node font={\small}, inner sep =0.2em]
\tikzstyle{gnegone}=[draw={rgb,255: red,204; green,112; blue,0}, text={rgb,255: red,204; green,112; blue,0},rounded rectangle, rounded rectangle left arc=none, node font={\small}, inner sep =0.2em, rotate=90]
\tikzstyle{if}=[trapezium, draw=black, fill=white, minimum width=6pt, minimum height=1pt, inner sep=1pt, rotate=270]
\tikzstyle{myif}=[trapezium, trapezium angle=60, draw,inner xsep=0pt,outer sep=0pt,minimum height=20pt, rotate=270, text width=10pt, fill=white,label={center:\rotatebox{270}{$\mathtt{if}$}}]
\tikzstyle{oplus}=[
    draw={rgb,255: red,204; green,112; blue,0}, 
    circle,
    minimum size=0.7em,
    inner sep=0.1em,  
    path picture={
        \draw[line width=0.5pt]
        (path picture bounding box.south) -- (path picture bounding box.center);
        \draw[black,line width=0.5pt]
        (path picture bounding box.west) -- (path picture bounding box.east);
    }
]

\tikzstyle{nn}=[fill=gray, draw=gray, regular polygon, regular polygon sides=3, minimum width=0.4em, shape border rotate=90, inner sep=1.5pt]

\tikzstyle{single}=[-, draw=black,line width=0.5pt,line cap=rect]
\tikzstyle{double}=[-, draw=black, line width=1.2pt, tikzit fill=cyan, tikzit draw=cyan,line cap=rect]
\tikzstyle{boxes}=[-, fill=white, densely dotted, tikzit draw={rgb,255: red,17; green,167; blue,0}]
\tikzstyle{green}=[-, line width=0.5pt, draw={rgb,255: red,204; green,112; blue,0}, tikzit draw={rgb,255: red,204; green,112; blue,0},line cap=butt]
\tikzstyle{greendouble}=[-, draw={rgb,255: red,204; green,112; blue,0}, line width=1.2pt, tikzit draw={rgb,255: red,150; green,72; blue,0},line cap=butt,shorten <=-0.25pt,shorten >=-0.25pt]

\knowledgenewmathcommand{\N}{\cmdkl{\mathbb{N}}}
\knowledgenewmathcommand{\R}{\cmdkl{\mathbb{R}}}
\knowledgenewmathcommand{\B}{\cmdkl{\mathbb{B}}}
\knowledgenewmathcommand{\G}{\cmdkl{\mathbb{G}}}
\newrobustcmd{\id}{\mathrm{id}}
\newrobustcmd{\Hom}{\mathrm{Hom}}
\newrobustcmd{\dom}{\mathrm{dom}}
\newrobustcmd{\FSet}{\mathbf{FinSet}}
\newrobustcmd{\FStoch}{\mathbf{FinStoch}}
\knowledgenewmathcommand{\FInj}{\cmdkl{\mathbf{FinInj}}}
\newrobustcmd{\OrdInj}{\mathbf{OrdInj}}
\newrobustcmd{\FStochSurj}{\mathbf{FinStoch}_{\mathrm{Surj}}}
\knowledgenewmathcommand{\BStoch}{\cmdkl{\mathbf{BinStoch}}}
\newrobustcmd{\Set}{\mathbf{Set}}
\newrobustcmd{\ImP}{\mathbf{ImP}}
\knowledgenewmathcommand{\BImP}{\cmdkl{\mathbf{BImP}}}
\knowledgenewmathcommand{\BInt}{\cmdkl{\mathbf{Int}_{\mathbb{B}}}}
\knowledgenewmathcommand{\EPBS}{\cmdkl{\mathbf{EPBS}}}
\knowledgenewmathcommand{\Circ}{\cmdkl{\mathsf{CausCirc}}}
\newrobustcmd{\PCirc}{\mathsf{ProbCirc}}
\newrobustcmd{\GCirc}{\Circ_{\scriptscriptstyle i}}
\newrobustcmd{\Kl}{\mathbf{KL}}

\newrobustcmd{\sig}{\Sigma}
\newrobustcmd{\sigv}{\sig^{\mathsf{v}}}
\newrobustcmd{\axv}{E^{\mathsf{v}}}
\newrobustcmd{\axh}{E^{\raisebox{-0.3ex}{$\scriptstyle\mathsf{h}$}}}
\newrobustcmd{\sigC}{\sig^{\catC}}
\newrobustcmd{\axC}{E^{\catC}}
\newrobustcmd{\thC}{\th^{\catC}}
\newrobustcmd{\closed}[1]{\overline{#1}}
\newrobustcmd{\Para}[2]{\withkl{\kl[para]}{\cmdkl{\mathbf{Para}}(#1,#2)}}
\knowledgenewmathcommand{\with}{\mathbin{\cmdkl{\&}}}
\newrobustcmd{\knightgen}{\vcenter{\hbox{\scalebox{0.65}{ $\begin{tikzpicture}
	\begin{pgfonlayer}{nodelayer}
		\node [style=none] (5) at (0, -1) {};
		\node [style=none] (6) at (0, 0) {};
		\node [style=none] (7) at (1, 0) {};
	\end{pgfonlayer}
	\begin{pgfonlayer}{edgelayer}
		\draw [style=green] (5.center) to (6.center);
		\draw [style=greendouble] (6.center) to (7.center);
	\end{pgfonlayer}
\end{tikzpicture}
$ }}}}
\newrobustcmd{\knightgenc}{\vcenter{\hbox{\scalebox{0.65}{ $\input{knightc.tikz}$ }}}}
\newrobustcmd{\knightbisgen}{\vcenter{\hbox{\scalebox{0.65}{ $\begin{tikzpicture}
	\begin{pgfonlayer}{nodelayer}
		\node [style=none] (5) at (0, -1) {};
		\node [style=none] (6) at (0, 0) {};
		\node [style=none] (7) at (1, 0) {};
	\end{pgfonlayer}
	\begin{pgfonlayer}{edgelayer}
		\draw [style=green] (5.center) to (6.center);
		\draw [style=green] (6.center) to (7.center);
	\end{pgfonlayer}
\end{tikzpicture}
$ }}}}

\knowledgenewmathcommand{\probs}{\cmdkl{\mathbb{I}}}
\knowledgenewmathcommand{\probstar}{\cmdkl{\mathbb{I}_*}}
\newlength{\ketketspacing}
\newlength{\nketketspacing}
\newrobustcmd{\ket}[1]{\mathord{\lvert #1 \rangle}\futurelet\next\ketspacing}
\newrobustcmd{\nket}[1]{\mathord{\lvert #1 ]}\futurelet\next\nketspacing}
\def\ketspacing{
    \ifx\next\ket
        \hspace{\ketketspacing}
    \else\ifx\next\nket
        \hspace{\ketketspacing}
    \else
        \relax
    \fi\fi
}
\def\nketspacing{
    \ifx\next\ket
        \hspace{\nketketspacing}
    \else\ifx\next\nket
        \hspace{\nketketspacing}
    \else
        \relax
    \fi\fi
}

\newrobustcmd{\IH}{\textsf{IH}}
\newrobustcmd{\FS}{\textsf{FS}}

\renewrobustcmd{\th}{\mathsf{T}}
\knowledgenewmathcommand{\thI}{\cmdkl{\th_{\!\scriptscriptstyle \mathsf{Int}}}}
\knowledgenewmathcommand{\thCC}{\cmdkl{\th_{\!\scriptscriptstyle \mathsf{CC}}}}
\knowledgenewmathcommand{\thCCi}{\cmdkl{\th_{\!\scriptscriptstyle \mathsf{CC}i}}}
\knowledgenewmathcommand{\flatt}[1]{#1^{\cmdkl{\flat}}}
\newrobustcmd{\sem}[1]{\left\llbracket{#1}\right\rrbracket}
\newrobustcmd{\op}[1]{{#1}^{\mathrm{op}}}
\renewrobustcmd{\b}[1]{\mathtt{#1}}
\makeatletter
\newrobustcmd*{\garrow}[1]{
  \relax\if@display
    \expandafter\xrightarrow{\ensuremath{#1}}
  \else
    \expandafter\xrightarrow{\raisebox{-0.2ex}{\scalebox{0.6}{\ensuremath{#1}}}}
  \fi
}
\makeatother

\LoopCommands\lettersUppercase[cat#1]
{\newmathcommand#2{\mathbf#1}}
\LoopCommands\lettersLowercase[tt#1]
{\newmathcommand#2{\mathtt#1}}
\LoopCommands\lettersLowercase[vv#1]
{\newmathcommand#2{\vec#1}}

\newrobustcmd{\Dset}{\mathcal{D}}
\newrobustcmd{\Dmet}{\overline{\mathcal{D}}}
\newrobustcmd{\dist}{\varphi}
\newrobustcmd{\distb}{\psi}
\newrobustcmd{\distc}{\tau}
\newrobustcmd{\Dist}{\Phi}
\newrobustcmd{\Distb}{\Psi}
\newrobustcmd{\dirac}[1]{\delta_{#1}}
\newrobustcmd{\supp}[1]{\mathrm{supp}(#1)}
\newrobustcmd{\tv}{\mathsf{tv}}
\newrobustcmd{\tvmax}{\mathsf{tv}_{\times}}
\newrobustcmd{\tvplus}{\mathsf{tv}_{\otimes}}
\newrobustcmd{\Kant}[1]{#1_{\mathrm{K}}}
\newrobustcmd{\Cpl}{\mathfrak{C}}
\newrobustcmd{\discrete}{d_{\top}}
\newrobustcmd{\fset}[1]{\underline{#1}}
\newrobustcmd{\drel}[1]{\Delta_{#1}}
\newrobustcmd{\ndrel}[1]{\Delta^{\!\mathsf{c}}_{#1}}
\newrobustcmd{\grade}{\gamma}
\newrobustcmd{\gradeb}{\delta}
\newrobustcmd{\regrade}{\withkl{\kl[regrade]}{\mathbin{\cmdkl{\star}}}}
\newrobustcmd{\act}{\withkl{\kl[action]}{\mathbin{\cmdkl{\bullet}}}}
\newrobustcmd{\comp}{\mathbin{\fatsemi}}
\newrobustcmd{\compsimple}{\mathbin{;}}
\newrobustcmd{\tensor}{\mathbin{\otimes}}
\newrobustcmd{\swap}[1]{\sigma_{\scriptscriptstyle #1}}
\newrobustcmd{\swapg}[1]{\sigma_{\scriptscriptstyle #1}}
\newrobustcmd{\unitG}{I_{\G}}
\newrobustcmd{\unitC}{I_{\catC}}
\newrobustcmd{\emat}{[]}
\newrobustcmd{\eword}{\varepsilon}
\newrobustcmd{\word}{u}
\newrobustcmd{\wordb}{v}
\newrobustcmd{\wordc}{w}
\newrobustcmd{\words}[1]{#1^*}
\newrobustcmd{\col}{\mathtt{c}}
\newrobustcmd{\colb}{\mathtt{d}}
\newrobustcmd{\ar}{\mathsf{ar}}
\newrobustcmd{\coar}{\mathsf{coar}}
\newrobustcmd{\grd}{\mathsf{grd}}
\newrobustcmd{\tsp}[1]{#1^{\!\mathsf T}}
\newrobustcmd{\inlinefig}[2]{\vcenter{\hbox{\scalebox{#1}{ $\input{#2.tikz}$ }}}}
\newrobustcmd{\copair}[1]{\left[ #1 \right]}
\newsavebox{\closebox}
\newlength{\closeheight}
\newlength{\closewidth}
\knowledgenewmathcommand{\transformer}[1]{ 
    \savebox{\closebox}{$ #1 $} 
    \settoheight{\closeheight}{\usebox{\closebox}} 
    \settowidth{\closewidth}{\usebox{\closebox}} 
    \ooalign{ 
        \raisebox{\closeheight+0.1em} 
                {\cmdkl{\rule{\closewidth}{0.5pt}}
            }\cr
        \usebox{\closebox}
    }
}

\newrobustcmd{\letter}{f}
\newrobustcmd{\bintp}[1]{\renewrobustcmd{\letter}{#1}\raisebox{0.3ex}{\scalebox{0.75}{ $\macrotikzfig{genbintp}{\letter}$ }}}
\newrobustcmd{\flip}[1]{\renewrobustcmd{\letter}{#1}\vcenter{\hbox{\scalebox{0.75}{ $\macrotikzfig{stateletter}{\letter}$ }}}}
\newrobustcmd{\sigCA}{\Sigma_{\mathsf{CA}}}
\newrobustcmd{\img}[1]{\vcenter{\hbox{\includegraphics[scale=0.15]{#1}}}}
\newrobustcmd{\one}{I}
\newrobustcmd{\RuleRefl}{\hyperlink{rulerefl}{\textsc{Refl}}}
\newrobustcmd{\RuleSymm}{\hyperlink{rulesymm}{\textsc{Symm}}}
\newrobustcmd{\RuleTrans}{\hyperlink{ruletrans}{\textsc{Trans}}}
\newrobustcmd{\RuleSeq}{\hyperlink{ruleseq}{\textsc{Seq}}}
\newrobustcmd{\RulePar}{\hyperlink{rulepar}{\textsc{Par}}}
\newrobustcmd{\RuleReg}{\hyperlink{rulereg}{\textsc{Reg}}}
\newrobustcmd{\RuleVar}{\hyperlink{rulevar}{\textsc{Var}}}
\newrobustcmd{\RuleWeak}{\hyperlink{ruleweak}{\textsc{Weak}}}
\newrobustcmd{\RulePair}{\hyperlink{rulepair}{\textsc{Pair}}}
\newrobustcmd{\RuleFst}{\hyperlink{rulefst}{\textsc{Fst}}}
\newrobustcmd{\RuleSnd}{\hyperlink{rulesnd}{\textsc{Snd}}}
\newrobustcmd{\RuleIF}{\hyperlink{ruleif}{\textsc{If}}}
\newrobustcmd{\RuleBern}{\hyperlink{rulebern}{\textsc{Flip}}}
\newrobustcmd{\RuleKnight}{\hyperlink{ruleknight}{\textsc{Knight}}}
\newrobustcmd{\RuleLet}{\hyperlink{rulelet}{\textsc{Let}}}
\newrobustcmd{\RuleObs}{\hyperlink{ruleobs}{\textsc{Obs}}}
\newrobustcmd{\RuleEx}{\hyperlink{ruleex}{\textsc{Ex}}}

\newrobustcmd{\pair}[2]{\langle #1, #2\rangle}
\newrobustcmd{\fst}[1]{\mathtt{fst}\ #1}
\newrobustcmd{\snd}[1]{\mathtt{snd}\ #1}
\newrobustcmd{\ifte}[3]{\mathtt{if}\ #1\ \mathtt{then}\ #2\ \mathtt{else} \ #3}
\newrobustcmd{\bern}[1]{\mathtt{flip}(#1)}
\newrobustcmd{\knight}{\mathtt{knight}}
\newrobustcmd{\letin}[3]{\mathtt{let}\ #1 = #2\ \mathtt{in}\ #3}
\newrobustcmd{\observe}[1]{\mathtt{observe}\ #1}
\makeatletter
\newrobustcmd{\bigplus}{
  \DOTSB\mathop{\mathpalette\mattos@bigplus\relax}\slimits@
}
\newcommand\mattos@bigplus[2]{
  \vcenter{\hbox{
    \sbox\z@{$#1\sum$}
    \resizebox{!}{0.9\dimexpr\ht\z@+\dp\z@}{\raisebox{\depth}{$\m@th#1+$}}
  }}
  \vphantom{\sum}
}
\makeatother

\newrobustcmd{\syncat}[2]{\mathcal{S}_{\scriptstyle #1,#2}}
\knowledgenewmathcommand{\sync}[1]{\cmdkl{\mathcal{S}}_{\scriptstyle #1}}
\knowledgenewmathcommand{\vync}[1]{\cmdkl{\mathcal{M}}_{\scriptstyle #1}}
\knowledgenewmathcommand{\synct}[1]{\cmdkl{\mathcal{S}}_{\scriptstyle #1}}
\knowledgenewmathcommand{\vynct}[1]{\cmdkl{\mathcal{M}}_{\scriptstyle #1}}
\newrobustcmd{\diagbox}[3]{
\begin{tikzpicture}
	\begin{pgfonlayer}{nodelayer}
		\node [style=basic box] (0) at (0, 0) {$#1$};
		\node [style=none] (1) at (1.5, 0) {};
		\node [style=none] (2) at (-1.5, 0) {};
		\node [style=none] (3) at (1.5, 0.5) {\scriptsize $#3$};
		\node [style=none] (4) at (-1.5, 0.5) {\scriptsize $#2$};
	\end{pgfonlayer}
	\begin{pgfonlayer}{edgelayer}
		\draw (2.center) to (0);
		\draw (0) to (1.center);
	\end{pgfonlayer}
\end{tikzpicture}
}

\newrobustcmd{\normd}{\mathcal{N}}
\newrobustcmd{\idmat}[1]{I_{#1}}
\newrobustcmd{\Aff}{\mathbf{Aff}}
\newrobustcmd{\CoIso}{\mathbf{CoIso}}
\newrobustcmd{\Gauss}{\mathbf{Gauss}}
\newrobustcmd{\GaussSharp}{\mathbf{Gauss}^{\sharp}}
\knowledgenewmathcommand{\Coll}{\cmdkl{\mathrm{Coll}}}
\newrobustcmd{\Cst}{\mathrm{Const}}
\newrobustcmd{\thAff}{\th_{\!\scriptscriptstyle\mathsf{Aff}}}
\newrobustcmd{\thCoIso}{\th_{\!\scriptscriptstyle\mathsf{CoIso}}}
\newrobustcmd{\thGauss}{\th_{\!\scriptscriptstyle\mathsf{Gauss}}}
\newrobustcmd{\thGaussSharp}{\th_{\!\scriptscriptstyle\mathsf{Gauss}^{\sharp}}}
\newrobustcmd{\noisegen}{\knightgen}
\newrobustcmd{\gaussgen}{\inlinefig{0.7}{gaussgen}}
\newrobustcmd{\affdiag}[1]{\renewrobustcmd{\letter}{#1}\vcenter{\hbox{{\scalebox{0.75}{ $\macrotikzfig{affdiag}{\letter}$ }}}}}
\newrobustcmd{\mean}{\mu}
\newrobustcmd{\meanb}{\nu}
\newrobustcmd{\cov}{\Sigma}
\newrobustcmd{\covb}{\Xi}
\newrobustcmd{\zmat}{\mathbf{0}}
\knowledgenewmathcommand{\collclass}[1]{\cmdkl{[}#1\cmdkl{]}}
\knowledgenewmathcommand{\collcong}{\mathrel{\cmdkl{\simeq}}}
\knowledgenewmathcommand{\utg}{\withkl{\kl[up to regrading]}{\mathrel{\cmdkl{\approx}}}}
\newrobustcmd{\colleq}{\equiv}
\newrobustcmd{\castgrade}[2]{{#2}^{#1}}

\newcounter{theqn}

\title[Graded String Diagrams]{Graded String Diagrams\\ with Applications in Probabilistic Reasoning}
\author[R.~Sarkis]{Ralph Sarkis \lmcsorcid{0000-0002-9037-2435}}
\address{Université de Lille, UMR 9189 CRIStAL, Lille, France}
\author[F.~Zanasi]{Fabio Zanasi \lmcsorcid{0000-0001-6457-1345}}
\address{University College London, United Kingdom}
\date{}

\thanks{This research was supported by the ARIA Safeguarded AI programme.}

\begin{document}
\begin{abstract}
  We introduce string diagrams for "graded symmetric monoidal categories". Our approach includes a definition of "graded monoidal theory" and the corresponding freely generated "syntactic category@@GRD". The "para construction" on monoidal "actegories" is a motivating example for our framework, and we show how an axiomatic presentation of the resulting category $\Para{\G}{\catC}$ may be modularly obtained from presentations for $\G$ and $\catC$. We demonstrate three different applications of our diagrammatic language under the same theme: probabilistic reasoning.

  First, we show how to axiomatise a variant of the "graded category" $\ImP$, recently introduced by Liell-Cock and Staton to model imprecise probability~\cite{LiellCock2025}. This culminates in a representation, as string diagrams with grading wires, of programs with primitives for nondeterministic and probabilistic choices and conditioning. Second, we exploit our modular presentation result to provide a new proof for a presentation of $\Gauss$, the category of Gaussian stochastic maps. Third, we use grading wires to represent causal interventions on Bayesian networks.
\end{abstract}

\maketitle

\section{Introduction}\label{sec:intro}

Graded categorical structures have been successfully applied in the study of effect systems~\cite{Mellies2012,Katsumata2014,Gaboardi2016,Mycroft2016,Sanada2023,Kura2026}, coalgebraic semantics~\cite{Milius2015,Dorsch2019,Ford2021c,Ford2022,Forster2024,DiLavore2026}, cellular automata~\cite{Capobianco2023}, program logics~\cite{Gaboardi2021}, and most recently to model imprecise probabilistic reasoning~\cite{LiellCock2025}.

In a "graded category" $\catC$, every morphism $f$ has a grade $\grade$ that is an object in a monoidal category $\G$. This grade can be understood as a type of parameters that $f$ uses, separate from its inputs and outputs. Grades have been used to represent, among other things, side effects of a program~\cite{Katsumata2014,Gaboardi2021}, a bound on the length of relevant traces~\cite{Milius2015}, the neighbourhood of a cell in a cellular automaton~\cite{Capobianco2023}, or the number of nondeterministic choices consumed by a process~\cite{LiellCock2025}. When composing graded morphisms, in sequence or in parallel, the grades combine according to the monoidal product in $\G$. The grade of a morphism can also be changed by "regrading" with a morphism in $\G$, which informally consists in acting on parameters before the graded morphism has access to them.

Somewhat independently of these developments, there is an increasingly popular research area using monoidal categories to model resource-sensitive computational processes, such as those found in probabilistic programming~\cite{ChoJacobs2019,Fritz2020}, quantum theory~\cite{Coecke2017}, concurrency~\cite{Bonchi2019b}, and machine learning~\cite{Cruttwell2022}. This family of approaches adopts \emph{string diagrams}, the graphical language of monoidal categories~\cite{Selinger_2010,PiedeleuZanasi2025}, as a visual and yet rigorous syntax for processes, providing an intuitive understanding of how information flows and is exchanged between components within a system. Particular emphasis is put on finding an \emph{axiomatic presentation} of monoidal categories of processes, in the form of a diagrammatic calculus whose equational theory is complete with respect to the semantic equivalence of interest. This perspective has been fruitful for instance in quantum theory~\cite{Hadzihasanovic2015,Jeandel2018,Jeandel2020,Backens2023,Poor2023}, concurrency~\cite{Bonchi2019b}, linear algebra~\cite{zanasi:tel-01218015,Baez2018,Boisseau2022}, automata theory~\cite{Piedeleu2023,Piedeleu2025}, and probabilistic programming~\cite{Piedeleu2025b,Stein2026}. It has determined the development of \emph{monoidal algebra}~\cite{Bonchi2018}, as the counterpart of traditional categorical algebra~\cite{Lawvere1963} adapted to a monoidal rather than cartesian setting.

Given how much graded category theory and monoidal algebra overlap in terms of applications, it is natural to ask how they can be combined in a uniform approach. Even though many extensions of monoidal algebra have been considered~\cite{DiLiberti2021,Bonchi2023,Bonchi2024b,Villoria2024,Bonchi2025b,Bonchi2025,Bonchi2025PolyJ,Lobbia2025}, none of these approaches readily incorporates a notion of grade. In essence, string diagrams always represent a morphism as a box with input and output wires $\inlinefig{0.7}{monoidal/montermwires}$, whereas morphisms in a "graded category" also interact through a grade, which cannot be conflated with the domain or codomain. The syntaxes for bicategories (see e.g.\@ \cite{Ponto2013}) and for double categories \cite{Myers2018} \textit{do} provide a second axis of sequential composition orthogonal to the domain/codomain one. However, this makes parallel composition---which is central to monoidal algebra---hard to represent, as it would require a third dimension graphically. In other words, graded categories could be handled within the existing 2-dimensional calculi, but our interest here is in structures that are graded and monoidal. We are not aware of a general diagrammatic syntax for monoidal double categories or monoidal bicategories, but some work in this direction exists \cite{Nester2023}.

\paragraph{Our Contributions}
To fill this gap, we propose a string diagrammatic language for "graded symmetric strict monoidal categories" (SMC). A morphism is pictured as a box with input wires, output wires, and grading wires: $\inlinefig{0.7}{fwires}$ or $\inlinefig{0.7}{f}$ after bunching up multiple wires into one. Morphisms can be composed in sequence and in parallel as usual, but they can also be "regraded". We represent this last operation by plugging a diagram representing the grading morphism to the grading wire as shown below on the left. Translating the laws of "graded SMCs" to this graphical syntax (\cref{fig:axiomssyncat}) shows that any reasonable deformation of a graded string diagram does not change the represented morphism, e.g.\@ the two diagrams below on the right are equal.
\[\scalebox{0.8}{\input{deformations.tikz}}\]

Given a signature where operation symbols have an arity, a coarity, and a grade, we show how to freely generate a "graded SMC" of string diagrams that obey these visually intuitive equations. We also show how to impose arbitrary additional equations, in the form of a \emph{"graded monoidal theory"}: imposing such equations is necessary to accurately model concrete settings, whenever one wants to reason diagrammatically about the modelled processes. This development closely follows the analogue for (classical) "monoidal theories", and in fact we prove (\cref{lem:zerogradedpart}) that diagrams with no grading wire (i.e.\@ graded by the empty word $\eword$) can actually be considered as classical string diagrams over the signature restricted to operations with no grading wire. We define a second construction taking graded string diagrams to classical ones by erasing grading wires (\cref{defn:flatth}). We show that it corresponds to the operation of identifying graded morphisms "up to regrading" (\cref{thm:collapsepresentation}).

Our syntax draws inspiration from similar pictorial languages for categories obtained with the "para construction"~\cite{CapucciGHR2022,Cruttwell2022,Gavranovic2024}, which however  were not formalised. In fact, our framework directly applies to the "para construction", and provides the following `modularity' principle (\cref{thm:gradedpresfrompres}): given axiomatisations of "SMCs" $\G$ and $\catC$, and a suitable action of $\G$ on $\catC$, we can derive an axiomatisation of the "graded category" $\Para{\G}{\catC}$ obtained via the "para construction", provided the monoidal unit in $\G$ is terminal with no outgoing morphism except the identity. We apply this principle twice in this paper.

Its first use underlies the axiomatisation of a graded string diagrammatic language for imprecise probability. Semantically, our starting point is the approach of Liell-Cock and Staton~\cite{LiellCock2025}, who model imprecise probability with a "graded category" whose base category consists of (precise) probabilistic processes and whose grades keep track of nondeterministic choices: their category $\ImP$ is obtained via the "para construction", from the action on $\FStoch$ of a subcategory of stochastic matrices representing choices. We follow the same recipe, but instantiate it differently: as base category we take the subcategory $\BStoch$ of $\FStoch$ that was recently axiomatised in~\cite{Piedeleu2025b}, and as grading category we take $\op{\FInj}$, the opposite of the category of finite sets and injections, which~\cite{LiellCock2025} observe to be sufficient for their purposes. The resulting "graded category" $\BImP$ is a binary variant of $\ImP$, to which \cref{thm:gradedpresfrompres} applies. Given $n,m,a \in \N$, a morphism $n \to m$ of grade $a$ in $\BImP$ is a list of $2^a$ stochastic matrices of dimension $2^m \times 2^n$, which we conveniently arrange as a single stochastic matrix divided into $2^a$ submatrices. Intuitively, such a morphism views each submatrix as a stochastic process and nondeterministically chooses one of them to execute. "Regrading" in $\BImP$ amounts to reordering and copying submatrices (not all such operations are possible).

We then turn to the string diagrammatic presentation of $\BImP$. The main result in~\cite{Piedeleu2025b} states that the category of stochastic matrices whose dimensions are powers of $2$ is presented by a theory they call $\Circ$. Following the modular recipe given in \cref{thm:gradedpresfrompres}, we add a generator $\inlinefig{0.75}{knightbis}$ of grade $1$ to $\Circ$ (its grading wire is the one pointing downwards) that we interpret as outputting a value which nondeterministically evaluates to true or false, and an equation saying that to delete this value is the same thing as not outputting it. We obtain the "graded theory" $\GCirc$ (imprecise $\Circ$), which axiomatises the "graded category" $\BImP$.

We can thus extend the programming language in~\cite{Piedeleu2025b}, which is interpreted in $\Circ$, with a construct for nondeterministic choice, denoted with $\knight$, and interpret it in $\GCirc$. We show (\cref{thm:desiderata}) that this language is commutative and affine, meaning that expressions that are independent of each other (e.g.\@ they concern different variables) can be rearranged. Furthermore, we demonstrate the usefulness of our modular approach by showing how to easily add the conditioning generator (and its corresponding programming construct) studied in~\cite{Piedeleu2025b}. We illustrate the features of this graphical programming language with a probabilistic riddle that is formalised with string diagrams (\cref{exmp:boyorgirl}).

The modularity principle is used a second time to give a new derivation of the string-diagrammatic presentation of Gaussian probability from~\cite{Stein2026}. We define a "graded SMC" $\GaussSharp$ as the "para construction" in which affine maps (the base category) with noisy parameters (coming from the grades) represent Gaussian maps. Two graded morphisms in $\GaussSharp$ may represent the same Gaussian map. We choose the grading category $\CoIso$, the category of matrices $A$ satisfying $A\tsp{A} = \id$, so that identifying morphisms of $\GaussSharp$ "up to regrading" removes this ambiguity, yielding the category $\Gauss$ of Gaussian maps. We then apply \cref{thm:gradedpresfrompres} to obtain a graded presentation of $\GaussSharp$ from presentations of $\CoIso$ and $\Aff$, and use \cref{thm:collapsepresentation} to obtain a presentation of $\Gauss$ by erasing the grading wires. The resulting theory coincides with the presentation of~\cite{Stein2026} (except for two redundant axioms). Our approach breaks the completeness proof of~\cite{Stein2026} into conceptually distinct steps, allowing the use of abstract diagrammatic machinery.

Our last contribution is a diagrammatic language for interventions on causal models. We upgrade the established formalisation of causal models as string diagrams (see e.g.\@ \cite{Coecke2012,Fong2013,Jacobs2021}) with grading to represent causal interventions. While previous works~\cite{JacobsKZ21,Yimu2021,Lorenz2023} view an intervention as a functor that transforms diagrams, our approach realises interventions as part of the diagrammatic syntax. We also show that the result of~\cite[Theorem~8.1]{JacobsKZ21} on computing the effect of an intervention on a single variable from observational data transfers to this setting (\cref{thm:gradedidentifiability}): here, the observational and the interventional models are "regradings" of a single graded diagram, and the original proof applies with one additional step accounting for the grading. This is similar in spirit to the preliminary report of Stein~\cite{Stein2025} on designing a programming language for causal models where interventions are first-class.

\textbf{Synopsis.} In \cref{sec:background}, we recall the necessary  background on "monoidal theories" and "actegories". \cref{sec:gradedtheories} contains the main definitions of this work. We introduce here "graded SMCs", "graded monoidal theories", and their string diagrammatic syntax. We also describe the construction of a syntactic "graded prop" generated by a "graded monoidal theory", thereby providing a sound and complete list of equational laws to reason with terms of such theories. In \cref{sec:collapse}, we describe the construction to get an "SMC" from a "graded SMC" by considering graded morphisms "up to regrading", and we provide a syntactic counterpart.
In \cref{sec:modularcompleteness}, we show how to combine two "monoidal theories" into a "graded monoidal theory" through the graded "para construction". In \cref{sec:application}, we define the "graded category" $\BImP$ and we apply the result of \cref{sec:modularcompleteness} to obtain a presentation. In \cref{sec:impprog}, we obtain a graphical interpretation of a programming language with nondeterministic and probabilistic constructs, illustrating our approach with two concrete examples. In \cref{sec:gauss}, we rederive a presentation of $\Gauss$ using the results of \cref{sec:collapse,sec:modularcompleteness}. Finally, in \cref{sec:interventions}, we define a theory of graded string diagrams which are interpreted as stochastic processes with annotated points of intervention.

\textbf{Comparison with conference version.} This paper is an extended version of our work published in~\cite{Sarkis2025}. It includes the missing proofs as well as three entirely new sections: \cref{sec:collapse,sec:gauss,sec:interventions}. We also generalised our main results from single-sorted "props" to coloured (or multisorted) ones, adjusting the preliminaries and definitions accordingly. Another addition is \cref{appendix:computing}, where we describe simple methods to compute with morphisms inside $\BImP$.

\textbf{Acknowledgements.} We thank the anonymous reviewers for their careful reading of previous versions and their very helpful comments.

\section{Background}\label{sec:background}

\subsection{Monoidal Theories}
We recall background on string diagrams and "monoidal theories", following the terminology of~\cite{Bonchi2018}. \AP First, a \textbf{""coloured prop""} (sometimes called multisorted "prop") is a symmetric strict monoidal category (""SMC"")\footnote{Throughout, we assume that "SMCs" are \textbf{strict}, namely, the associators and unitors are identities.} whose objects are finite words over a fixed set of colours and whose monoidal product on objects is concatenation. "Props" are the basic framework to study algebraic objects with a monoidal structure. They may be freely obtained from generators and equations. This requires introducing a notion of symmetric monoidal signature and symmetric "monoidal theory". In the sequel, we will omit the adjectives `symmetric' and `coloured' for brevity, but they are assumed.

\begin{defi}[Monoidal term]\label{defn:monterms}
  \AP A \textbf{""monoidal signature""} $\sig$ is a set $\sig_0$ of \textbf{colours} and a set $\sig_1$ of \textbf{generators}, together with functions $\ar,\coar: \sig_1 \rightarrow \words{\sig_0}$ assigning to each generator an \textbf{arity} $\word$ and \textbf{coarity} $\wordb$ in $\words{\sig_0}$, the set of finite words over $\sig_0$. We use the notation $g: \word \rightarrow \wordb \in \sig$ to compactly represent the data of a generator in the signature. We inductively define the set $\intro*\vynct{\sig}$ of \textbf{"monoidal $\sig$-terms@monoidal term"} (and their (co)arities) built with the generators in $\sig$ via the rules \cref{eqn:monterms:trivial,eqn:monterms:SMC,eqn:monterms:seqcomp,eqn:monterms:parcomp}, where $\eword$ denotes the empty word in $\sig_0$.\\
  \begin{minipage}{0.27\textwidth}
    \begin{gather}
      \label{eqn:monterms:trivial}
      \scalebox{0.85}{\begin{bprooftree}
          \AxiomC{$g: \word \rightarrow  \wordb \in \sig$}
          \UnaryInfC{$g: \word \rightarrow \wordb \in \vynct{\sig}$}
        \end{bprooftree}}
    \end{gather}
  \end{minipage}
  \begin{minipage}{0.72\textwidth}
    \begin{gather}
      \label{eqn:monterms:SMC}
      \scalebox{0.85}{\begin{bprooftree}
          \AxiomC{$\phantom{\sig}$}
          \UnaryInfC{$\id_{\eword}: \eword \rightarrow \eword \in \vynct{\sig}$}
        \end{bprooftree}}
      \scalebox{0.85}{\begin{bprooftree}
          \AxiomC{$\col \in \sig_0$}
          \UnaryInfC{$\id_{\col}: \col \rightarrow \col \in \vynct{\sig}$}
        \end{bprooftree}}
      \scalebox{0.85}{\begin{bprooftree}
          \AxiomC{$\col, \colb \in \sig_0$}
          \UnaryInfC{$\swap{\col,\colb}: \col\colb \rightarrow \colb\col \in \vynct{\sig}$}
        \end{bprooftree}}
    \end{gather}
  \end{minipage}\\[0.3em]
  \begin{minipage}{0.48\textwidth}
    \begin{gather}
      \label{eqn:monterms:seqcomp}
      \scalebox{0.85}{\begin{bprooftree}
          \AxiomC{$f: \word \rightarrow \wordb \in \vynct{\sig}$}
          \AxiomC{$g: \wordb \rightarrow \wordc \in \vynct{\sig}$}
          \BinaryInfC{$f\compsimple  g : \word \rightarrow \wordc \in \vynct{\sig}$}
        \end{bprooftree}}
    \end{gather}
  \end{minipage}
  \begin{minipage}{0.51\textwidth}
    \begin{gather}
      \label{eqn:monterms:parcomp}
      \scalebox{0.85}{\begin{bprooftree}
          \AxiomC{$f: \word \rightarrow \wordb \in \vynct{\sig}$}
          \AxiomC{$f': \word' \rightarrow \wordb' \in \vynct{\sig}$}
          \BinaryInfC{$f\tensor f' : \word\word' \rightarrow \wordb\wordb' \in \vynct{\sig}$}
        \end{bprooftree}}
    \end{gather}
  \end{minipage}\\[0.6em]
  These formation rules of "monoidal terms" correspond straightforwardly with the structure of an "SMC". We can see any generator as a term via \cref{eqn:monterms:trivial}; we can compose "monoidal terms" in sequence and in parallel (monoidal product) with \cref{eqn:monterms:seqcomp} and \cref{eqn:monterms:parcomp} respectively; and there are identity "monoidal terms" and symmetries for all words: \cref{eqn:monterms:SMC} provides them for colours and for the empty word, and those for longer words are obtained as monoidal products \cref{eqn:monterms:parcomp} of these basic terms (combined with sequential compositions \cref{eqn:monterms:seqcomp} in the case of symmetries).

  As customary, we use \emph{string diagrams} (see e.g.\@ \cite{Selinger_2010,PiedeleuZanasi2025}) to visually represent monoidal terms. A generic term $f: \word \rightarrow \wordb \in \vynct{\sig}$ is depicted as a box with a wire on the left labelled with $\word$ and a wire on the right labelled with $\wordb$: $\inlinefig{0.6}{monoidal/monterm}$. We write the terms $\id_{\eword}$ as $\inlinefig{0.8}{idzer}$, $\id_{\col}$ as $\inlinefig{0.8}{wirecol}$, and $\swap{\col,\colb}$ as $\inlinefig{0.6}{swapcol}$. We will sometimes omit the labels of wires when they are irrelevant or clear from the context. The two forms of composition are respectively \[\inlinefig{0.69}{monoidal/fnm} \compsimple \inlinefig{0.69}{monoidal/gml} = \inlinefig{0.69}{monoidal/seqcomp} \text{ and }  \inlinefig{0.69}{monoidal/fnm} \tensor \inlinefig{0.69}{monoidal/fnmprime} = \inlinefig{0.69}{monoidal/parcomp}\]
\end{defi}
\begin{defi}[Theory]\label{defn:month}
  \AP A \textbf{""monoidal theory""} is a tuple $\th = (\sig,E)$, where $\sig$ is a "monoidal signature", and $E$ is a set of equations, called \textbf{axioms}, between "monoidal $\sig$-terms@monoidal term" whose arities and coarities coincide, i.e.\@ pairs $(s,t) \in \vynct{\sig}\times \vynct{\sig}$ with $\ar(s) = \ar(t)$ and $\coar(s) = \coar(t)$. We denote such equations by $s=t$.
\end{defi}
\begin{defi}[Model]\label{defn:monthmodel}
  Let $\th = (\sig,E)$ be a "monoidal theory" and $\catC$ an "SMC". Given an assignment $M: \sig \rightarrow \catC$ sending each colour $\col \in \sig_0$ to an object $M\col \in \catC$, and each generator $g: \word \rightarrow \wordb \in \sig$ to a morphism in $\catC(\bigotimes_{\col\in \word}M\col, \bigotimes_{\col\in\wordb}M\col)$, we can canonically extend $M$ to all terms in $\vynct{\sig}$, because each term formation rule \cref{eqn:monterms:SMC,eqn:monterms:seqcomp,eqn:monterms:parcomp} corresponds to part of the structure in an "SMC". A \textbf{model} of $\th$ valued in $\catC$ is such an assignment $M$ satisfying $M(s) = M(t)$ for all axioms $s=t \in E$.
\end{defi}

Models of "monoidal theories" are usually viewed through the lens of functorial semantics, \emph{\`{a} la} Lawvere~\cite{Lawvere1963}. In a nutshell, we can construct a "syntactic category@@MON" $\vync{\th}$ for a "theory@@MON" $\th$, such that models of $\th$ correspond to symmetric strict monoidal functors with domain $\vync{\th}$. The morphisms of $\vync{\th}$ are monoidal terms quotiented by some congruence generated by the axioms in $\th$ and the necessary equations to satisfy the laws of an "SMC".
\begin{defi}[Closure]\label{defn:monthclosure}
  Given a "theory@@MON" $\th = (\sig,E)$, the \textbf{closure} of $E$, denoted $\closed{E}$, is the smallest congruence on $\vynct{\sig}$ containing $E$ and all the laws of "SMCs", see~\cref{fig:axiomsvyncat}.
\end{defi}
\begin{defi}[Syntactic category]\label{defn:monthsyncat}
  \AP The \textbf{""syntactic category@@MON""} of a "monoidal theory" $\th = (\sig,E)$ is a "prop" $\intro*\vync{\th}$ where the set of morphisms $\vync{\th}(\word,\wordb)$ consists of "monoidal $\sig$-terms@monoidal term" of arity $\word$ and coarity $\wordb$, quotiented by the congruence $\closed{E}$.
  Identities, symmetries, sequential composition, and monoidal product are defined via formation rules \cref{eqn:monterms:SMC,eqn:monterms:seqcomp,eqn:monterms:parcomp}. The axioms of \cref{fig:axiomsvyncat} ensure that all the laws of "SMCs" hold.

  \setcounter{theqn}{0}
  \begin{figure}
    \begin{minipage}{0.5\textwidth}
      \begin{gather}\label{eqn:month:assocseq}
        \inlinefig{0.55}{monoidal/assocseqcomp} \stepcounter{theqn}\tag{M\arabic{theqn}}\\[0.2em]
        \label{eqn:month:neutralid}
        \inlinefig{0.6}{monoidal/neutralid} \stepcounter{theqn}\tag{M\arabic{theqn}}\\[0.3em]
        \label{eqn:month:neutralidzer}
        \inlinefig{0.6}{monoidal/neutralidzer} \stepcounter{theqn}\tag{M\arabic{theqn}}\\[0.4em]\label{eqn:month:assocpar}
        \inlinefig{0.6}{monoidal/assocparcomp} \stepcounter{theqn}\tag{M\arabic{theqn}}\\[0.4em]\label{eqn:month:overthewire}
        \inlinefig{0.6}{monoidal/overthewire} \stepcounter{theqn}\tag{M\arabic{theqn}}
      \end{gather}
    \end{minipage}
    \begin{minipage}{0.49\textwidth}
      \begin{gather}\label{eqn:month:swapswap}
        \inlinefig{0.6}{monoidal/swapswap} \stepcounter{theqn}\tag{M\arabic{theqn}}\\[0.4em]\label{eqn:month:interchange}
        \inlinefig{0.6}{monoidal/interchange}  \stepcounter{theqn}\tag{M\arabic{theqn}}\\[0.4em]\label{eqn:month:coherentid}
        \inlinefig{0.6}{coherentid} \stepcounter{theqn}\tag{M\arabic{theqn}}\\[0.4em]\label{eqn:month:coherentswaps}
        \inlinefig{0.6}{monoidal/coherentswaps} \stepcounter{theqn}\tag{M\arabic{theqn}}
      \end{gather}
    \end{minipage}
    \caption{Laws of symmetric strict monoidal categories.}\label{fig:axiomsvyncat}
  \end{figure}

  We say that $\th$ is a \textbf{presentation} of a "prop" $\catC$ (or that $\th$ \textbf{presents} $\catC$) if $\vync{\th}$ and $\catC$ are isomorphic as "props", or equivalently, if there is a bijective-on-objects, fully faithful, symmetric monoidal functor between them. Note that a symmetric strict monoidal functor between "props" is bijective on objects precisely when it is the identity on objects modulo a renaming of colours.
\end{defi}
Functorial semantics is founded upon the fact that models of $\th$ valued in $\catC$ are in 1-to-1 correspondence with symmetric monoidal functors $\vync{\th} \rightarrow \catC$.

\subsection{Actegories}
The "para construction", which plays a role in our developments, is based on "actegories". "Actegories", an abstraction of the concept of monoid/group action, are categories being acted on by monoidal categories. We recall the definition from~\cite{Capucci22,CapucciGHR2022}.
\begin{defi}[Actegory]\label{defn:act}
  \AP Let $\G$ be an "SMC" and $\catC$ be a category. An \textbf{""action""} of $\G$ on $\catC$ is a functor $\act: \G \times \catC \rightarrow \catC$ equipped with two natural families of isomorphisms, $\eta_X : I \act X \cong X$ and $\mu_{A,B,X} : A \act (B \act X) \cong (A \otimes B) \act X$,
  that satisfy the coherence laws listed in~\cite[Definition 3.1.1]{Capucci22}. We say that $\catC$ equipped with this "action" is a \textbf{$\G$-"actegory"}.

  We call $\act$ a \textbf{"symmetric monoidal action@action"} and $\catC$ a \textbf{symmetric monoidal $\G$-"actegory"} if $\catC$ is an "SMC", $\act$ is a symmetric monoidal functor (with the usual coordinatewise "SMC" structure on $\G \times \catC$), and there is a natural family of isomorphisms,
  $\kappa_{A,X,Y} : A \act (X \otimes Y) \cong X \otimes (A \act Y)$,
  that satisfy the additional coherence laws listed in \cref{fig:coherencesymmact} (\cref{appendix:background}).

  In a symmetric monoidal "actegory", the \textbf{mixed associator} ($\kappa'$) and \textbf{interchanger} ($\iota$) are defined as follows (we omit the subscripts as they can be inferred from the types):
  \[
    \kappa' = \begin{tikzcd}[scale cd=0.8, cramped, ampersand replacement=\&]
      {A \act (X \otimes Y)} \& {A \act (Y \otimes X)} \& {Y \otimes (A \act X) } \& {(A \act X) \otimes Y}
      \arrow["{A \act \sigma}", from=1-1, to=1-2]
      \arrow["\kappa", from=1-2, to=1-3]
      \arrow["\sigma", from=1-3, to=1-4]
    \end{tikzcd}                                                                \]
  \[\iota   =   \begin{tikzcd}[scale cd=0.72, sep=small, cramped, ampersand replacement=\&]
      {(A\otimes B) \act (X \otimes Y)} \& {X \otimes ((A\otimes B) \act Y)} \& {X\otimes (A \act (B \act Y))} \& {A \act (X \otimes (B \act Y))} \& {(A\act X) \otimes (B \act Y)}
      \arrow["\kappa", from=1-1, to=1-2]
      \arrow["{X \otimes \mu^{-1}}", from=1-2, to=1-3]
      \arrow["{\kappa^{-1}}", from=1-3, to=1-4]
      \arrow["{\kappa'}", from=1-4, to=1-5]
    \end{tikzcd}.\]
\end{defi}

\begin{exa}[Self-action]\label{exmp:selfaction}
  Any "SMC" $\catC$ acts on itself in a canonical way: $\act: \catC \times \catC \rightarrow \catC$ is defined as the monoidal product $\otimes$~\cite[Example~3.2.4]{Capucci22}. This is a "symmetric monoidal action" where the natural isomorphisms $\eta$, $\mu$, $\kappa$ are compositions of structure maps (associators, unitors, and symmetries). Any subcategory of $\catC$ with the inherited "SMC" structure also acts on $\catC$ in the same way.
\end{exa}

\section{Graded Monoidal Theories and Graded Props}\label{sec:gradedtheories}

\subsection{Graded categories}
\AP Let $\intro*\G$ be a symmetric strict monoidal category. 
\knowledge{\G}{autoref,scope=subsection,also now}
We call the objects of $\G$ \textbf{grades} and denote them with $\grade$, $\gradeb$, etc. We denote composition in $\G$ with $\compsimple$, monoidal product with juxtaposition, the unit with $\unitG$, and the symmetries with $\swapg{\grade\grade'}: \grade\grade' \rightarrow \grade'\grade$.
Following $\G$-graded categories~\cite{Wood1976,LucyshynWright25}, we introduce the notion of "SMC@graded SMC" graded over $\G$. In a nutshell, every morphism of a $\G$-"graded SMC" comes with a grade inside $\G$, and whenever one composes morphisms (in sequence or in parallel), grades combine via the monoidal product in $\G$.
\begin{defi}[Graded SMC]\label{defn:gradedsmc}
  \AP A \textbf{$\G$-""graded symmetric strict monoidal category""} $\catC$ (also called $\G$-"graded SMC" or $\G$"-SMC") consists of the following data.
  \begin{itemize}
    \item A class of \textbf{objects} $\mathrm{Ob}(\catC)$. We write $X \in \catC$ instead of $X \in \mathrm{Ob}(\catC)$.
    \item For all objects $X,Y \in \catC$ and grade $\grade \in \G$, a set $\catC_{\grade}(X,Y)$ of \textbf{morphisms of grade $\grade$}. We often write $f: X \garrow{\grade} Y$ to mean $f \in \catC_{\grade}(X,Y)$.
          \itemAP For all objects $X,Y \in \catC$ and morphism between grades $t: \grade \rightarrow \gradeb$, a function $t \regrade - : \catC_{\gradeb}(X,Y) \rightarrow \catC_{\grade}(X,Y)$ called \textbf{""regrading"" by $t$}.
    \item For all $X,Y,Z \in \catC$ and $\grade,\gradeb \in \G$, a \textbf{composition} ${\comp} : \catC_{\grade}(X,Y) \times \catC_{\gradeb}(Y,Z) \rightarrow \catC_{\grade\gradeb}(X,Z)$.
    \item For all objects $X \in \catC$, an \textbf{identity} morphism $\id_X \in \catC_{\unitG}(X,X)$.
    \item A \textbf{monoidal product} operation $\tensor : \mathrm{Ob}(\catC) \times \mathrm{Ob}(\catC) \rightarrow \mathrm{Ob}(\catC)$ and a \textbf{unit} $\unitC \in \mathrm{Ob}(\catC)$.    \item For all $X, X', Y, Y' \in \catC$ and $\grade, \gradeb \in \G$, a monoidal product of morphisms defined as a function $\tensor : \catC_{\grade}(X,Y) \times \catC_{\gradeb}(X',Y') \rightarrow \catC_{\grade\gradeb}(X\tensor X', Y \tensor Y')$.
    \item For all $X,X' \in \catC$, a morphism $\swap{X,X'}: X \tensor X' \garrow{\unitG} X' \tensor X$ of grade $\unitG$ called the \textbf{symmetry}.
    \item Composition, monoidal product, and symmetries must satisfy properties akin to those of "SMCs", and be compatible with "regrading". The full list of conditions is in \cref{appendix:gradedtheories}, but our diagrammatic syntax will help to visualise them (see \cref{fig:axiomssyncat}).
  \end{itemize}
\end{defi}
\begin{rem}\label{rem:interchange}
  The interchange law \cref{eqn:month:interchange} cannot hold as usual because the monoidal product in $\G$ is not commutative \emph{on the nose}: the grades of $(f \tensor f') \comp (g \tensor g')$ and $(f \comp g) \tensor (f' \comp g')$ are different objects of $\G$ (related by a symmetry of $\G$, but not equal), so the two morphisms cannot be equal. The law is therefore required to hold only "up to regrading" by that symmetry, see \cref{eqn:th:interchangeupto}. The same reasoning applies to the sliding equation: \cref{eqn:month:overthewire} becomes \cref{eqn:th:overthewire}.
\end{rem}

In this work, we say that two morphisms $f,g$ in a "graded category" are \textbf{""equal up to regrading""} if there exists a morphism between grades $t$ such that $t \regrade f = g$ or $t \regrade g = f$. We write $f \utg g$.

\begin{rem}[Terminology]
  The adjective `graded' is overloaded in mathematics, but it usually means that each element of some structure is equipped with an element of another structure (usually a monoid). Graded categories were also called `locally graded' by Levy~\cite{Levy2019}. A folklore result (see e.g.\@ \cite{Wood1976,Gavranovic2024,LucyshynWright25}) establishes a correspondence between $\G$-graded categories and $[\op{\G},\Set]$-enriched categories, where the monoidal structure $[\op{\G},\Set]$ is Day convolution~\cite{Day1970}. This extends to the monoidal case and beyond as is explained in~\cite[Proposition~3.5]{LiellCock2025}.

  Other graded monoidal categorical structures have been studied in the literature, and they differ in how grades combine under composition and monoidal product. In the group-graded monoidal categories of~\cite[Section~3]{Frohlich1974}, grades combine under composition, while the monoidal product is only defined for morphisms of the same grade and keeps that grade. In the PCM-graded monoidal categories of~\cite[Definition~3.1]{Earnshaw2026}, the monoidal product combines grades and sequential composition is restricted to morphisms of the same grade. In the duoid-graded monoidal categories of~\cite[Definition~3]{DiLavore2025} both composition and monoidal product combine grades, but with the two different monoid structures of the duoid in which grades lie. A more detailed investigation of the differences is out of the scope of this work as we will focus on a diagrammatic syntax tailored to the "graded SMCs" of \cref{defn:gradedsmc}.
\end{rem}
A morphism between "graded SMCs" is a functor that preserves all the additional structure.
\begin{defi}[Graded functor]\label{defn:gfunc}
  \AP Given two "graded SMCs" $\catC$ and $\catD$, a \textbf{$\G$-""graded symmetric strict monoidal functor""} $F: \catC \rightarrow \catD$ is a function $F: \mathrm{Ob}(\catC) \rightarrow \mathrm{Ob}(\catD)$ along with a family of functions $F_{\grade}: \catC_{\grade}(X,Y) \rightarrow \catD_{\grade}(FX,FY)$ indexed by $X,Y \in \catC$ and $\grade \in \G$ (although we omit indices in the notation) satisfying
  \begin{gather*}\setlength\arraycolsep{1pt}
    \begin{array}{rlcrlcrlcrl}
      F(t \regrade f) & = t \regrade F(f) & \ \  & F(X \tensor Y) & = FX \tensor FY     & \ \  & F(\id_X)     & = \id_{FX}  & \ \  & F(\swap{X,Y}) & = \swap{FX,FY}. \\
      F(f\comp g)     & = F(f) \comp F(g) & \ \  & F(f\tensor g)  & = F(f) \tensor F(g) & \ \  & F(I_{\catC}) & = I_{\catD} & \ \  &               &
    \end{array}
  \end{gather*}
  In the sequel, we will simply call $F$ a "graded functor". We say $F$ is \textbf{full} (respectively \textbf{faithful}) if $F_{\grade}$ is surjective (respectively injective) for all $X$, $Y$, and $\grade$.
\end{defi}

It will be useful to consider the "SMC" of $\unitG$-graded morphisms inside a "graded SMC" as in e.g.\@ \cite[Definition 8.1]{Katsumata2022} and~\cite[Definition 12]{McDermott2022}.
\begin{prop}\label{prop:underlyingcat}
  \AP Any $\G$-"graded SMC" $\catC$ has an ""underlying SMC"", which we denote with $\catC_{\unitG}$, whose objects are $\mathrm{Ob}(\catC)$ and morphisms are $\unitG$-graded morphisms of $\catC$.
\end{prop}

\begin{exa}[Graded para]\label{exmp:gradedpara}
  \AP A known method to build a "graded category" is a variant of the ""para construction"". We are particularly interested in the "para construction" for symmetric monoidal "actegories" described in~\cite[Sections 2 and 2.1]{CapucciGHR2022}, but instead of building a bicategory, we build a "graded category".

  Given a "symmetric monoidal action" $\act: \G \times \catC \rightarrow \catC$  (\cref{defn:act}), there is a $\G$-"graded SMC" $\Para{\G}{\catC}$ whose "underlying SMC" is $\catC$ and morphisms $X \garrow{\grade} Y$ are $\catC$-morphisms $\grade \act X \rightarrow Y$. It is defined as follows.
  \begin{itemize}
    \item Objects of $\Para{\G}{\catC}$ are those of $\catC$.
    \item The set of $\grade$-graded morphisms from $X$ to $Y$ is $\Para{\G}{\catC}_{\grade}(X,Y) = \catC(\grade \act X, Y)$.
    \item For every morphism between grades $t: \gradeb \rightarrow \grade$, and graded morphism $f: X \garrow{\grade} Y$,
          \begin{equation}\label{eqn:defn:parareg}
            t \regrade f \coloneq \gradeb \act X \xrightarrow{t \act \id_X} \grade \act X \xrightarrow{f} Y.
          \end{equation}
    \item The composition of two graded morphisms $f: X \garrow{\grade} Y$ and $g: Y \garrow{\gradeb} Z$ is given by (we use the symmetry in $\G$ and the isomorphism $\mu$ from \cref{defn:act})\footnote{The symmetry of $\G$ is only used here to swap the two grades; one could avoid it by grading over the reverse monoidal category of $\G$, whose monoidal product is $(\grade,\gradeb) \mapsto \gradeb\grade$. However, we would then require a symmetry to define the monoidal product of graded morphisms below (through the interchanger $\iota$). Thus, we keep $\G$ symmetric since we are interested in the monoidal structure of $\Para{\G}{\catC}$.}
          \begin{equation}\label{eqn:defn:paracomp}
            f \comp g \coloneq ( \grade\gradeb) \act X \cong (\gradeb\grade) \act X \cong \gradeb \act (\grade \act X) \xrightarrow{\id_\gradeb \act f} \gradeb \act Y \xrightarrow{g} Z.
          \end{equation}
    \item The identity morphism $\id_X : X \xrightarrow{\unitG} X$ is $\unitG \act X \cong X \xrightarrow{\id_X} X$.
    \item The monoidal product of objects in $\Para{\G}{\catC}$ is inherited from the monoidal product of objects in $\catC$.
    \item The monoidal product of $f : X \garrow{\grade} Y$ and $g: X' \garrow{\gradeb} Y'$ is given by (we use the interchanger $\iota$ from \cref{defn:act})
          \begin{equation}\label{eqn:defn:paratens}
            f \tensor g \coloneq ( \grade\gradeb) \act (X \otimes X') \cong (\grade \act X) \otimes (\gradeb \act X') \xrightarrow{f \otimes g} Y \otimes Y'.
          \end{equation}
    \item The symmetries are inherited from the symmetries of $\catC$.
  \end{itemize}
  By our hypotheses on the "action" $\act$, the grade $\unitG$ acts trivially (via $\unitG \act X \cong X$), so the "underlying SMC" of $\Para{\G}{\catC}$ (i.e.\@ the "SMC" of $\unitG$-graded morphisms) is isomorphic to $\catC$.

  The main concrete example of this paper is the "graded category" $\BImP$ introduced in \cref{sec:bimp}. It is an instance of this construction: $\BImP = \Para{\op{\FInj}}{\BStoch}$.
\end{exa}

Following the development of monoidal algebra (\cref{sec:background}), we restrict our focus to "graded SMCs" whose objects are finite words, which are more amenable to axiomatisation.
\begin{defi}[Graded prop]\label{defn:gprop}
  \AP A $\G$"-SMC" $\catC$ is a \textbf{$\G$-""graded (coloured) prop@graded prop""} if both $\G$ and $\catC_{\unitG}$ are "props". Note that we require the grading category $\G$ to be a "prop" as well, and not merely the "underlying SMC" $\catC_{\unitG}$: grades are then words over a set of colours, and the morphisms of $\G$ are amenable to the same (direct) syntactic treatment as those of $\catC$, which is essential for the modular axiomatisation of \cref{sec:modularcompleteness}. To emphasise that we work in a "graded prop", we use the symbols $a$, $b$, $c$ for grades, $\word$, $\wordb$, $\wordc$ for objects, concatenation for monoidal products, and $\eword$ for units.
\end{defi}

\begin{exa}
  When $\G$ and $\catC$ are "props", $\Para{\G}{\catC}$ is a "graded prop".
\end{exa}

\subsection{Graded Theories} We introduce graded counterparts for signatures, terms, equations, theories, models, and syntactic categories. We assume throughout that $\intro*\G$ is a "prop" with colours $\G_0$.
\knowledge{\G}{autoref,scope=subsection,also now}

\begin{defi}[Graded term]
  \AP A \textbf{($\G$-)""graded (monoidal) signature@graded signature""} $\sig$ is a set $\sig_0$ of \textbf{colours} and a set $\sig_1$ of \textbf{generators}, together with functions $\ar,\coar:\sig_1 \rightarrow \words{\sig_0}$ and $\grd:\sig_1 \rightarrow \words{\G_0}$ assigning to each generator an \textbf{arity} $\word \in \words{\sig_0}$, a \textbf{coarity} $\wordb\in\words{\sig_0}$, and a \textbf{grade} $a$ in $\words{\G_0}$. We use the notation $g: \word \garrow{a} \wordb \in \sig$ to compactly represent the data of a generator in the signature. We inductively define the set $\intro*\synct{\sig}$ of \textbf{"graded (monoidal) $\sig$-terms@graded term"} built with the generators in $\sig$ and $\G$-morphisms via the rules \cref{eqn:sigterms:trivial,eqn:sigterms:SMC,eqn:sigterms:seqcomp,eqn:sigterms:parcomp,eqn:sigterms:regrade}. The rules also define the (co)arities and grades of such terms.\\
  \begin{minipage}{0.49\textwidth}
    \begin{gather}
      \label{eqn:sigterms:trivial}
      \scalebox{0.85}{\begin{bprooftree}
          \AxiomC{$g: \word \garrow{a} \wordb \in \sig$}
          \UnaryInfC{$g: \word \garrow{a} \wordb \in \synct{\sig}$}
        \end{bprooftree}}\\
      \label{eqn:sigterms:seqcomp}
      \scalebox{0.85}{\begin{bprooftree}
          \AxiomC{$f: \word \garrow{a} \wordb \in \synct{\sig}$}
          \AxiomC{$\!\!\!\!\!\!g: \wordb \garrow{b} \ell \in \synct{\sig}$}
          \BinaryInfC{$f\comp g : \word \garrow{ab} \ell \in \synct{\sig}$}
        \end{bprooftree}}
    \end{gather}
  \end{minipage}
  \begin{minipage}{0.5\textwidth}
    \begin{gather}\label{eqn:sigterms:regrade}
      \scalebox{0.85}{\begin{bprooftree}
          \AxiomC{$t: b \rightarrow a \in \G$}
          \AxiomC{$\!\!\!\!\!\!f: \word \garrow{a} \wordb \in \synct{\sig}$}
          \BinaryInfC{$t \regrade f : \word \garrow{b} \wordb \in \synct{\sig}$}
        \end{bprooftree}}\\\label{eqn:sigterms:parcomp}
      \scalebox{0.85}{\begin{bprooftree}
          \AxiomC{$f: \word \garrow{a} \wordb \in \synct{\sig}$}
          \AxiomC{$\!\!\!\!\!\!\!\!f': \word' \garrow{b} \wordb' \in \synct{\sig}$}
          \BinaryInfC{$f\tensor f' : \word\word' \garrow{ab} \wordb\wordb' \in \synct{\sig}$}
        \end{bprooftree}}
    \end{gather}
  \end{minipage}
  \begin{equation}
    \label{eqn:sigterms:SMC}
    \scalebox{0.85}{\begin{bprooftree}
        \AxiomC{$\phantom{\sig_0}$}
        \UnaryInfC{$\id_{\eword}: \eword \garrow{\eword} \eword \in \synct{\sig}$}
      \end{bprooftree}}\quad
    \scalebox{0.85}{\begin{bprooftree}
        \AxiomC{$\col \in \sig_0$}
        \UnaryInfC{$\id_{\col}: \col \garrow{\eword} \col \in \synct{\sig}$}
      \end{bprooftree}}\quad
    \scalebox{0.85}{\begin{bprooftree}
        \AxiomC{$\col,\colb \in \sig_0$}
        \UnaryInfC{$\swap{\col,\colb}: \col\colb \garrow{\eword} \colb\col \in \synct{\sig}$}
      \end{bprooftree}}
  \end{equation}
  These operations on "graded terms" (simply called "terms@@GRD" in the sequel) correspond straightforwardly with the structure of a "graded SMC". We can see any generator as a "term@@GRD" via \cref{eqn:sigterms:trivial}; there are identity "terms@@GRD" and symmetries thanks to \cref{eqn:sigterms:SMC}; we can compose "terms@@GRD" in sequence and in parallel with \cref{eqn:sigterms:seqcomp} and \cref{eqn:sigterms:parcomp} respectively; and we can "regrade" a term by any $\G$-morphism with \cref{eqn:sigterms:regrade}.
\end{defi}
We introduce a two-dimensional graphical representation for "graded terms" as string diagrams. We denote a term $f: \word \garrow{a} \wordb$ by a box with a wire on the left labelled $\word$ and a wire on the right labelled $\wordb$, and a wire on the bottom labelled $a$:$\scalebox{0.75}{\input{fnma.tikz}}$. As before, we write the "terms@@GRD" $\id_{\eword}$ as $\inlinefig{0.8}{idzer}$, $\id_{\col}$ as $\inlinefig{0.8}{wirecol}$, and $\swap{\col,\colb}$ as $\inlinefig{0.6}{swapcol}$. We write sequential composition, monoidal product, and "regrading" as below. The colouring on $\G$-morphisms is only a reading aid, and it can safely be ignored. "Graded terms" are read from left to right, then top to bottom, while $\G$-morphisms are read from bottom to top then left to right. When the grade of a term is $\eword$, we omit the grading wire altogether, so that $\eword$-"graded terms" are drawn exactly as the string diagrams of \cref{sec:background}.
\begin{minipage}{0.55\textwidth}
  \begin{gather*}
    \inlinefig{0.6}{fnma} \comp \inlinefig{0.6}{gmlb} = \inlinefig{0.6}{seqcompfg}\\
    \inlinefig{0.6}{fnma} \tensor \inlinefig{0.6}{gnpmpb} = \inlinefig{0.6}{parcompfg}
  \end{gather*}
\end{minipage}
\begin{minipage}{0.44\textwidth}
  \begin{gather*}
    \inlinefig{0.6}{tba} \regrade \inlinefig{0.6}{fnma} = \inlinefig{0.6}{regradedfnmb}
  \end{gather*}
\end{minipage}

\begin{rem}
  The \textit{grading} wire coming orthogonally to the (co)arity wire is inspired by the string diagrams used for the "para construction" in e.g.\@ \cite{CapucciGHR2022,Cruttwell2022,Gavranovic2024}, although our grading wires come from the bottom rather than the top (\cref{rem:slidingfrombelow} explains why this convention is more natural in our setting).\end{rem}

\begin{defi}[Theory]
  \AP A \textbf{$\G$-""graded monoidal theory""} is a tuple $\th = (\sig,\axh)$, where $\sig$ is a "graded signature", and $\axh$ is a set of equations,\footnote{The superscript $\mathsf{h}$ stands for `horizontal': these are equations between "graded terms", which we draw horizontally. Later on, $(\sigv,\axv)$ will denote a presentation of the grading category $\G$, whose terms we draw vertically ($\mathsf{v}$).} called \textbf{axioms}, between "graded $\sig$-terms@graded term" whose arities, coarities, and grades coincide, i.e.\@ pairs $(s,t) \in \synct{\sig}\times \synct{\sig}$ with $\ar(s) = \ar(t)$, $\coar(s) = \coar(t)$, and $\grd(s) = \grd(t)$. We denote such equations by $s=t$.
\end{defi}

\begin{defi}[Model]
  Let $\th = (\sig,\axh)$ be a $\G$-"graded monoidal theory", $\catC$ a $\G$"-SMC". Given an assignment $M: \sig \rightarrow \catC$ sending each colour $\col \in \sig_0$ to an object $M\col \in \catC$ and each generator $g: \word \garrow{a} \wordb \in \sig$ to a morphism in $\catC_a(\bigotimes_{\col\in\word}M\col, \bigotimes_{\col\in\wordb}M\col)$, we can canonically extend $M$ to all $\synct{\sig}$ because each term formation rule \cref{eqn:sigterms:SMC,eqn:sigterms:seqcomp,eqn:sigterms:parcomp,eqn:sigterms:regrade} corresponds to part of the structure in a "graded SMC". A \textbf{model} of $\th$ valued in $\catC$ is such an assignment $M$ satisfying $M(s) = M(t)$ for all axioms $s = t \in \axh$.\end{defi}

As for (ungraded) monoidal algebra in \cref{sec:background}, we will see models as functors, so we need to define a "syntactic category@@GRD" for a "theory@@GRD" $\th$. As usual, its morphisms will be "terms@@GRD" quotiented by a congruence. In this case, the congruence is generated by the axioms in $\axh$ and the necessary equations to satisfy the laws of a $\G$"-SMC".

\begin{defi}[Closure]\label{defn:closure}
  Given a "theory@@GRD" $\th = (\sig,\axh)$, we define the \textbf{closure} of $\axh$, denoted $\closed{\axh}$, to be the smallest congruence on $\synct{\sig}$ (i.e.\@ equivalence relation compatible with $\comp$, $\tensor$, and $\regrade$, see \cref{fig:infrulessyncat}) containing $\axh$ and equations \cref{eqn:th:assocseq,eqn:th:neutralid,eqn:th:assocpar,eqn:th:neutralidzer,eqn:th:seqregrade,eqn:th:parregrade,eqn:th:parregradepar,eqn:th:trivialregrade,eqn:th:overthewire,eqn:th:interchangeupto,eqn:th:swapswap,eqn:th:coherentid,eqn:th:coherentswaps} in \cref{fig:axiomssyncat}, where terms are universally quantified.

  \setcounter{theqn}{0}
  \begin{figure}
    \begin{subfigure}{0.41\textwidth}
      \begin{gather}\label{eqn:th:assocseq}
        \inlinefig{0.54}{assocseqcomp} \stepcounter{theqn}\tag{G\arabic{theqn}}\\
        \label{eqn:th:neutralid}
        \inlinefig{0.54}{neutralid} \stepcounter{theqn}\tag{G\arabic{theqn}}\\
        \label{eqn:th:neutralidzer}
        \inlinefig{0.54}{neutralidzer} \stepcounter{theqn}\tag{G\arabic{theqn}}
      \end{gather}
    \end{subfigure}
    \begin{subfigure}{0.26\textwidth}
      \begin{gather}\label{eqn:th:seqregrade}
        \inlinefig{0.54}{seqregrade} \stepcounter{theqn}\tag{G\arabic{theqn}}
      \end{gather}
    \end{subfigure}
    \begin{subfigure}{0.31\textwidth}
      \begin{gather}\label{eqn:th:parregradepar}
        \inlinefig{0.54}{parregradepar} \stepcounter{theqn}\tag{G\arabic{theqn}}
      \end{gather}
    \end{subfigure}\\[0.2em]
    \begin{subfigure}{0.41\textwidth}
      \begin{equation}\label{eqn:th:assocpar}
        \inlinefig{0.54}{assocparcomp} \stepcounter{theqn}\tag{G\arabic{theqn}}
      \end{equation}
    \end{subfigure}
    \begin{subfigure}{0.32\textwidth}
      \begin{equation}\label{eqn:th:parregrade}
        \inlinefig{0.54}{parregrade} \stepcounter{theqn}\tag{G\arabic{theqn}}
      \end{equation}
    \end{subfigure}
    \begin{subfigure}{0.24\textwidth}
      \begin{equation}\label{eqn:th:trivialregrade}
        \inlinefig{0.54}{trivialregrade} \stepcounter{theqn}\tag{G\arabic{theqn}}
      \end{equation}
    \end{subfigure}\\[0.2em]
    \begin{subfigure}{0.38\textwidth}
      \begin{equation}\label{eqn:th:overthewire}
        \inlinefig{0.54}{overthewire} \stepcounter{theqn}\tag{G\arabic{theqn}}
      \end{equation}
    \end{subfigure}
    \begin{subfigure}{0.45\textwidth}
      \begin{equation}\label{eqn:th:interchangeupto}
        \inlinefig{0.54}{interchangeupto} \stepcounter{theqn}\tag{G\arabic{theqn}}
      \end{equation}
    \end{subfigure}\\[0.2em]
    \begin{subfigure}{0.3\textwidth}
      \begin{equation}\label{eqn:th:swapswap}
        \inlinefig{0.54}{swapswap2} \stepcounter{theqn}\tag{G\arabic{theqn}}
      \end{equation}
    \end{subfigure}
    \begin{subfigure}{0.3\textwidth}
      \begin{equation}\label{eqn:th:coherentid}
        \inlinefig{0.54}{coherentid} \stepcounter{theqn}\tag{G\arabic{theqn}}
      \end{equation}
    \end{subfigure}
    \begin{subfigure}{0.35\textwidth}
      \begin{equation}\label{eqn:th:coherentswaps}
        \inlinefig{0.54}{coherentswaps} \stepcounter{theqn}\tag{G\arabic{theqn}}
      \end{equation}
    \end{subfigure}

    \caption{Laws of "graded symmetric strict monoidal categories".}\label{fig:axiomssyncat}
  \end{figure}
\end{defi}

\begin{defi}[Syntactic category]\label{defn:syntacticcat}
  \AP The \textbf{""syntactic category@@GRD""} of a $\G$-"graded monoidal theory" $\th = (\sig,\axh)$ is a $\G$-"graded prop" $\intro*\sync{\th}$ defined as follows.
  The set of morphisms ${\sync{\th}}_{a}(\word,\wordb)$ is the set of $\sig$-terms of arity $\word$, coarity $\wordb$, and grade $a$, quotiented by the congruence $\closed{\axh}$. Identities, symmetries, sequential composition, monoidal product, and "regrading" are defined via the corresponding term formation rules~\cref{eqn:sigterms:trivial,eqn:sigterms:SMC,eqn:sigterms:seqcomp,eqn:sigterms:parcomp,eqn:sigterms:regrade}. The axioms in \cref{fig:axiomssyncat} ensure that all the laws of "graded SMCs" hold.

  We say that $\th$ is a \textbf{(graded) presentation} of a $\G$-"graded prop" $\catC$ (or that $\th$ presents $\catC$) if $\sync{\th}$ and $\catC$ are isomorphic as "graded props", or equivalently, if there is a bijective-on-objects, fully faithful, "graded functor" between them.
\end{defi}

\begin{rem}\label{rem:verticalsyntax}
  Strictly speaking, our notion of graded term is not purely syntactic, since the "regrading" rule~\cref{eqn:sigterms:regrade} involves arbitrary $\G$-morphisms (which may or may not be freely obtained from a signature). This means that some equations are automatically imposed, e.g.\@ if $s\compsimple t = r$ in $\G$, then   $\inlinefig{0.65}{fininjeqns2}\stepcounter{theqn}$ holds in $\sync{\th}$ by reflexivity. We provide this formulation for the sake of generality. However, if $\G$ were presented by a "monoidal theory" $(\sigv,\axv)$, we could have treated graded terms as pure syntax, by requiring that $t$ appearing in the "regrading" rule~\cref{eqn:sigterms:regrade} is a $\sigv$-term rather than a morphism (= an equivalence class of $\sigv$-terms). With this setup, we still obtain a definition of syntactic category for a "graded theory" equivalent to the one of \cref{defn:syntacticcat}, by requiring that graded $\sig$-terms are further quotiented by the equations of $\axv$ and the laws of "SMCs" applied to the $\sigv$-terms appearing in the $\sig$-terms as grading. In our leading examples (\cref{sec:bimp,sec:diagraminterventions}), $\G$ is indeed obtained from a "monoidal theory". \end{rem}

The axioms in \cref{fig:axiomssyncat} are direct diagrammatic translations of properties of "graded SMCs". Thus, they are \textit{sound} in any "graded SMC", and we get the following result.
\begin{thm}[Soundness]\label{lem:soundness}
  Models of $\th$ valued in $\catC$ are in 1-to-1 correspondence with $\G$-"graded functors" $\sync{\th} \rightarrow \catC$.
\end{thm}
\begin{proof}
  Let $M$ be a model of $\th$ valued in $\catC$.
  Recall that, by hypothesis of $M$ being a model, the equations in $\axh$ hold after applying $M$. Moreover, all the equations in \cref{fig:axiomssyncat} hold after applying $M$ because $\catC$ is a "graded SMC". Since equality of morphisms in $\catC$ is a congruence, we conclude that all the equations in $\closed{\axh}$ hold after applying $M$, because $\closed{\axh}$ was defined as the smallest congruence containing $\axh$ and the equations in \cref{fig:axiomssyncat}. Hence, the extended assignment $M : \synct{\sig} \rightarrow \catC$ factors through a  $\G$-"graded functor" $\sync{\th} \rightarrow \catC$.

  To summarise the first part of the proof, we have that any model satisfying the axioms in $\axh$ satisfies all the equations in $\closed{\axh}$. This justifies the name of the theorem: we proved soundness of the closure in all models.

  Conversely, given a $\G$-"graded functor" $M: \sync{\th} \rightarrow \catC$, it is easy to check that the assignment $\sig \rightarrow \catC$ defined by restricting $M$ to equivalence classes of generators yields a model of $\th$.
\end{proof}

Completeness is shown by considering the model of $\th$ that corresponds to $\id: \sync{\th} \rightarrow \sync{\th}$.
\begin{thm}[Completeness]\label{thm:completeness}
  If $s=t$ is satisfied in all models of $(\sig,\axh)$, then $s=t \in \closed{\axh}$.
\end{thm}
\begin{proof}
  We define the syntactic model of $\th$ to be the identity functor $\sync{\th} \rightarrow \sync{\th}$. It corresponds to an assignment sending a generator in $\sig$ to the equivalence class of that generator, seen as a morphism in $\sync{\th}$. The extension of this assignment to $\synct{\sig}$ sends a $\sig$-term to its equivalence class. Therefore, an equation between $\sig$-terms $s=t$ is satisfied in the syntactic model if and only if $s$ and $t$ are equal modulo $\closed{\axh}$. In particular, if all models (including the syntactic model) satisfy an equation $s=t$, then $s=t \in \closed{\axh}$.
\end{proof}

The definitions and results of this section culminate in a conservative extension of (classical) "monoidal theories" and "props". Indeed, let us suppose that $\eword$ is \textbf{strictly terminal} in $\G$, namely, $\eword$ is a terminal object of $\G$ and there are no morphisms $\eword \rightarrow \word$ for $\word \neq \eword$ (in particular, $\id_{\eword}$ is the only morphism $\eword \rightarrow \eword$). We note that $\eword$-graded terms in $\sync{\th}$ are precisely those built out of $\eword$-graded generators, identities, and symmetries, possibly "regraded" by $\id_{\eword}$. Indeed, in the formation rules \cref{eqn:sigterms:trivial,eqn:sigterms:SMC,eqn:sigterms:seqcomp,eqn:sigterms:parcomp,eqn:sigterms:regrade}, the grade of a term can go from $\word \neq \eword$ to $\eword$ only via rule \cref{eqn:sigterms:regrade} with a $\G$-morphism $t: \eword \rightarrow \word$, and we assumed no such $t$ exists; moreover, the only possible "regrading" of an $\eword$-graded term is by $\id_{\eword}$, which disappears by \cref{eqn:th:trivialregrade}. Furthermore, if we restrict \cref{fig:axiomssyncat} to $\eword$-graded terms (those that do not have a dangling wire at the bottom, once "regradings" by $\id_{\eword}$ are removed), we recover the usual laws of SMCs: the "regradings" by symmetries in \cref{eqn:th:overthewire,eqn:th:interchangeupto} become trivial, the axioms \cref{eqn:th:seqregrade}, \cref{eqn:th:parregradepar}, \cref{eqn:th:parregrade}, and \cref{eqn:th:trivialregrade} become vacuous, and the remaining axioms are those of \cref{fig:axiomsvyncat}. We infer that the "underlying prop" of a "syntactic category@@GRD" $\sync{\th}$ is the "syntactic category@@MON" (\cref{defn:monthsyncat}) for the "monoidal theory" given by the $\eword$-graded generators and $\eword$-graded equations in $\th$.

We formally prove this result so that we can rely on it when defining the "graded theory" for the "para construction" in \cref{sec:modularcompleteness}.
\begin{lem}\label{lem:zerogradedpart}
  Suppose that $\eword$ is strictly terminal in $\G$, and let $\th^{\eword} \coloneq (\sig^{\eword},E^{\eword})$ be the (non-graded) "monoidal theory" obtained by taking only the $\eword$-graded generators and equations in a "graded theory" $\th = (\sig,\axh)$. The "underlying prop" of $\sync{\th}$ is the "syntactic category@@MON" of $\th^{\eword}$.
\end{lem}
\begin{proof}
  Formally, $\th^{\eword}$ is defined by
  \[\sig^{\eword} = \{ g: \word \rightarrow \wordb \mid g: \word \garrow{\eword} \wordb \in \sig \} \quad E^{\eword} = \{ f = f' \mid f = f' \in \axh \text{ and } f,f' : \word \garrow{\eword} \wordb\}.\]
  Let $\vync{\th^{\eword}}$ be the "prop" generated by $\th^{\eword}$ and $U\sync{\th}$ denote the "underlying prop" of $\sync{\th}$. Also, let $i$ be the map sending every generator in $\sig^{\eword}$ to its $\eword$-graded counterpart in $U\sync{\th}$. We will show that $i$ canonically extends to an isomorphism of "props" $\vync{\th^{\eword}} \cong U\sync{\th}$.

  We also denote with $i: \vynct{\sig^{\eword}} \rightarrow U\sync{\th}$ the extension of $i$ to monoidal $\sig^{\eword}$-terms. It sends a term to its equivalence class modulo $\closed{\axh}$. To check that $i$ factors through $\vync{\th^{\eword}}$, it suffices to show that the equations in $E^{\eword}$ hold in $U\sync{\th}$ after applying $i$. This is immediate from the fact that the equations in $E^{\eword}$ are contained in $\axh$.

  Next, in order to show the factored functor $i^*: \vync{\th^{\eword}} \rightarrow U\sync{\th}$ is an isomorphism, we need to show that any $\eword$-graded term of $\sync{\th}$ is in the image of $i^*$ (fullness), and that whenever $f = g$ is in the closure $\closed{\axh}$ for $\eword$-graded terms $f$ and $g$, then $f=g$ is in $\closed{E^{\eword}}$ (faithfulness).

  Recall that sequential and parallel compositions and "regrading" of terms of grades other than $\eword$ cannot yield $\eword$-graded terms, thus all $\eword$-graded terms are obtained by sequential and parallel composition of $\eword$-graded generators, identities, and symmetries. This implies that $i^*$ is full. Moreover, it implies that equational reasoning for $\eword$-graded terms is entirely independent of other terms. Explicitly, the set of equations between $\eword$-graded terms in $\closed{\axh}$ is the smallest congruence (only compatible with sequential and parallel composition) that contains the axioms in $\closed{\axh}$ between $\eword$-graded terms and the equations in \cref{fig:axiomssyncat} between $\eword$-graded terms. Finally, we compare the restriction of \cref{fig:axiomssyncat} to $\eword$-graded terms with the axioms of "SMCs" \cref{fig:axiomsvyncat} to find that they generate the same congruence. Therefore, for any $\eword$-graded terms $f$ and $g$, $f=g \in \closed{\axh}$ if and only if $f= g \in \closed{E^{\eword}}$. This implies $i^*$ is faithful.
\end{proof}

\subsection{Erasing Grades}\label{sec:collapse}
The construction of the "underlying SMC" (\cref{prop:underlyingcat}) erases grades from a $\G$"-SMC" by forgetting about all morphisms whose grades are not $\unitG$, and thus loses all information about the grading category. In this section, we describe another way to obtain an "SMC" from a $\G$"-SMC" by keeping graded morphisms, but forgetting their grades and identifying morphisms that can be "regraded" to each other. This construction (\cref{defn:collapse}) is similar to the one described in \cite[Proposition~5.6]{Earnshaw2026}.

Our goal here is to quotient a $\G$"-SMC" by considering its morphisms "up to regrading". Recall that for $f: X \garrow{\grade} Y$ and $g : X \garrow{\gradeb} Y$, we have $f \utg g$ whenever $f = t \regrade g$ or $g = t \regrade f$. This relation is not an equivalence because the grading category $\G$ might not have inverses, so if $g = t \regrade f$ and $g = t' \regrade f'$, nothing guarantees that we can find $t''$ such that $f' = t'' \regrade f$ or vice versa. We will use the equivalence closure of $\utg$.
\begin{defi}[$\collcong$]\label{defn:collcong}
  Given a $\G$"-SMC" $\catC$, for each pair of objects $X,Y$ in $\catC$, $\intro*\collcong$ is the smallest equivalence relation on the set $\coprod_{\grade \in \mathbb{G}}\catC_{\grade}(X,Y)$ such that $f \collcong g$ whenever $g = t \regrade f$. Equivalently, $\collcong$ is the transitive closure of equality "up to regrading" (reflexivity already holds by "regrading" with an identity and symmetry by definition of $\utg$).
\end{defi}
In order to quotient by $\collcong$ and obtain an "SMC" (with composition and monoidal product), we need to show that $\collcong$ is compatible with respect to $\comp$ and $\tensor$.
\begin{lem}\label{lem:collcongcongruence}
  The family $\collcong$ is a congruence, namely,
  \begin{enumerate}
    \item if $f \collcong f'$ and $g \collcong g'$, then $f \comp g \collcong f' \comp g'$, and
    \item if $f \collcong f'$ and $g \collcong g'$, then $f \tensor g \collcong f' \tensor g'$.
  \end{enumerate}
\end{lem}
\begin{proof}
  Let us only show compatibility for sequential composition, as we can follow exactly the same proof for monoidal product. It is enough to show that $\utg$ is compatible because the transitive closure of a compatible reflexive relation is compatible (see e.g.\@ \cite[Theorem~6]{Chajda1975}).

  Without loss of generality, suppose that $f = t \regrade f'$. We consider two cases for $g$. If $g = s \regrade g'$, then $f \comp g = (t \regrade f') \comp (s \regrade g') = (t \otimes s) \regrade (f' \comp g')$, so $f \comp g \collcong f' \comp g'$. Otherwise, if $g' = s \regrade g$, then we can use the previous case (with $g = \id \regrade g$) to get $f \comp g \collcong f' \comp g$, then a symmetric argument to get $f' \comp g \collcong f' \comp g'$. We conclude by transitivity of $\collcong$.
\end{proof}
\begin{defi}[$\Coll$]\label{defn:collapse}
  Given a $\G$"-SMC" $\catC$, we define an "SMC" called the \textbf{""regrading collapse""} of $\catC$, denoted $\intro*\Coll(\catC)$, as follows:
  \begin{itemize}
    \item Objects of $\Coll(\catC)$ are those of $\catC$.
    \item For any objects $X, Y$ in $\catC$, $\Coll(\catC)(X,Y) = (\coprod_{\grade \in \G} \catC_{\grade}(X,Y))/\collcong$. We will write $\intro*\collclass{f} \in \Coll(\catC)(X,Y)$ for the equivalence class of a graded morphism $f: X \garrow{\grade} Y$.
    \item Composition is defined by $\collclass{f}\compsimple \collclass{g} = \collclass{f\comp g}$. Monoidal product is defined on objects like in $\catC$ (with the same unit) and on morphisms by $\collclass{f}\otimes \collclass{g} = \collclass{f\tensor g}$. These definitions are independent of the choices of representatives because $\collcong$ is a congruence (\cref{lem:collcongcongruence}).
    \item Identities and symmetries are induced from $\catC$. Namely, the identity on $X$ is $\collclass{\id_X}$ and the symmetry between $X$ and $Y$ is $\collclass{\swap{X,Y}}$.
  \end{itemize}
\end{defi}
Following \cref{lem:zerogradedpart} where we gave a "monoidal theory" to present the "underlying prop" of a "graded SMC" with a given presentation, we now describe a syntactic counterpart to the "regrading collapse". In short, given a "graded monoidal theory", we erase the grading wires on \emph{all} generators, retaining their input and output wires, instead of keeping only the generators with no grading wires.

\begin{defi}[Flattened theory]\label{defn:flatth}
  Let $\th=(\sig,\axh)$ be a $\G$-"graded monoidal theory". Its \textbf{""flattened theory""} is the classical "monoidal theory" $\intro*\flatt{\th}=(\flatt{\sig},\flatt{E})$ defined as follows. The colours of $\flatt{\sig}$ are those of $\sig$ and there is a "flattened generator" in $\flatt{\sig}$ with no grade for each (graded) generator of $\th$:
  \[ \flatt{\sig}_1 \coloneq \{\flatt{g}:\word\rightarrow\wordb \mid g:\word\garrow{a}\wordb\in\sig \text{ for some } a\}. \]
  To list the axioms of $\flatt{E}$, we recursively extend "flattening" to "graded terms" of $\th$ by
  \begin{gather}\label{eqn:flattenterms}
    \flatt{\id_{\word}} = \id_{\word}, \quad
    \flatt{\swap{\word,\wordb}} = \swap{\word,\wordb}, \quad
    \flatt{(f\comp g)}=\flatt{f}\compsimple \flatt{g}, \quad
    \flatt{(f\tensor g)}=\flatt{f}\tensor\flatt{g}, \quad
    \flatt{(t\regrade f)}=\flatt{f}.
  \end{gather}
  Then, $\flatt{E}$ contains the same axioms as $\axh$ where "graded terms" are replaced by their "flattening":
  \[ \flatt{E} \coloneq \{\flatt{f}=\flatt{g} \mid f=g\in\axh\}.\]
\end{defi}

\begin{thm}\label{thm:collapsepresentation}
  If a $\G$-"graded monoidal theory" $\th$ presents a $\G$-"graded prop" $\catC$, then the "flattened theory" $\flatt{\th}$ presents $\Coll(\catC)$.
\end{thm}
\begin{proof}
  By hypothesis, we have an isomorphism of "graded props" $\sem{-} : \sync{\th} \rightarrow \catC$, and we will construct an isomorphism of "props" $F: \vync{\flatt{\th}} \rightarrow \Coll(\catC)$. We use functorial semantics by defining a model of $\flatt{\th}$ (\cref{defn:monthmodel}) and showing the induced functor is an isomorphism.

  Let the mapping $F$ send a generator $\flatt{g} \in \flatt{\sig}$ to $\collclass{\sem{g}}$, the equivalence class of $\sem{g}$ in the "regrading collapse" of $\catC$. We can show by structural induction that for any "graded term" $f$ in $\th$, we have $F(\flatt{f}) = \collclass{\sem{f}}$, where $F$ is recursively extended to terms in $\flatt{\th}$. The base case of generators is by definition of $F$. The base cases of identities and symmetries are dealt with by the fact that $\collclass{\id_\word}$ is the identity on $\word$ and $\collclass{\swap{\word,\wordb}}$ is the symmetry between $\word$ and $\wordb$ in $\Coll(\catC)$. For the inductive cases, we have the following derivations (\IH\ and \FS\ denote steps using the inductive hypothesis and functoriality of $\sem{-}$ respectively):
  \begin{align*}
    F(\flatt{(f \comp g)})    & \stackrel{\cref{eqn:flattenterms}}{=} F(\flatt{f} \compsimple \flatt{g}) = F(\flatt{f}) \compsimple F(\flatt{g}) \stackrel{\IH}{=} \collclass{\sem{f}} \compsimple \collclass{\sem{g}} \stackrel{\text{Lem.\@ } \ref{lem:collcongcongruence}}{=} \collclass{\sem{f}\comp \sem{g}} \stackrel{\FS}{=} \collclass{\sem{f \comp g}} \\
    F(\flatt{(f \tensor g)})  & \stackrel{\cref{eqn:flattenterms}}{=} F(\flatt{f} \otimes \flatt{g}) = F(\flatt{f}) \otimes F(\flatt{g}) \stackrel{\IH}{=} \collclass{\sem{f}} \otimes \collclass{\sem{g}} \stackrel{\text{Lem.\@ } \ref{lem:collcongcongruence}}{=} \collclass{\sem{f}\tensor \sem{g}} \stackrel{\FS}{=} \collclass{\sem{f \tensor g}}         \\
    F(\flatt{(t \regrade f)}) & \stackrel{\cref{eqn:flattenterms}}{=} F(\flatt{f}) \stackrel{\IH}{=} \collclass{\sem{f}} \stackrel{\text{Def.\@ }\ref{defn:collcong}}{=} \collclass{t \regrade \sem{f}} \stackrel{\FS}{=} \collclass{\sem{t \regrade f}}
  \end{align*}
  Therefore, whenever $\flatt{f} = \flatt{g}$ is an axiom in $\flatt{E}$, then $F(\flatt{f}) = F(\flatt{g})$ because $f = g \in \axh$ implies $\collclass{\sem{f}} = \collclass{\sem{g}}$. This means that $F$ is a model and we use the same notation for the corresponding symmetric strict monoidal functor $F:\vync{\flatt{\th}} \rightarrow \Coll(\catC)$. It remains to show it is an isomorphism of "props".

  First note that it is identity-on-objects because $\th$ and $\flatt{\th}$ have the same colours, hence their "syntactic categories@@MON" have the same objects, and $\catC$ also has the same objects by the isomorphism $\sem{-}$ and so does $\Coll(\catC)$ by construction. Next, $F$ is full because any morphism in $\Coll(\catC)$ is of the form $\collclass{f'}$ for a $\catC$-morphism $f'$ and since $\sem{-}$ is full, we can find a "graded term" $f$ such that $F(\flatt{f}) = \collclass{\sem{f}} = \collclass{f'}$.

  Finally, for faithfulness, we rely on the fact that any term in $\flatt{\th}$ can be seen as the "flattening" of a "graded term" in $\th$. Indeed, one can reproduce all steps for constructing a given term in $\flatt{\th}$ inside the "graded theory" $\th$, and the "flattening" of the resulting "graded term" will be equal to the given term. Then, without loss of generality, we can suppose $F(\flatt{f}) = F(\flatt{g})$ for two "graded terms" $f$ and $g$. We obtain $\collclass{\sem{f}} = \collclass{\sem{g}}$, namely $\sem{f} \collcong \sem{g}$. It follows since $\sem{-}$ is an isomorphism that preserves "regrading" that $f \collcong g$ as well. We conclude that $\flatt{f} = \flatt{g}$ by (repeatedly) using the equation $\flatt{(t \regrade h)} = \flatt{h}$ from \cref{eqn:flattenterms} and the fact that $\collcong$ is the transitive closure of $\utg$.
\end{proof}

\section{Modular Completeness for Graded Theories}\label{sec:modularcompleteness}

In monoidal algebra it is typical to study a certain "prop" $\catC$ of interest through the lens of string diagrammatic calculi. Within this perspective, a key result is finding a "monoidal theory" $\th$ presenting $\catC$ (in the sense of \cref{defn:monthsyncat}). That would mean that we can reason about $\catC$-morphisms entirely in algebraic terms, by manipulation of the string diagrams of $\vync{\th}$ via the axioms of $\th$.

In the case of a "graded prop", there are two categories at stake: the base category $\catC$, and the grading category $\G$. It is natural to ask the following question: can we systematically derive a presentation of the "graded prop" from a known presentation of both $\catC$ and $\G$? In this section we outline conditions under which such a question may be answered positively.

\AP As above, we fix a "prop" $\intro*\G$, and we assume it has a presentation $(\sigv,\axv)$.
\knowledge{\G}{autoref,scope=section,also now}
We further suppose that $\eword$ is strictly terminal, i.e.\@ $\eword$ is terminal in $\G$ and there is no term $\scalebox{0.6}{\begin{tikzpicture}
	\begin{pgfonlayer}{nodelayer}
		\node [style=none] (8) at (0, 0.75) {};
		\node [style=gletter] (9) at (0, -0.25) {$t$};
	\end{pgfonlayer}
	\begin{pgfonlayer}{edgelayer}
		\draw [style=green] (9) to (8.center);
	\end{pgfonlayer}
\end{tikzpicture}
}$ of type $\eword \rightarrow \word$ with $\word \neq \eword$. Note that this assumption only concerns the grading category $\G$: the base category $\catC$ may well have morphisms out of its unit $\eword$, such as the states $0 \rightarrow n$ of $\BStoch$ in \cref{sec:application}. We will show that, given a symmetric monoidal action $\act: \G \times \catC \rightarrow \catC$, a presentation of a "prop" $\catC$ can be adapted into a graded presentation of $\Para{\G}{\catC}$. Fix a monoidal presentation $\thC = (\sigC,\axC)$ of $\catC$.

Since both $\G$ and $\catC$ are props---in particular, objects are generated by concatenation of colours in $\sigv_0$ and $\sigC_0$ respectively---the action $\act$ is completely determined on objects by the values of $\col \act \eword$ for each $\col \in \G_0$. Let $k_{\col} \coloneq \col \act \eword$. We show that $a \act \word = \word \otimes \bigotimes_{\col \in a}k_{\col}$.
For brevity, we will write $k^a$ for $\bigotimes_{\col \in a}k_{\col}$.

\begin{lem}\label{lem:actiondetermined}
  Let $\act: \G \times \catC \rightarrow \catC$ be a symmetric monoidal action of "props" and $k_{\col} = \col \act \eword$ for each $\col \in \G_0$. For any $a \in \words{{\sigv_0}}$ and $\word \in \words{{\sigC_0}}$, $a \act \word = \word \otimes k^a$.
\end{lem}
\begin{proof}
  This relies on the following three properties that are true in any symmetric monoidal action (see~\cite{CapucciGHR2022}); they are instances of the isomorphisms $\kappa$, $\iota$, and $\eta$ of \cref{defn:act}, respectively.
  \[ A \act (X \otimes Y) \stackrel{*}{\cong} X \otimes (A \act Y) \qquad \quad (A\otimes B) \act (X\otimes Y) \stackrel{**}{\cong} (A \act X) \otimes (B \act Y) \qquad \quad I \act X \stackrel{***}{\cong} X\]
  The proof then proceeds in three steps.
  \begin{enumerate}
    \item For any $a \in \words{{\sigv_0}}$ and $\word \in \words{{\sigC_0}}$, we have $a \act \word = a \act (\word \otimes \eword) \stackrel{*}{=} \word \otimes (a\act \eword)$.
    \item For any $a \in \words{{\sigv_0}}$ and $\col \in \sigv_0$, we have $(a\col) \act \eword = (a\otimes \col) \act (\eword\otimes\eword) \stackrel{**}{=} (a \act \eword) \otimes (\col \act \eword)$. Thus, by initialising an induction with the base cases $a=\col$ with $\col \act \eword = k_{\col}$, we find that $a \act \eword = k^a$ for all $a \neq \eword$.
    \item Since $\eword \act \word \stackrel{***}{=} \word$, we conclude that for all $a \in \words{{\sigv_0}}$ and $\word \in \words{{\sigC_0}}$, $a \act \word = \word \otimes k^a$.\qedhere
  \end{enumerate}
\end{proof}
Now, a graded morphism $f:\word \garrow{a} \wordb$ in $\Para{\G}{\catC}$ corresponds, by definition, to a morphism $f_{\eword}:a \act \word \rightarrow \wordb$ in $\catC$. The presentation of $\catC$ allows us to represent $f_{\eword}$ as a string diagram with two input wires labelled $\word$ and $k^a$ and an output wire labelled $\wordb$. To obtain a graded string diagram that represents $f$, we will simply bend down the second input wire of $f_{\eword}$ into a grading wire labelled $a$ as shown below.
\[\scalebox{0.8}{\input{transformfzero.tikz}}\]
We make this visually intuitive process formal by defining a graded presentation of $\Para{\G}{\catC}$.

The graded signature $\sig$ contains all the generators in $\sigC$ seen as $\eword$-graded generators, plus an additional $\col$-graded generator $\knightgenc : \eword \garrow{\col} k_{\col}$ for each $\col \in \sigv_0$, where $k_{\col} = \col \act \eword$. As explained before in \cref{lem:zerogradedpart}, there is a correspondence between monoidal $\sigC$-terms and $\eword$-graded $\sig$-terms because we did not add $\eword$-graded generators. Thus, all the axioms in $\axC$ can be interpreted as equations between graded terms, and we put those in $\axh$. These will be the only equations between $\eword$-graded terms so that $\eword$-graded string diagrams correspond to morphisms in $\catC$.

The remaining equations in $\axh$ will reflect the definition of "regrading" in $\Para{\G}{\catC}$ \cref{eqn:defn:parareg}. For each generator $g$ in the presentation of $\G$ (in $\sigv_1$), we add the equation on the left in \cref{eqn:slidinggenerator}, where $g$ is drawn vertically ($\inlinefig{0.5}{verticalgen}$) as it acts on grading wires and $g_{\eword}$ is a fixed $\sigC$-term in the equivalence class of the morphism $g \act \id_\eword \in \catC$, seen as an $\eword$-graded $\sig$-term. We also add the same equation for the symmetries $\inlinefig{0.6}{swapgradecol}$, which is drawn on the right in \cref{eqn:slidinggenerator}, using that $\swap{\col,\colb} \act \id_{\eword} = \swap{k_{\col},k_{\colb}}$.\begin{equation}\label{eqn:slidinggenerator}
  \scalebox{0.8}{\input{axslidinggen.tikz}} \qquad \qquad \scalebox{0.8}{\input{slidingswap.tikz}}.
\end{equation}
\begin{nota}
  Given $a \in \words{{\sigv_0}}$, we write $\inlinefig{0.5}{knightathick}$ for the term of type $\eword \garrow{a} k^a$ defined as the monoidal product $\bigotimes_{\col \in a} \knightgenc$ (we implicitly use the associativity of $\otimes$ \cref{eqn:th:assocpar}). A specific case of this being $a=\eword$, which gives $\inlinefig{0.5}{knighteword} = \id_{\eword}: \eword \garrow{\eword} \eword$ (represented by the empty diagram). We often omit the label $a$ and write $\knightgen$. Moreover, while we use thick wires to emphasise that the labels of the grading wire and the output wire are different, we will often omit this distinction and write $\knightbisgen$.
\end{nota}

Concisely, $\th \coloneq (\sig,\axh)$, where
\begin{equation}\label{defn:gthfromth}
  \begin{array}{rl}
    \sig & \coloneq \{ g : \word \garrow{\eword} \wordb \mid g : \word \rightarrow \wordb \in \sigC \} \sqcup \{ \knightgenc: \eword \garrow{\col} k_{\col} \mid \col \in \sigv_0 \}, \text{ and } \\[0.3em]
    \axh & \coloneq \{ f = f' \mid f = f' \in \axC \} \cup \{ \inlinefig{0.6}{axslidinggennolab} \mid g \in \sigv_1\} \cup \{\inlinefig{0.6}{slidingswapnolab}\},
  \end{array}
\end{equation}
where $\sigC$-terms are implicitly seen as $\eword$-graded $\sig$-terms in $\axh$ and $g_{\eword}$ denotes the chosen $\sigC$-term representing $g \act \id_{\eword}$ (\cref{lem:slidingterm} will show that this choice is immaterial).

It is important to note that, by construction and using \cref{lem:zerogradedpart} (which applies because $\eword$ is strictly terminal in $\G$; this assumption underlies the whole section), the "underlying prop" of $\sync{\th}$ is isomorphic to $\sync{\thC}$, and hence to $\catC$. Moreover, since the "underlying prop" of $\Para{\G}{\catC}$ is $\catC$, this establishes an isomorphism between the "underlying props" of $\sync{\th}$ and $\Para{\G}{\catC}$. It sees an $\eword$-graded string diagram as a diagram in $\sync{\thC}$, then sends it to the corresponding morphism in $\catC$ through the presentation, and finally to that same morphism seen as an $\eword$-graded morphism of $\Para{\G}{\catC}$. Our goal is now to make this an isomorphism between the whole "graded props".

Our argument relies on finding a particular shape that terms can be put into (akin to a normal form). Here, we want to factor any graded term as the composition of a very simple graded term and a $\sigC$-term, showing that any graded term is obtained by bending wires of an $\eword$-graded term. We need a technical lemma that generalises \cref{eqn:slidinggenerator} to all $\G$-morphisms.
\begin{lem}\label{lem:slidingterm}
  For any morphism $\inlinefig{0.5}{verticalterm}$ in $\G$, and any monoidal $\sigC$-term $t_{\eword}$ in the equivalence class of $t \act \id_{\eword} \in \catC$, the closure $\closed{\axh}$ contains the equation
  \begin{equation}\label{eqn:slidingterm}
    \scalebox{0.8}{\input{axslidingterm.tikz}}.
  \end{equation}
\end{lem}
\begin{proof}
  Since $\G$ is presented by a "monoidal theory" $(\sigv,\axv)$, we can proceed by induction on a representative monoidal $\sigv$-term of $\inlinefig{0.5}{verticalterm}$. When $t$ is a generator (resp.\@ symmetry) and $t_{\eword}$ is the representative chosen in \cref{defn:gthfromth}, \cref{eqn:slidingterm} is simply the axiom on the left (resp.\@ right) of \cref{eqn:slidinggenerator}, which belongs to $\axh$. When $t$ is $\id_{\col}$ (resp.~$\id_{\eword}$) and $t_{\eword} = \id_{k_{\col}}$ (resp.~$t_{\eword} = \id_{\eword}$), \cref{eqn:slidingterm} holds by \cref{eqn:th:trivialregrade,eqn:th:neutralid}. In both cases, if $u_{\eword}$ is any other $\sigC$-term in the same equivalence class as $t_{\eword}$, then $t_{\eword} = u_{\eword}$ is provable in $\thC$, hence it belongs to $\closed{\axh}$ by \cref{lem:zerogradedpart}, and \cref{eqn:slidingterm} holds for $u_{\eword}$ as well. This handles all base cases.

  For sequential composition, we first use the induction hypothesis (indicated by \textsf{IH}) twice in the following derivation.
  \[\scalebox{0.7}{\input{proofslidingtermseq.tikz}}\]
  We then show that $\inlinefig{0.7}{s0t0}$ is in the equivalence class of $(s\compsimple t) \act \id_{\eword}$. Indeed, by functoriality of $\act$, we have $(s\compsimple t) \act \id_{\eword} = (s\compsimple t) \act (\id_{\eword} \compsimple  \id_{\eword}) = (s \act \id_{\eword}) \compsimple  (t\act \id_{\eword})$. Thus, if $s_{\eword}$ and $t_{\eword}$ are representatives of $s \act \id_{\eword}$ and $t \act \id_{\eword}$ respectively, then $s_{\eword}\compsimple t_{\eword}$ is a representative of $(s\compsimple t) \act \id_{\eword}$. Finally, for any other term $u_{\eword}$ in that equivalence class, a proof of $s_{\eword}\compsimple t_{\eword} = u_{\eword}$ in the "monoidal theory" presenting $\catC$ can be adapted to a proof that $s_{\eword}\compsimple t_{\eword} = u_{\eword}$ belongs to $\closed{\axh}$.

  For parallel composition, we also apply the induction hypothesis twice.
  \[\scalebox{0.7}{\input{proofslidingtermpar.tikz}}\]
  Then, the rest of the argument follows as above (but instead of functoriality of $\act$, we use the fact that it is monoidal).
\end{proof}
We can now factor any term. One may visualise this process as sliding all the $\G$-morphisms through $\knightgen$, then pulling back the vertical dangling wires all the way to the left.
\begin{lem}\label{lem:decompositiongradedterm}
  Given a graded $\sig$-term $\inlinefig{0.5}{fnma}$, there is a monoidal $\sigC$-term $\inlinefig{0.5}{fzero}$ (seen as an $\eword$-graded $\sig$-term) such that
  \begin{equation}\label{eqn:decompositionfzero}
    \scalebox{0.8}{\input{decompositionfzero.tikz}}.
  \end{equation}
\end{lem}
\begin{proof}
  We construct $f_{\eword}$ explicitly by induction on $f$. The base cases are trivial:
  \begin{itemize}
    \item if $f$ is a generator of grade $\eword$ or an identity or a symmetry, then $a=\eword$ and the $\eword$-product of $\knightgen$ is $\id_{\eword}$, so $f_{\eword} = f$ makes \cref{eqn:decompositionfzero} hold, and
    \item if $f$ is $\knightgenc$, then we can set $f_{\eword} = \begin{tikzpicture}
	\begin{pgfonlayer}{nodelayer}
		\node [style=none] (0) at (-1, 0) {};
		\node [style=none] (1) at (0.5, 0) {};
		\node [style=none] (2) at (-1.5, 0) {$k_\col$};
		\node [style=none] (3) at (1, 0) {$k_\col$};
	\end{pgfonlayer}
	\begin{pgfonlayer}{edgelayer}
		\draw [style=single] (0.center) to (1.center);
	\end{pgfonlayer}
\end{tikzpicture}
$.
  \end{itemize}
  For the inductive steps, suppose we have the following terms.
  \[\scalebox{0.65}{\input{inductionnormalform.tikz}}\]
  In the derivations below, the dotted horizontal lines are identities on $\eword$ (drawn as in \cref{defn:monterms}) that make explicit which subterms are composed in sequence or in parallel at each step. The dotted boxes in the last diagrams of each derivation identify the shape of \cref{eqn:decompositionfzero}.

  For sequential composition, we have the following derivation identifying the term $(f\comp g)_{\eword}$ as $(f_{\eword} \tensor \id_{k^b})\compsimple g_{\eword}$.
  \[\scalebox{0.85}{\input{fzeroconseq.tikz}}\]

  For parallel composition, we have the following derivation identifying the term $(f \tensor f')_{\eword}$ as $(\id_{\word} \otimes \swap{\word',k^a} \otimes \id_{k^b})\compsimple(f_{\eword} \otimes f'_{\eword})$.
  \[\scalebox{0.85}{\input{derivationnormalformpar.tikz}}\]

  For "regrading", we rely on \cref{lem:slidingterm}. Given a vertical term $t$, the following derivation identifies $(t \regrade f)_{\eword}$ as $(\id_{\word} \otimes t_{\eword})\compsimple f_{\eword}$.
  \[\scalebox{0.85}{\input{derivationnormalformreg.tikz}}\qedhere\]
\end{proof}

\begin{thm}\label{thm:gradedpresfrompres}
  Let $\catC$ be a "prop" presented by $\thC$ and $\act: \G \times \catC \rightarrow \catC$ a symmetric monoidal action. 
  The "theory@@GRD" $\th = (\sig,\axh)$ defined in \cref{defn:gthfromth} is a graded presentation of $\Para{\G}{\catC}$.
\end{thm}
\begin{proof}
  We need to construct an isomorphism $\sem{-}: \sync{\th} \rightarrow \Para{\G}{\catC}$.

  We will use \cref{lem:soundness} to recursively define $\sem{-}$ by giving its value on generators. Following the discussion above, $\sem{-}$ sends each $\eword$-graded generator to the corresponding $\eword$-graded morphism in $\Para{\G}{\catC}$ via
  \({\sync{\th}}_{\eword}(\word,\wordb) \cong \sync{\thC}(\word,\wordb) \cong \catC(\word,\wordb) \cong \Para{\G}{\catC}_{\eword}(\word,\wordb)\).
  For the remaining generators $\knightgenc: \eword \garrow{\col} k_{\col}$, we define $\sem{\knightgenc}$ to be the morphism resulting from the action of $\id_{\col}$ in $\G$ on $\id_{\eword}$ in $\catC$, namely, $\sem{\knightgenc} \coloneq \id_{\col} \act \id_{\eword}: \eword \garrow{\col} k_{\col}$.

  The assignment $\sem{-}$ is canonically extended to graded $\sig$-terms, and it is easy to see that it acts exactly like the isomorphism $\sync{\thC} \rightarrow \catC$ on $\eword$-graded terms. Therefore, for all axioms $f=g \in \axC$, we have $\sem{f} = \sem{g}$. For the remaining axioms in $\axh$, those coming from \cref{eqn:slidinggenerator}, we have for any morphism $t: b \rightarrow a \in \G$,
  \begin{align*}
    \sem{\!\inlinefig{0.6}{tregradeknightbis}}
     & = t \regrade \sem{\!\inlinefig{0.5}{knighta}}                  &  & \text{recursive def.~of $\sem{-}$}                                                   \\
     & = t \regrade (\id_a \act \id_{\eword})                         &  & \text{def.~of $\sem{\!\inlinefig{0.5}{knighta}}$}                                    \\
     & = (t \act \id_{\eword}) \compsimple  (\id_a \act \id_{\eword}) &  & \text{def.~of $\regrade$ in $\Para{\G}{\catC}$ \cref{eqn:defn:parareg}}              \\
     & = t \act \id_{\eword}                                          &  & \text{functoriality of $\act$}                                                       \\
     & = (\id_b \act \id_{\eword}) \compsimple  (t \act \id_{\eword}) &  & \text{functoriality of $\act$}                                                       \\
     & = (\id_b \act \id_{\eword}) \comp (t \act \id_{\eword})        &  & \text{def.~of $\comp$ in $\Para{\G}{\catC}$ \cref{eqn:defn:paracomp} and strictness} \\
     & = \sem{\!\inlinefig{0.5}{knightb}} \comp \sem{t_{\eword}}      &  & \text{def.~of $\sem{\!\inlinefig{0.5}{knightb}}$ and $t_{\eword}$}                   \\
     & = \sem{\!\inlinefig{0.6}{knightt0bis}}                         &  & \text{recursive def.~of $\sem{-}$}
  \end{align*}
  We have shown that all axioms in $\axh$ hold in $\Para{\G}{\catC}$, thus we have defined a model of $\th$, and by \cref{lem:soundness}, $\sem{-}$ factors through an identity-on-objects functor $\sem{-}: \sync{\th} \rightarrow \Para{\G}{\catC}$.

  To show that $\sem{-}$ is in fact an isomorphism, it suffices to prove it is fully faithful.  This relies on the fact that $\sem{-}$ is an isomorphism (hence fully faithful) when restricted to $\eword$-graded terms, and the decomposition in \cref{lem:decompositiongradedterm}.

  For faithfulness, suppose that $f,g: \word \garrow{a} \wordb$ are $\sig$-terms such that $\sem{f} = \sem{g}$. We will show that $f = g \in \closed{\axh}$. Observe that the proofs of \cref{lem:slidingterm,lem:decompositiongradedterm} did not rely on the equations in $\axh$, only on equations in \cref{fig:axiomssyncat} that are valid in any "graded SMC". Therefore, we automatically have \[\sem{\!\inlinefig{0.5}{idtensknight}} \comp \sem{f_{\eword}} = \sem{f} = \sem{g} = \sem{\!\inlinefig{0.5}{idtensknight}} \comp \sem{g_{\eword}}.\]
  Let us unroll some definitions to study precomposition by $\sem{\!\inlinefig{0.5}{idtensknight}}$:
  \begin{align*}
    \sem{\!\inlinefig{0.5}{idtensknight}} \comp - & = (\id_{\word} \otimes (\id_a \act \id_{\eword})) \comp -                           &  & \text{def.~of $\sem{-}$}                                              \\
                                                  & = (\id_{\eword} \act (\id_{\word} \otimes (\id_a \act \id_{\eword})))\compsimple  - &  & \text{def.~of $\comp$ in $\Para{\G}{\catC}$ \cref{eqn:defn:paracomp}} \\
                                                  & = \id_{\word\otimes k^a}\compsimple  -                                              &  & \text{functoriality of $\act$ and $\otimes$}
  \end{align*}
  This shows that precomposition by $\sem{\!\inlinefig{0.5}{idtensknight}}$ is injective, and hence, with the previous equation, we obtain $\sem{f_{\eword}} = \sem{g_{\eword}}$. Since $\sem{-}$ is faithful on $\eword$-graded terms, we conclude that $f_{\eword} = g_{\eword} \in \closed{\axh}$, thus $f = g \in \closed{\axh}$.

  For fullness, let $o: \word \garrow{a} \wordb$ be a morphism in $\Para{\G}{\catC}$. We will show there is a term $f_o$ satisfying $\sem{f_o} = o$. Viewing it inside $\catC$, $o$ is a morphism $\word \otimes k^a \rightarrow \wordb$, which we can see as an $\eword$-graded morphism $o_{\eword}: \word \otimes k^a \garrow{\eword} \wordb$ in $\Para{\G}{\catC}$. Since $\sem{-}$ is full on $\eword$-graded morphisms, there is an $\eword$-graded $\sig$-term $f_{o_{\eword}}$ such that $\sem{f_{o_{\eword}}} = o_{\eword}$. Now, we define $f_o$ by precomposing $f_{o_{\eword}}$ with $\inlinefig{0.5}{idtensknight}$, and we find (as desired)
  \[\sem{f_o} = \sem{\!\inlinefig{0.5}{idtensknight}} \comp \sem{f_{o_{\eword}}} = \id_{\word \otimes k^a} \compsimple  o = o.\qedhere\]
\end{proof}

\section{Case Study: Imprecise Probability}\label{sec:application}
In this section, we use the previously introduced framework (in particular, \cref{thm:gradedpresfrompres}) to completely axiomatise a "graded category" $\BImP$ of imprecise probabilistic processes. $\BImP$ may be regarded as a variant of the category $\ImP$ studied in~\cite{LiellCock2025}, where we restrict to powers of the two-element set as objects. The categories $\G$ and $\catC$ will be instantiated respectively by $\op{\FInj}$, the opposite category of finite sets and injections, and  $\BStoch$, the category of stochastic matrices whose dimensions are powers of two. We will see $\op{\FInj}$ as a subcategory of $\BStoch$ acting on it with the monoidal product. The result of applying the "para construction" is $\BImP$, a category of nondeterministic stochastic matrices. Both $\op{\FInj}$ and $\BStoch$ have already been presented diagrammatically in the literature, so \cref{thm:gradedpresfrompres} readily applies to find a graded presentation of $\BImP$.

\subsection{\texorpdfstring{$\BStoch$}{BStoch}: Binary Stochastic Matrices}\label{sec:finstoch}
We first recall the definition of $\BStoch$ and its presentation in~\cite{Piedeleu2025b}.

\begin{defi}[$\BStoch$]
  \AP An $m\times n$ \textbf{stochastic matrix} is a matrix with $m$ rows and $n$ columns with entries in the interval $\intro*\probs \coloneq [0,1]$ such that, in each column, the sum of all the entries is $1$.
  The "SMC" $\FStoch$ has objects the natural numbers and as morphisms $n \rightarrow m$ the $m\times n$ stochastic matrices, with sequential composition given by matrix multiplication. On objects, $n \otimes m$ is the usual product of integers. On matrices, given $A: n \rightarrow  m$ and $A' : n' \rightarrow  m'$, $A \otimes A' : nn' \rightarrow mm'$ is defined in block matrix form as the Kronecker product
  \[A \otimes A' \coloneq \scalebox{0.8}{$\begin{bmatrix} a_{11}A' & \dots  & a_{1n}A' \\
                \vdots   & \ddots & \vdots   \\
                a_{m1}A' & \dots  & a_{mn}A'\end{bmatrix}$}.\]

  The category $\intro*\BStoch$ has objects the natural numbers, and as morphisms $n \rightarrow m$ the $2^m \times 2^n$ stochastic matrices: $\BStoch(n,m) = \FStoch(2^n,2^m)$.
  All the structure of $\BStoch$ is obtained by seeing it as a full subcategory of $\FStoch$ with the inclusion $\BStoch \hookrightarrow \FStoch$ sending $n$ to $2^n$. In particular, the monoidal product on objects is addition since $2^n \otimes 2^{m} = 2^n2^m = 2^{n+m}$, thus $\BStoch$ is a "prop" with a single colour.\footnote{We identify the set of finite words over a single colour with the set of natural numbers by taking the length of a word, and concatenation becomes addition of lengths.} The object $0$ is the unit of the monoidal product as $0+n = n$.
\end{defi}

In short, one may regard $\BStoch$ as a binary version of $\FStoch$. It is pertinent in our developments because, unlike $\FStoch$, a "monoidal theory" presenting $\BStoch$ is known.\footnote{To be precise, a presentation of $\FStoch$ is known when the monoidal product is given by direct sum of matrices~\cite{Fritz09}, but not with the Kronecker product. This difference is highly relevant in applications, because parallel probabilistic processes may be correlated only with the latter.} We recall the "theory@@MON" from~\cite{Piedeleu2025b}, as we will leverage it for our axiomatisation. 

\begin{defi}\label{defn:ThCC}
  \AP The "monoidal theory" $\intro*\thCC$ has generators $\inlinefig{0.9}{deletegen}$, $\inlinefig{0.9}{copygen}$, $\inlinefig{0.87}{andgen}$, $\inlinefig{0.87}{notgen}$, and $\inlinefig{0.9}{statep}$ for all $p \in \probs$, and its axioms are listed in~\cref{fig:binstocheqns} (reproducing~\cite[Figure~4]{Piedeleu2025b}).
  \AP Following~\cite{Piedeleu2025b}, we write $\intro*\Circ$ for the "syntactic category@@MON" $\sync{\thCC}$.
\end{defi}

Before stating the presentation result, it is worth providing some intuition on how string diagrams of $\Circ$ are interpreted in $\BStoch$, through an assignment $\sem{-}$ of stochastic matrices to the generators of $\thCC$. We may regard $\sem{-}$ as a semantics, giving an operational meaning to string diagrams of $\Circ$ as circuits.

First, a single wire $\inlinefig{0.7}{lonewire}$ (the identity $1 \to 1$) carries \textbf{probabilistic bits}, that is, either the bit $\ket{\b1}$ which also represents the column vector $\raisebox{0.3ex}{\scalebox{0.55}{$\begin{bmatrix} 1\\0 \end{bmatrix}$}}$, or the bit $\ket{\b0}=\raisebox{0.3ex}{\scalebox{0.55}{$\begin{bmatrix} 0\\1 \end{bmatrix}$}}$, or a convex combination $p \ket{\b1} + (1-p)\ket{\b0} = \raisebox{0.3ex}{\scalebox{0.55}{$\begin{bmatrix} p\\1-p \end{bmatrix}$}}$ with $p \in \probs$. We interpret a probabilistic bit $p \ket{\b1} + (1-p)\ket{\b0}$ as a probability distribution $\dist$ on $\intro*\B \coloneq \{\b1,\b0\}$ with $\dist(\b1) = p$ and  $\dist(\b0) = 1-p$.
The semantics of the identity wire is the identity $2 \times 2$ matrix: it leaves the bit unchanged.

A double wire $\inlinefig{0.7}{doublewire}$ carries two probabilistic bits. However, probabilistic bits can be \emph{correlated}, not unlike entangled qubits in quantum circuits. Formally, a double wire carries distributions on $\B\otimes \B = \{\b1\b1,\b1\b0,\b0\b1,\b0\b0\}$. If $\dist$ and $\distb$ are two probabilistic bits, then their Kronecker product is a distribution on $\B \otimes \B$---we identify $\ket{\b1\b1}$ with the Kronecker product $\ket{\b1}\ket{\b1}=\ket{\b1}\otimes \ket{\b1}$, and similarly for $\ket{\b1\b0}$, $\ket{\b0\b1}$, and $\ket{\b0\b0}$. However, a distribution on $\B \otimes \B$ does not always arise as the independent product of two probabilistic bits, e.g.\@ $\sfrac{1}{2}\ket{\b1\b1}+\sfrac{1}{2}\ket{\b0\b0}$.

Any morphism in $\BStoch(0,n)$ is a \textbf{state} that `emits' $n$ probabilistic bits, and it corresponds to a column vector of dimension $2^n$ whose entries sum up to $1$---the entry in the $i$th row will be the weight of the distribution at $\ket{\mathtt{bin}(2^n-i)}$, where $\mathtt{bin}(2^n-i)$ is the zero-padded binary representation of $2^n-i$. For $n=1$, and given $p \in \probs$, the generator $\inlinefig{0.9}{statep}$ emits $p\ket{\b1}+(1-p)\ket{\b0}= \raisebox{0.3ex}{\scalebox{0.55}{$\begin{bmatrix} p\\1-p \end{bmatrix}$}}$. In particular, $\inlinefig{0.9}{stateone}$ emits $\ket{\b1}$ and $\inlinefig{0.9}{statezero}$ emits $\ket{\b0}$. The generator $\inlinefig{0.9}{deletegen}$ behaves as \textbf{discard}: it accepts any bit coming in and outputs nothing. The generator $\inlinefig{0.9}{copygen}$ behaves as \textbf{copy}: it takes an input $p\ket{\b1}+(1-p)\ket{\b0}$ and creates two correlated bits $p\ket{\b1\b1}+(1-p)\ket{\b0\b0}$. The remaining generators $\inlinefig{0.9}{andgen}$ and $\inlinefig{0.9}{notgen}$ are \textbf{and} and \textbf{not}: they behave like their classical circuits counterpart sending $\ket{x}\ket{y}$ to $\ket{x\ \mathtt{and}\ y}$ and $\ket{x}$ to $\ket{\mathtt{not}\ x}$ respectively.\footnote{A stochastic matrix $f \in \BStoch(n,m)$ must preserve convex combinations, thus its value on all $n$-dimensional products of $\ket{\b1}$ and $\ket{\b0}$, equivalently on $\ket{\vvx}$ for all $\vvx\in \B^n$, completely determines it.}

The symmetries, represented diagrammatically by crossing wires, swap the bits carried in each wires. For instance, $\swap{1,1} = \inlinefig{0.6}{swap}$ sends $\ket{xy}$ to $\ket{yx}$ for any $x,y \in \B$.
Here are a few symmetries with the corresponding permutation matrix (dotted wires represent the identity on $2^0= 1$ and they carry no probabilistic bit).
\[\inlinefig{0.65}{idswapid} = \scalebox{0.65}{$\begin{bmatrix}
        1
      \end{bmatrix}$} \quad
  \inlinefig{0.65}{idswap} = \scalebox{0.65}{$\begin{bmatrix}
        1 & 0 \\ 0 & 1
      \end{bmatrix}$}\quad
  \inlinefig{0.65}{swapid} = \scalebox{0.65}{$\begin{bmatrix}
        1 & 0 \\ 0 & 1
      \end{bmatrix}$} \quad
  \inlinefig{0.65}{swap} = \scalebox{0.65}{$\begin{bmatrix}
        1 & 0 & 0 & 0 \\
        0 & 0 & 1 & 0 \\
        0 & 1 & 0 & 0 \\
        0 & 0 & 0 & 1
      \end{bmatrix}$}.
\]
\begin{thmC}[{\cite[Theorem 3.13]{Piedeleu2025b}}]\label{thm:axbstoch}
  $\thCC$ presents $\BStoch$.
\end{thmC}
More precisely, the assignment $\sem{-}$ described above extends to an isomorphism of "props" $\sem{-}: \Circ \cong \BStoch$.

\subsection{\texorpdfstring{$\BImP$}{BImP}: Binary Imprecise Probability}\label{sec:bimp}

The recent work~\cite{LiellCock2025} studies a "graded category" $\ImP$ as a model for probabilistic programming languages based on imprecise probability. Formally, $\ImP$ is defined as the "graded category" resulting from the "para construction" $\Para{\FStochSurj}{\FStoch}$ (recall \cref{exmp:gradedpara}). That is, $\ImP$-morphisms are stochastic matrices with grading from the $\FStoch$-subcategory of surjective stochastic matrices. Intuitively, the grading indicates nondeterministic choices of a probabilistic process.

Our goal is to identify the "monoidal theory" presenting a variant of $\ImP$. The reason for not tackling $\ImP$ directly is that, as mentioned in \cref{sec:finstoch}, a presentation of $\FStoch$ is not known. However, as we will see in~\cref{sec:impprog}, our variant suffices to provide a semantic model for imprecise probabilistic programming. We tweak $\ImP$ along two axes.

First, for the base category, we pick $\BStoch$, the binary version of $\FStoch$, which enjoys a monoidal presentation by $\Circ$ (\cref{thm:axbstoch}). Second, we also restrict the grading category. As remarked in~\cite[p.\@ 17]{LiellCock2025}, $\FStochSurj$ is not the only possible choice of grading category: for their developments, it suffices to consider only the deterministic surjections $2^B \rightarrow 2^A$ that are images of injections $A \rightarrow B$ under the contravariant powerset functor $2^{-}: \Set \rightarrow \Set$.
\AP The "prop" corresponding to this class of morphisms is $\op{\intro*\FInj}$, the opposite of the skeleton category of finite sets and injections.
In other words, $\op{\FInj}$ has objects the natural numbers and a morphism $n \rightarrow m$ is an injection $\fset{m} \rightarrow \fset{n}$, where $\fset{n}\coloneq \{0,\dots, n-1\}$. Note that $0$ is strictly terminal in $\op{\FInj}$, as required in \cref{sec:modularcompleteness}: for every $n$, the empty function is the unique injection $\fset{0} \rightarrow \fset{n}$, hence the unique morphism $n \rightarrow 0$, and there is no injection $\fset{n} \rightarrow \fset{0}$ for $n > 0$, hence no morphism $0 \rightarrow n$.

Just as $\FStochSurj$ is a subcategory of $\FStoch$, we can view $\op{\FInj}$ as a subprop of $\BStoch$ with the functor $I: \op{\FInj} \hookrightarrow \BStoch$ defined as follows. It is an identity-on-objects functor $I(n) = n$, and it sends an injective function $t: \fset{m} \rightarrow \fset{n}$ (a morphism $n \rightarrow m$ in $\op{\FInj}$) to the stochastic matrix $I(t)$ that sends $\ket{\ttx_0\cdots \ttx_{n-1}}$ to $\ket{\ttx_{\scriptscriptstyle t(0)}\cdots \ttx_{\scriptscriptstyle t(m-1)}}$. In particular, if $t: \fset{2} \rightarrow \fset{2}$ is the function swapping $0$ and $1$, then $I(t) = \swap{1,1} = \inlinefig{0.6}{swap}$.

We may now obtain a symmetric monoidal action $\act : \op{\FInj} \times \BStoch \rightarrow \BStoch$ defined on objects by $(n,m) \mapsto n+m$, and on morphisms by $(t,A) \mapsto I(t) \otimes A$. This is an example of the self-action of a subcategory in \cref{exmp:selfaction} since we view $\op{\FInj}$ as a subcategory of $\BStoch$ via $I$.

Applying the "para construction" (\cref{exmp:gradedpara}), we obtain $\Para{\op{\FInj}}{\BStoch}$, which we call $\BImP$. It is a "graded prop" with a single colour. Just as $\BStoch$ compared to $\FStoch$, we regard $\BImP$ as a binary version of $\ImP$.
\begin{defi}[$\BImP$]
  \AP We unroll the definition of $\Para{\op{\FInj}}{\BStoch}$ to make the structure of $\intro*\BImP$ explicit.
  \begin{itemize}
    \itemAP Objects of $\BImP$ are natural numbers (elements of $\intro*\N$).
    \item For any $n,m,a \in \N$, the set of morphisms $n \rightarrow m$ with grade $a$ is $\BImP_a(n,m) \coloneq \BStoch(a \otimes n,m)= \FStoch(2^{a+n},2^m)$.
    \item For any injective function $t: \fset{a} \rightarrow \fset{b}$ (seen as a morphism $b \rightarrow a$ in $\op{\FInj}$), the "regrading" $t \regrade -: \BImP_{a}(n, m) \rightarrow \BImP_{b}(n, m)$ is given by
          \[t\regrade g \coloneq b \otimes n \garrow{I(t) \otimes \id_{n}} a \otimes n \garrow{g} m.\]
    \item For any two morphisms $f \in \BImP_a(n,m)$ and $g \in \BImP_b(m, \ell)$, their composition is
          \[f \comp g \coloneq a \otimes b \otimes n \garrow{\swap{a,b}\otimes \id_n} b \otimes a \otimes n \garrow{\id_{b} \otimes f} b \otimes m \garrow{g} \ell.\]
    \item For any two morphisms $f \in \BImP_a({n},{m})$ and $f' \in \BImP_b({n'},{m'})$, their tensor is 
          \[f\tensor f' \coloneq a \otimes b \otimes {n} \otimes {n'} \garrow{\id_{a}\otimes \swap{b,n} \otimes \id_{{n'}}}  a \otimes {n} \otimes b \otimes {n'} \garrow{f \otimes f'} {m}\otimes {m'}.\]
  \end{itemize}
\end{defi}

Now, how is $\BImP$ interpreted as a category of imprecise probabilistic processes? A morphism $f \in \BImP_a(n,m)$ corresponds to a $2^m \times 2^a2^n$ stochastic matrix. We may regard $f$ as decomposed into $2^a$ stochastic \textbf{submatrices} $f\ket{\vvx}\ket{-}$ of dimension $2^m\times 2^n$, with $\vvx$ ranging in $\B^a$ and $f\ket{\vvx}\ket{-}$ denoting the morphism in $\BStoch$ that sends $\ket{\vvu}$ to $f\ket{\vvx}\ket{\vvu}$.
To unravel $f$, we write all the submatrices in decreasing order of the number represented by $\vvx$ and separate them with vertical lines. For example, \renewcommand{\arraystretch}{0.8}
\[f = \left[\!\! \begin{array}{c|c}
      f\ket{\tt1}\ket{-} & f\ket{\tt0}\ket{-}
    \end{array}\!\! \right]  = \scalebox{0.7}{$\Big[\!\! \begin{array}{c|c}
          1 & 0.5 \\
          0 & 0.5
        \end{array}\!\! \Big]$}\]
represents a morphism $f: 0 \garrow{1}  1$ in $\BImP$. It is decomposed as $2^1 = 2$ stochastic matrices of dimension $2^1 \times 2^0$ separated by a vertical line. Similarly to the interpretation of $\ImP$ given in~\cite{LiellCock2025}, we view the separations as nondeterministic choices between the matrices on either side. In the example of $f = \scalebox{0.7}{$\Big[\!\! \begin{array}{c|c}
          1 & 0.5 \\
          0 & 0.5
        \end{array}\!\! \Big]$}$, each choice is a column vector, thus a distribution. Compared with the algebraic syntax used in e.g.\@ \cite{Bonchi22}, $f$ may also be written as $\dirac{1} \oplus (\dirac{1} +_{0.5} \dirac{0})$.\footnote{$\dirac{i}$ denotes the Dirac distribution on $\ket{i}$. There is a caveat because $\oplus$ models a commutative nondeterministic choice in~\cite{Bonchi22}, while in $\BImP$,\renewcommand{\arraystretch}{0.8} $\scalebox{0.7}{$\Big[ \!\!\begin{array}{c|c}
            1 & 0.5 \\
            0 & 0.5
          \end{array}\!\! \Big]$}  \neq \scalebox{0.7}{$\Big[ \!\!\begin{array}{c|c}
            0.5 & 1 \\
            0.5 & 0
          \end{array}\!\! \Big]$}$.
  This is further discussed in~\cite{LiellCock2025}.}

\subsection{Graded Presentation of \texorpdfstring{$\BImP$}{BImP}}
We now have all the ingredients to give an axiomatisation of $\BImP$. First, it is well known that the "prop" $\op{\FInj}$ is presented by a "monoidal theory" containing a single generator $\inlinefig{0.75}{delgradegenerator}$ and no axioms, see e.g.\@ \cite[Example 2.28(a)]{zanasi:tel-01218015}. As discussed, the "prop" $\BStoch$ is presented by the "monoidal theory" $\thCC$ (\cref{thm:axbstoch}). The "graded prop" $\BImP$ is defined as $\Para{\op{\FInj}}{\BStoch}$. Therefore, by \cref{thm:gradedpresfrompres}, we modularly obtain a presentation of $\BImP$ by treating $\thCC$ as a $0$-"graded theory" and adding a generator $\knightbisgen$, together with the equations \cref{eqn:slidinggenerator} for the generator $g = \inlinefig{0.75}{delgradegenerator}$ and the symmetry $g= \inlinefig{0.75}{swapgradegenerator}$. There are natural choices for the terms representing $g_{\eword}$ in both cases, the discard generator in ${\thCC}$ and the swapping of wires respectively. Thus the additional axioms are (no need for thick wires since $1 \act 0 = 1$ in this instance)
\begin{equation}\label{eqn:slidingdel}
  \scalebox{0.8}{\input{slidingdelbis.tikz}} \quad \text{ and } \quad \scalebox{0.8}{\input{slidingswapbis.tikz}}.
\end{equation}

\begin{thm}\label{thm:axbimp}
  \AP The "graded theory" obtained by \cref{thm:gradedpresfrompres} for $\Para{\op{\FInj}}{\BStoch}$ described above, denote it $\intro*\thCCi$, is a graded presentation of $\BImP$.
\end{thm}
We write $\GCirc$ (imprecise $\Circ$) for the "syntactic category@@GRD" of $\thCCi$. We can read string diagrams of $\GCirc$ with the same interpretation as those of $\Circ$, except we view the vertical wires as carrying nondeterministic bits carrying one of $\ket{\b1}$ or $\ket{\b0}$.

\begin{rem}\label{rem:slidingfrombelow}
  We can now reveal how the choice of putting the grading wire on the bottom makes more visual sense. In \cref{lem:slidingterm}, we prove a generalisation of the axiom \cref{eqn:slidinggenerator} for all $\G$-morphisms. In the "theory@@GRD" presenting $\BImP$, this result becomes visually intuitive. Any monoidal term representing a morphism in $\op{\FInj}$ can slide through the generator $\knightbisgen$ (or $a$-wise monoidal products of it) after a doing a $90$-degree rotation as drawn below. If the grading wires were coming from the top, this generalised axiom would involve a reflection and a rotation. This is simply an artefact of the choice of reading convention.  \[\scalebox{0.8}{\input{slideexmp2.tikz}}\]
\end{rem}

\section{Imprecise Probabilistic Programming}\label{sec:impprog}

The primary goal of~\cite{LiellCock2025} is to provide theoretical foundations for a compositional language modelling imprecise probability. They prove that their framework of choice, the category $\ImP$, is a convenient setting to interpret such a language. In this section, we interpret in $\BImP$ a simpler language in the spirit of theirs. Our aim is not to recover the language of~\cite{LiellCock2025} in full (the differences are discussed in \cref{sec:progbimp}), but rather to illustrate what the string diagrammatic presentation of \cref{thm:axbimp} brings to the table: programs become graded string diagrams, and properties of the language can be established by equational reasoning on these diagrams.

In \cref{sec:progbimp}, we define a programming language with nondeterministic and probabilistic choices that can be translated into string diagrams of $\BImP$. We also show that this language is commutative and affine, i.e.\@ the equations in \cref{thm:desiderata} hold, hence it satisfies \textit{Desideratum 1} of~\cite{LiellCock2025}. In \cref{sec:conditioning}, we extend the language with a construct for exact conditioning by `augmenting' $\BImP$, mirroring how $\Circ$ is augmented to $\mathsf{ProbCirc}$ in~\cite{Piedeleu2025b}. This is immediate thanks to the modular completeness result \cref{thm:gradedpresfrompres}. We end the section by presenting two examples of graphical programs that illustrate the potential for modelling imprecise probabilities.

\subsection{Programming in \texorpdfstring{$\BImP$}{BImP}}\label{sec:progbimp}
We extend the language studied in~\cite{Piedeleu2025b} (without $\observe{\!\!}$ for now) with the $\knight$ construct for nondeterministic choices discussed in~\cite{LiellCock2025}. We inductively define the syntax where $x$ ranges over a countable set of variables and $p$ over $\probs=[0,1]$:\begin{align*}
  e, e_1, e_2          & \Coloneq x \ \mid\ \bern{p}\ \mid\ \knight\ \mid\ \pair{e_1}{e_2} \ \mid \ \fst{e}\ \mid \ \snd{e}\                           \\*
                       & \quad \dots\ \mid \ifte{e}{e_1}{e_2}\ \mid \ \letin{x}{e_1}{e_2}                                                    &  & \text{(expressions)} \\*
  \tau, \tau_1, \tau_2 & \Coloneq \B \ \mid \ \tau_1 \otimes \tau_2                                                               &  & \text{(types)}       \\*
  \Gamma               & \Coloneq \emptyset\  \mid \ \Gamma,  x : \tau                                                           &  & \text{(contexts)},
\end{align*}
A context is \textbf{well-formed} if no variable occurs twice in it; we only consider well-formed contexts, and we regard them as lists rather than sets, because the position of a variable in the context determines the wire it is assigned to in \cref{fig:semantics}. Weakening and exchange of variables are handled by the explicit rules \RuleWeak{} and \RuleEx{} in \cref{fig:typerules}. The intended semantics for this language is rather natural: $\bern{p}$ generates a probabilistic bit $p\ket{\b1}+(1-p)\ket{\b0}$, $\knight$ generates a nondeterministic bit that is $\ket{\b1}$ or $\ket{\b0}$, the pairing, projections, and $\ifte{\!\!}{\!\!}{\!\!}$ have their usual semantics, and similarly for the $\mathtt{let}$ primitive. That last primitive is essential to allow reusing the result of a process (\textit{cf.}\@ named choices in~\cite{LiellCock2025}).

\AP The type and grade of an expression can be inferred with the rules in \cref{fig:typerules} (\textit{cf.}\@ \cite[Figure~2]{Piedeleu2025b} and~\cite[Section 2.3]{LiellCock2025}), where $\Gamma \vdash e: \tau \intro*\with a$ should be read as `$e$ has type $\tau$ and grade $a$ in the context $\Gamma$'. We note three differences from~\cite{LiellCock2025}.
\begin{enumerate}
  \item The grades of the two branches of an $\ifte{\!\!}{\!\!}{\!\!}$ do not need to coincide, and they add up in the grade of the resulting expression. This stems from the fact that we interpret conditionals with a morphism inside $\BImP$---the multiplexer defined in~\cite{Piedeleu2025b}---which is composed with the morphisms representing the condition and both branches. Thus, the grades of $e$, $e_1$ and $e_2$ are all added in $\ifte{e}{e_1}{e_2}$. In contrast, \cite{LiellCock2025} treat conditional statements with a coproduct (that distributes with the monoidal product) which combines two morphisms of grade $a$ into one morphism of grade $a$.
  \item We do not have a construct for "regrading" inside the language. In the context of $\BImP$, "regrading" would essentially allow one to increase the grade of an expression:
        \[\Gamma \vdash e: \tau \with a \implies \Gamma \vdash e: \tau \with b\qquad a \leq b.\]
        This is needed in~\cite{LiellCock2025} to increase the grade of a branch in a conditional to match the grade of the other branch.
  \item We do not have \emph{named} nondeterministic choices. In~\cite{LiellCock2025}, nondeterministic choices carry names, and different occurrences of the same name denote the same choice; this allows correlating choices made in different parts of a program, for instance in the two branches of a conditional. In our language, every occurrence of $\knight$ is a fresh, independent choice, and the only way to reuse its outcome is to bind it with $\mathtt{let}$.
\end{enumerate}
As a consequence, our language is less expressive than that of~\cite{LiellCock2025}. \cref{rem:nohoisting} explains the formal reason behind these differences: $\BImP$ lacks the coproducts that~\cite{LiellCock2025} use to interpret conditional statements.
\begin{figure}
  \begin{gather*}
    \begin{bprooftree}
      \AxiomC{$x \notin \Gamma$}
      \RightLabel{\hypertarget{rulevar}{\scriptsize \textsc{Var}}}
      \UnaryInfC{$\Gamma, x: \tau \vdash x : \tau \with 0$}
    \end{bprooftree}
    \begin{bprooftree}
      \AxiomC{$x \notin \Gamma$}
      \AxiomC{$\Gamma \vdash e : \tau_1 \with a$}
      \RightLabel{\hypertarget{ruleweak}{\scriptsize \textsc{Weak}}}
      \BinaryInfC{$\Gamma, x: \tau_0 \vdash e : \tau_1 \with a$}
    \end{bprooftree}\begin{bprooftree}
      \AxiomC{$\Gamma \vdash e : \tau_1 \tensor \tau_2\with a$}
      \RightLabel{\hypertarget{rulefst}{\scriptsize \textsc{Fst}}}
      \UnaryInfC{$\Gamma \vdash \fst{e}: \tau_1\with a$}
    \end{bprooftree}\\
    \begin{bprooftree}
      \AxiomC{$p \in \probs$}
      \RightLabel{\hypertarget{rulebern}{\scriptsize \textsc{Flip}}}
      \UnaryInfC{$\Gamma \vdash \bern{p} : \B\with 0$}
    \end{bprooftree}
    \begin{bprooftree}
      \AxiomC{$\phantom{[]}$}
      \RightLabel{\hypertarget{ruleknight}{\scriptsize \textsc{Knight}}}
      \UnaryInfC{$\Gamma \vdash \knight : \B\with 1$}
    \end{bprooftree}\begin{bprooftree}
      \AxiomC{$\Gamma \vdash e : \tau_1 \tensor \tau_2\with a$}
      \RightLabel{\hypertarget{rulesnd}{\scriptsize \textsc{Snd}}}
      \UnaryInfC{$\Gamma \vdash \snd{e}: \tau_2\with a$}
    \end{bprooftree}
    \\
    \begin{bprooftree}
      \AxiomC{$\Gamma \vdash e_1 : \tau_1\with a$}
      \AxiomC{$\Gamma \vdash e_2 : \tau_2\with b$}
      \RightLabel{\hypertarget{rulepair}{\scriptsize \textsc{Pair}}}
      \BinaryInfC{$\Gamma \vdash \pair{e_1}{e_2} : \tau_1 \tensor \tau_2\with a+b$}
    \end{bprooftree}\!
    \begin{bprooftree}
      \AxiomC{$\Gamma \vdash e_1 : \tau_1\with a$}
      \AxiomC{$\Gamma, x : \tau_1 \vdash e_2 : \tau\with b$}
      \RightLabel{\hypertarget{rulelet}{\scriptsize \textsc{Let}}}
      \BinaryInfC{$\Gamma \vdash \letin{x}{e_1}{e_2} : \tau\with a+b$}
    \end{bprooftree}
    \\
    \begin{bprooftree}
      \AxiomC{$\Gamma \vdash e : \B\with a$}
      \AxiomC{$\Gamma \vdash e_1 : \tau\with b$}
      \AxiomC{$\Gamma \vdash e_2 : \tau\with c$}
      \RightLabel{\hypertarget{ruleif}{\scriptsize \textsc{If}}}
      \TrinaryInfC{$\Gamma \vdash \ifte{e}{e_1}{e_2}: \tau\with a+b+c$}
    \end{bprooftree}
    \begin{bprooftree}
      \AxiomC{$\Gamma, x: \tau_1, y : \tau_2, \Gamma' \vdash e : \tau\with a$}
      \RightLabel{\hypertarget{ruleex}{\scriptsize \textsc{Ex}}}
      \UnaryInfC{$\Gamma, y: \tau_2, x : \tau_1, \Gamma' \vdash e : \tau\with a$}
    \end{bprooftree}
  \end{gather*}
  \caption{Typing rules for a discrete imprecise probabilistic programming language. The context $\Gamma$ is assumed to be well-formed (no repeated variables) in \RuleVar, \RuleBern, and \RuleKnight; an easy induction shows that the context of every derivable judgement is well-formed.}\label{fig:typerules}
\end{figure}
We interpret types and contexts as objects in $\BImP$ (natural numbers) with the following inductive definition (recall that addition is the monoidal product in $\BImP$):
\[\sem{\B} = 1 \quad \sem{\tau_1 \otimes \tau_2} = \sem{\tau_1} + \sem{\tau_2} \qquad \sem{\emptyset} = 0 \quad \sem{\Gamma,x:\tau} = \sem{\Gamma} + \sem{\tau}.\]
A typeable expression $\Gamma \vdash e: \tau \with a$ is interpreted as a morphism $\sem{\Gamma} \garrow{a} \sem{\tau}$ in $\BImP$. By \cref{thm:axbimp}, we can represent such morphisms with graded string diagrams built out of the generators of $\Circ$ and $\knightbisgen$. In \cref{fig:semantics}, we define the semantics of typeable expressions as such diagrams by induction on the typing derivation (\textit{cf.}\@ \cite[Figure~3]{Piedeleu2025b}).
\begin{figure}
  \begin{subfigure}{0.30\textwidth}
    \RuleVar : \scalebox{0.7}{\input{semvar.tikz}}
  \end{subfigure}
  \begin{subfigure}{0.31\textwidth}
    \RuleWeak : \scalebox{0.7}{\input{semweak.tikz}}
  \end{subfigure}
  \begin{subfigure}{0.32\textwidth}
    \RuleFst : \scalebox{0.7}{\input{semfst.tikz}}
  \end{subfigure}\\
  \begin{subfigure}{0.30\textwidth}
    \RuleBern : \scalebox{0.7}{\input{sembern.tikz}}
  \end{subfigure}
  \begin{subfigure}{0.31\textwidth}
    \RuleKnight : \scalebox{0.7}{\input{semknight.tikz}}
  \end{subfigure}
  \begin{subfigure}{0.32\textwidth}
    \RuleSnd : \scalebox{0.7}{\input{semsnd.tikz}}
  \end{subfigure}\\
  \begin{subfigure}{0.30\textwidth}
    \RulePair : \scalebox{0.7}{\input{sempair.tikz}}
  \end{subfigure}
  \begin{subfigure}{0.31\textwidth}
    \RuleLet : \scalebox{0.7}{\input{semlet.tikz}}
  \end{subfigure}
  \begin{subfigure}{0.31\textwidth}
    \RuleEx : \scalebox{0.7}{\input{semex.tikz}}
  \end{subfigure}\\
  \begin{subfigure}{0.38\textwidth}
    \RuleIF : \scalebox{0.7}{\input{semif.tikz}}
  \end{subfigure}
  \caption{Inductive definition of the semantics of our language. The trapezium node labelled $\mathtt{if}$ is syntactic sugar for a diagram of type $1 + n + n \garrow{0} n$ for any $n\in \N$. Intuitively, it reads the input of the first wire and when that is true, it outputs the inputs of the next $n$ wires, and otherwise, it outputs the inputs of the last $n$ wires. It is defined as the multiplexer in~\cite[Examples 2.1 and 2.2]{Piedeleu2025b} (see also \cref{appendix:bimp} for the case $n=1$).}\label{fig:semantics}
\end{figure}

The upshot of giving a diagrammatic semantics to these programs is that we can show some properties of the language using diagrammatic reasoning. We prove that our language fulfils \textit{Desideratum~1} of~\cite{LiellCock2025}, namely, it is commutative and affine.
\begin{thm}\label{thm:desiderata}
  The programs below, within appropriate contexts, are "equal up to regrading" (denoted $\utg$ here). For associativity, equality holds without "regrading".
  Note that, for commutativity and weakening, $u$ does not use $x$, as is explicit from the typing.
  \begin{align*}
    \Gamma \vdash t : \tau_1 \with a \qquad \Gamma, x: \tau_1         & \vdash u : \tau_2 \with b \qquad \Gamma, y : \tau_2 \vdash v: \tau \with c                             \\*
    \letin{x}{t}{(\letin{y}{u}{v})}                                 & \utg \letin{y}{(\letin{x}{t}{u})}{v}                                    &  & \text{(associativity)} \\[0.3em]
    \Gamma \vdash t : \tau_1 \with a \qquad \Gamma \vdash u: \tau_2 & \with b \qquad \Gamma, x : \tau_1, y : \tau_2 \vdash v:\tau \with c                                    \\*
    \letin{x}{t}{(\letin{y}{u}{v})}                                 & \utg \letin{y}{u}{(\letin{x}{t}{v})}                                    &  & \text{(commutativity)} \\*
    \letin{x}{t}{u}                                                 & \utg u                                                                  &  & \text{(weakening)}
  \end{align*}
\end{thm}
\begin{proof}
  The proof is purely an exercise in diagrammatic reasoning. We omit the recurring semantics brackets $\sem{-}$ in the diagrams to improve readability.

  Associativity holds by associativity of the copy generator (and its $n$-ary version), together with the interchange law \cref{eqn:th:interchangeupto} to move the copy generator past $t$. The "regrading" by a symmetry introduced by \cref{eqn:th:interchangeupto} involves the grade $\eword$ of the copy generator, hence it is the identity, and associativity holds on the nose.
  \[\scalebox{0.7}{\input{derivationassoc.tikz}}\]
  Commutativity holds ("up to regrading") by the following derivation relying on commutativity of the copy generator (and its $n$-ary version).
  \[\scalebox{0.7}{\input{derivationcomm.tikz}}\]
  Weakening holds ("up to regrading") by the following derivation relying on the properties of the discard generator (and its $n$-ary version): it is unital for the copy generator, it annihilates any $0$-graded term (i.e.\@ $\inlinefig{0.5}{annihilates}$), and it slides through the generator $\knightbisgen$ by \cref{eqn:slidingdel}.
  \[\inlinefig{0.7}{derivationweak}\qedhere\]
\end{proof}

\begin{rem}\label{rem:nohoisting}
  The hoisting equation that corresponds to \textit{Desideratum 2} in~\cite{LiellCock2025} does not hold in our semantics because the following two diagrams are not equal in general.
  \[\scalebox{0.6}{\input{nohoisting.tikz}}\]
  Indeed, the grades of the two diagrams do not match ($a + a \neq a$ unless $a=0$). What is more, the copy generator is not natural, meaning that not all boxes can be copied. Boxes that can be copied are called \textbf{deterministic} in the literature on Markov categories~\cite{Fritz2020}. In $\BImP$, a morphism is deterministic if it does not contain the generators $\inlinefig{0.6}{statep}$ with $p \neq 0,1$ nor~$\knightbisgen$. Abstractly, the reason for the failure of hoisting is that $\ifte{\!\!}{\!\!}{\!\!}$ is not interpreted via a coproduct. In fact, $\BImP$ does not support the full structure of a distributive graded Markov category that $\ImP$ has because, while $\BImP_0= \BStoch$ is a Markov category, there is no coproduct. We mention a potential solution in the conclusion.
\end{rem}

\subsection{Adding Exact Conditioning}\label{sec:conditioning}
The authors of~\cite{Piedeleu2025b} provide another presentation result that builds upon \cref{thm:axbstoch}: they define a "prop" $\PCirc$ by adding a \textbf{conditioning} generator $\inlinefig{0.7}{conditiongen}$ and axioms to $\thCC$, then they show that $\PCirc$ is isomorphic to the "prop" of substochastic matrices~\cite[Section 4]{Piedeleu2025b}---the entries of a column sum up to at most $1$. The semantics of $\inlinefig{0.7}{conditiongen}$ is a substochastic map that sends $\ket{x}\ket{y}$ to $\ket{x}$ if $x = y$ and has no output otherwise (i.e.\@ a distribution with total weight $0$).
In matrix form,
\[\sem{\!\inlinefig{0.7}{conditiongen}} = \scalebox{0.85}{$\begin{bmatrix}
        1 & 0 & 0 & 0 \\
        0 & 0 & 0 & 1
      \end{bmatrix}$}.\]
Intuitively, two probabilistic bits can flow through this generator only when they are equal. Thus, the probability of outputting $\ket{\b1}$ is the probability that $\ket{\b1\b1}$ goes in and similarly for $\ket{\b0}$, while the inputs $\ket{\b1\b0}$ and $\ket{\b0\b1}$ are ignored. The total weight of the distribution in the output can be lower than $1$, hence the semantics is now in substochastic matrices.

Now, $\op{\FInj}$ also embeds inside $\PCirc$, exactly like it does in $\Circ \cong \BStoch$, thus we can apply \cref{thm:gradedpresfrompres} to get an analogue of \cref{thm:axbimp}. Namely, we have a graded presentation of the substochastic variant of $\BImP$. We do not go over the technical details, which differ very little from the presentation of $\BImP$.

When programming in this new "graded category", the $\observe{\!\!}$ primitive is available. It obeys the following typing rule and diagrammatic translation.
\[\Gamma, x: \tau \vdash  \observe{x}: \tau\with 0\quad  \rightsquigarrow \quad\inlinefig{0.8}{observe}\]
According to the semantics of the conditioning generator, $\observe{x}$ puts a condition on the program to require that $x$ is true (i.e.\@ equal to a product of $\ket{\b1}$s). This is an essential component of probabilistic programs. With this functionality, we can showcase two examples that underscore the role of nondeterminism in imprecise probability.
\AtBeginEnvironment{quotation}{\small}
\begin{exa}[Boy or Girl]\label{exmp:boyorgirl}
  The \textit{Boy or Girl paradox} is a well-known probabilistic riddle that seems to originate from Martin Gardner's \textit{Mathematical Games} column~\cite{Gardner1959}. We can state it as follows.
  \begin{quotation}
    Your neighbours have two kids and at least one of them is a girl. Assuming that the sex of each kid is assigned with a fair coin flip ($50\%$ chance of boy and $50\%$ chance of girl) independently of each other, what is the probability that both kids are girls?
  \end{quotation}
  We model this riddle with a probabilistic program, which we immediately translate to a diagram in $\BImP$ (\cref{fig:boyorgirl} on the left).
  This conveys an operational interpretation of the riddle. First, the sex of both kids is stored in two probabilistic bits determined by two independent coin flips represented with the generator $\inlinefig{0.6}{statehalf}$. Then, in parallel, we require the $\mathtt{or}$ of these bits to be $\b1$, meaning that at least one bit is a $\b1$ (read \textit{a girl}), and we output the bits. Applying the functor $\sem{-}$ yields a $4 \times 1$ matrix corresponding to the subdistribution $\sfrac{1}{4}\ket{\b1\b1} + \sfrac{1}{4}\ket{\b1\b0} + \sfrac{1}{4}\ket{\b0\b1}$. To answer the riddle, we take the probability of having two girls relative to the total probability: $\frac{\sfrac{1}{4}}{\sfrac{3}{4}}= \sfrac{1}{3}$.
  \begin{figure}
    \scalebox{0.85}{\input{boyorgirl1.tikz}}\hspace{10em} \scalebox{0.85}{\input{boyorgirl2.tikz}}
    \caption{Graphical representation of two versions of the Boy or Girl paradox.}\label{fig:boyorgirl}
  \end{figure}

  A variant of this riddle can be stated as follows.
  \begin{quotation}
    Your neighbours have two kids: a toddler and a teen. You meet one of them, and she is a girl. What is the probability that both kids are girls?
  \end{quotation}
  In this version, there are two possible events, either you learn that the toddler is a girl, or you learn that the teen is a girl. Both events give you more information than in the original riddle. This riddle is modelled in \cref{fig:boyorgirl} on the right.
  As before, two bits represent the sex of each kid. We condition on the first bit (say the toddler) being a $\b1$ and, in parallel, we condition on the second bit being a $\b1$. We pick one of these two possibilities with an $\mathtt{if}$ block whose \textit{condition} wire takes a nondeterministic bit. Applying the functor $\sem{-}$ yields two $4 \times 1$ submatrices corresponding to the subdistributions $\sfrac{1}{4}\ket{\b1\b1} + \sfrac{1}{4}\ket{\b1\b0}$ and $\sfrac{1}{4}\ket{\b1\b1} + \sfrac{1}{4}\ket{\b0\b1}$. In both cases, the relative probability of the neighbours having two girls is $\sfrac{1}{2}$.
\end{exa}

\begin{rem}
  Adding the $\observe{\!\!}$ construct to our language breaks affineness, namely, the weakening property of \cref{thm:desiderata} does not hold any more. Indeed, our proof of that property uses the naturality property ($\inlinefig{0.75}{annihilates}$) which does not hold for the new generator $f = \inlinefig{0.75}{conditiongen}$. Concretely, the equation $\inlinefig{0.75}{conditiondel}$ cannot be true because if you condition on two bits being equal before discarding them, you are still putting some constraint on the rest of the process---that it led to those two bits being equal. The other two properties in \cref{thm:desiderata}, associativity and commutativity, still hold.
\end{rem}

\section{Case Study: Modular Presentation for Gaussian Maps}\label{sec:gauss}
Our second application of graded string diagrams concerns continuous probability, specifically Gaussian probability. In \cite{Stein2026}, the authors give a string-diagrammatic presentation of $\Gauss$, the category of Gaussian maps, whose morphisms are a special sort of stochastic process. We show that their presentation can be recovered using the "para construction" (\cref{exmp:gradedpara}) together with our modular completeness result (\cref{thm:gradedpresfrompres}). We begin by recalling their result and the necessary background.

\subsection{Affine and Gaussian Maps}\label{sec:affinegauss}

Gaussian maps are essentially affine maps with additive noise modelled by a (multivariate) normal distribution. More precisely, they are random functions between Euclidean spaces (powers of $\intro*\R$, the real numbers) of the form $f(x)=Ax+z$, where $A$ is an affine map and $z$ is a noise vector distributed according to a standard normal distribution. Gaussian maps form a "prop" $\Gauss$ built similarly to the "prop" $\Aff$ of affine maps. To define $\Gauss$, we first recall the definitions of $\Aff$ and normal distributions.

\begin{defi}[$\Aff$]\label{defn:Aff}
  The "prop" $\Aff$ has natural numbers as objects and \textbf{affine maps} between Euclidean spaces as morphisms. A morphism $m \rightarrow n$ is a pair $(A,\mean)$, where $A\in\R^{n\times m}$ and $\mean\in\R^n$.\footnote{When $n$ or $m$ is $0$, there is exactly one empty matrix $\emat \in \R^{n\times m}$ and one empty vector $\emat \in \R^0$.} It is interpreted as the function $x\mapsto Ax+\mean$. Composition and monoidal product are given by 
  \[
    (A,\mean)\compsimple(B,\meanb) = (BA,\,B\mean+\meanb)
    \qquad\text{and}\qquad
    (A,\mean)\tensor(B,\meanb) = (A\oplus B,\,\mean\oplus\meanb),
  \]
  where $\oplus$ denotes direct sum, i.e.\@ $A\oplus B$ is the block matrix $\scalebox{0.85}{$\begin{bmatrix}A&\zmat\\\zmat&B\end{bmatrix}$}$ and $\mean\oplus\meanb = \scalebox{0.85}{$\begin{bmatrix}\mean\\\meanb\end{bmatrix}$}$.

  Letting $\idmat{n}$ denote the identity $n\times n$ matrix, the identity affine map $n \rightarrow n$ is $(\idmat{n},0)$. The symmetry $\sigma_{n,m}:n+m\rightarrow m+n$ is $(\scalebox{0.85}{$\begin{bmatrix}\zmat& \idmat{n}\\\idmat{m}&\zmat\end{bmatrix}$}, 0)$, which is interpreted as the function sending $\scalebox{0.85}{$\begin{bmatrix}x\\y\end{bmatrix}$}$ to $\scalebox{0.85}{$\begin{bmatrix}y\\x\end{bmatrix}$}$ for all $x\in\R^n$ and $y\in\R^m$.
\end{defi}

\begin{defi}[$\normd$]\label{defn:normal}
  Given $n\in\N$, $\mean \in \R^n$, and $\cov \in \R^{n\times n}$ positive semidefinite\footnote{A matrix $\cov \in \R^{n\times n}$ is called \textbf{positive semidefinite} if $\tsp{x}\cov x \geq 0$ for all $x \in \R^n$.}, the \textbf{(multivariate) normal distribution} $\normd(\mean,\cov)$ on $\R^n$ is the probability distribution of a random vector $z\in\R^n$ with mean $\mathbb{E}[z]=\mean$ and covariance $\mathrm{Cov}(z)=\mathbb{E}[(z-\mean)\tsp{(z-\mean)}]=\cov$. The \textbf{standard normal distribution} on $\R^n$ is $\normd(0,\idmat{n})$. We write $z \sim \normd(\mean,\cov)$ to say that $z\in \R^n$ is distributed according to the distribution $\normd(\mean,\cov)$.
\end{defi}

\begin{defi}[$\Gauss$]\label{defn:gauss}
  The "prop" $\Gauss$ has natural numbers as objects and \textbf{Gaussian maps} as morphisms. A Gaussian map $m\rightarrow n$ is a triple $(A,\mean,\cov)$, where $A\in\R^{n\times m}$, $\mean\in\R^n$, and $\cov\in\R^{n\times n}$ is positive semidefinite. It is interpreted as the random function $x\mapsto Ax+z$, where $z \sim \normd(\mean,\cov)$. Composition and monoidal product are given by
  \[
    (A,\mean,\cov)\compsimple(B,\meanb,\covb)=(BA,B\mean+\meanb,B\cov \tsp{B}+\covb)
    \text{ and }
    (A,\mean,\cov)\tensor(B,\meanb,\covb)=(A\oplus B,\mean\oplus\meanb,\cov\oplus\covb).
  \]
  A Gaussian map is called \textbf{deterministic} if its covariance is the zero matrix $\zmat$. Deterministic Gaussian maps are precisely affine maps via the inclusion $\Aff \hookrightarrow \Gauss$ defined by $(A,\mean)\mapsto (A,\mean,\zmat)$. The identities and symmetries of $\Gauss$ are inherited from $\Aff$ via this inclusion: the identity Gaussian map $n\to n$ is $(\idmat{n},0,\zmat)$, and the symmetry $\sigma_{n,m}:n+m\to m+n$ is $(\scalebox{0.85}{$\begin{bmatrix}\zmat& \idmat{n}\\\idmat{m}&\zmat\end{bmatrix}$}, 0,\zmat)$. This makes the inclusion above into a symmetric strict monoidal functor.
\end{defi}
We recall the string diagrammatic presentations of $\Aff$ and $\Gauss$ given in \cite{Bonchi2019c} and \cite{Stein2026} respectively.
\begin{defi}[$\thAff$]
  The "monoidal theory" $\thAff$ has generators $\inlinefig{0.75}{copygen}$, $\inlinefig{0.75}{deletegen}$, $\inlinefig{0.75}{add}$, $\inlinefig{0.75}{zero}$, $\inlinefig{0.75}{one}$, and $\inlinefig{0.75}{scale}$ for all $k \in \R$, and its axioms are listed in \cref{fig:affeqns} (reproducing \cite[\textit{Affine fragment} of Figure~1]{Stein2026}).
\end{defi}
\begin{thmC}[{\cite{Bonchi2019c}}]\label{thm:Affpresentation}
  The "monoidal theory" $\thAff$ presents $\Aff$.
\end{thmC}
In particular, for any matrix $A \in \R^{n\times m}$, we can choose a term in $\thAff$ that is mapped to the affine map $(A,0)$ under the isomorphism $\sync{\thAff} \cong \Aff$. In the developments below, we denote this term $\affdiag{A}$. Note that the specific choice of term will be immaterial since any two terms mapped to $(A,0)$ must be equal as diagrams in $\thAff$.

The presentation of $\Gauss$ extends that of $\thAff$ with a single generator representing a source of random noise, together with axioms that ensure this noise follows a standard normal distribution.
\begin{defi}[$\thGauss$]
  The "monoidal theory" $\thGauss$ has the generators of $\thAff$ as well as $\gaussgen$ and its axioms are those of $\thAff$ plus the following axioms where $\theta$ ranges over $[0,2\pi)$.\\\setcounter{theqn}{0}
  \begin{minipage}{0.3\textwidth}
    \begin{equation}\label{eqn:normaldelete}
      \input{normaldelete.tikz} \stepcounter{theqn}\tag{N\arabic{theqn}}
    \end{equation}
  \end{minipage}
  \begin{minipage}{0.69\textwidth}
    \begin{equation}\label{eqn:normalrotate}
      \input{normalrotate.tikz} \quad \text{ with } R_{\theta} \coloneq \begin{bmatrix}
        \cos\theta&-\sin\theta\\\sin\theta&\cos\theta\end{bmatrix} \stepcounter{theqn}\tag{N\arabic{theqn}}
    \end{equation}
  \end{minipage}\\
\end{defi}
\begin{thmC}[{\cite[Theorem~1]{Stein2026}}]\label{thm:existinggausspresentation}
  The "monoidal theory" $\thGauss$ presents $\Gauss$.
\end{thmC}

\subsection{Noise as Parameters}\label{sec:gaussnoise}
The completeness proof of \cref{thm:existinggausspresentation} relies on a key insight that resembles our proof of modular completeness (\cref{thm:gradedpresfrompres}). Indeed, in our modular presentation of $\Para{\G}{\catC}$, we can drag the extra graded generators $\knightgenc$ all the way to the left (\cref{lem:decompositiongradedterm}), and in $\thGauss$, we can also drag the generator $\gaussgen$ all the way to the left (simply using the laws of "SMCs" in \cref{fig:axiomsvyncat}). Formally, for any term $f$ in $\thGauss$, there is a term $f'$ with no occurrence of $\gaussgen$ (i.e.\@ a term in $\thAff$) such that
\begin{equation}\label{eqn:gaussdecompose}
  \input{decomposegauss.tikz}
\end{equation}
Here, $\gaussgen$ may represent a power of the generator (i.e.\@ $\gaussgen \otimes \stackrel{n}{\cdots} \otimes \gaussgen$ for some $n$).

In terms of semantics, this says that a Gaussian map can be implemented as an affine computation with extra noisy parameters. Let us unroll the details. Let $f: \R^m \rightarrow \R^n$ be defined by $f(x) = Ax + z$ with $z \sim \normd(\mean,\cov)$. We can first decompose the noisy part by letting $z' \sim \normd(0,\cov)$ so that $z$ and $\mean + z'$ are identically distributed. Then, we note that for any matrix $L \in \R^{n\times k}$ such that $L\tsp{L}= \cov$, 
$Lz''$ and $z'$ are identically distributed when $z'' \sim \normd(0,\idmat{k})$. Therefore, the affine map $f':\R^k\oplus \R^m \rightarrow \R^n$ defined by $f'(z,x) = Ax+ \mean + Lz$ implements the Gaussian map $f$ through the extra parameters $z$. Namely, $f'(z,x)$ and $f(x)$ are identically distributed when $z \sim \normd(0,\idmat{k})$.

The construction of $f'$ clearly separates the deterministic affine computation in $f$ and its source of noise. It also introduces a choice of implementation: different matrices $L$, and even different numbers of extra parameters (values of $k$), may determine the same Gaussian map. More concretely, any affine map that takes in standard normally distributed noise and produces standard normally distributed noise can be precomposed to $f'$ and the result will still implement $f$. Using the composition of $\Gauss$, we can characterise such maps:
\begin{equation}\label{eqn:charcoiso}
  (\emat,0,\idmat{m}) \compsimple (A,\mean,\zmat) = (\emat,0,\idmat{n}) \Longleftrightarrow \mean = 0 \text{ and }A\tsp{A} = \idmat{n}.
\end{equation}
Such affine maps are called \textbf{coisometries}. They form a subcategory of $\Aff$ that will play the role of the grading category.

\begin{defi}[$\CoIso$]\label{defn:coiso}
  We denote with $\CoIso$ the subprop of $\Aff$ containing all natural numbers and only the affine maps $(A,0): m \rightarrow n$ where $A$ is a coisometry. It is straightforward to check that identities, symmetries, and compositions and direct sums of coisometries are coisometries.
\end{defi}

We note that the object $0$ is strictly terminal in $\CoIso$. Indeed, for any $k$, the empty matrix $\emat \in \R^{0\times k}$ is a coisometry because $\emat\tsp{\emat} \in \R^{0\times 0}$ must be the identity matrix on $0$, and the empty matrix $\emat \in \R^{k\times 0}$ is not a coisometry because $\emat\tsp{\emat} = \zmat \in \R^{k \times k}$.
\begin{rem}
  This last equality ($\emat\tsp{\emat} = \zmat$) is better justified diagrammatically. A folklore result in graphical linear algebra (see e.g.\@ \cite[Remark~3.12]{zanasi:tel-01218015}), which also applies in $\thAff$, is that whenever a term $f$ represents a matrix $A$ (a morphism $(A,0)$ in $\Aff$), the term $\tsp{f}$ obtained by flipping $f$ horizontally and swapping the colours of white and black generators (e.g.\@ $\inlinefig{0.75}{deletegen} \leftrightarrow \inlinefig{0.75}{zero}$) represents the transpose $\tsp{A}$. Therefore, we have the following diagrammatic representations for the empty matrix $\emat \in \R^{k\times 0}$, its transpose $\tsp{\emat} \in \R^{0\times k}$, and the composition $\emat\tsp{\emat}$, and one can then check that the semantics of the latter is indeed the zero matrix $\zmat$.
  \[\input{emattranspose.tikz}\]
\end{rem}

We now describe the "graded prop" $\GaussSharp$ whose morphisms will be explicit implementations of Gaussian maps as affine maps with noisy parameters. It arises as the "actegory" induced by the inclusion $\CoIso \hookrightarrow \Aff$ (recall \cref{exmp:selfaction}).

\begin{defi}[$\GaussSharp$]\label{defn:gausssharp}
  We define the $\CoIso$-"graded prop" $\GaussSharp\coloneq\Para{\CoIso}{\Aff}$, where the action $\act: \CoIso \times \Aff \rightarrow \Aff$ is given by $k \act n = k+n$ on objects and $(A,0) \act (B,\mean) = (A \oplus B, 0 \oplus \mean)$ on morphisms.
  Its objects are natural numbers, and a morphism $m\garrow{k}n$ is specified by a $(k+m) \times n$ matrix and an $n$-dimensional vector. In order to better match the intuition of grading as extra parameters, we choose to decompose the matrix into $A \in \R^{m \times n}$ and $L \in \R^{k\times n}$, so that a morphism $n \garrow{k} m$ is specified by $(A,\mean,L)$, which can be interpreted as the affine map $(z,x)\mapsto Ax+\mean + Lz$.

  For $(t,0):k'\rightarrow k$ in $\CoIso$, "regrading" is
  \[ (t,0) \regrade(A,\mean,L)=(A,\mean,Lt). \]
  Composition and monoidal product are defined as follows:
  \[(A,\mean,L) \comp (B,\meanb,M) = (BA,B\mean+\meanb,\begin{bmatrix}BL&M\end{bmatrix}) \text{ and } (A,\mean,L) \otimes (B,\meanb,M) = (A\oplus B,\mean\oplus \meanb, L\oplus M)\]
  Identities and symmetries have grade $0$ and are inherited from $\Aff$.
\end{defi}

The decomposition of a Gaussian map $f$ into $f'$ described above justifies viewing a morphism $(A,\mean,L)$ in $\GaussSharp$ as an implementation of the Gaussian map $(A,\mean, L\tsp{L})$. While several morphisms in $\GaussSharp$ can implement the same Gaussian map, we will show that any two such morphisms must be "regradings" of each other by a coisometry. In other words, identifying different implementations of the same map is equivalent to identifying different "regradings" of a morphism. Consequently, $\Gauss$ is isomorphic to the "regrading collapse" of $\GaussSharp$. The rest of this subsection is dedicated to showing this isomorphism.
\begin{lem}\label{lem:charcoisoregrade}
  Let $L \in \R^{n\times k}$ and $M \in \R^{n \times k'}$. The following are equivalent,
  \begin{enumerate}
    \item $L\tsp{L} = M\tsp{M}$
    \item $Lt = Ms$ for coisometries $t: m \rightarrow k$ and $s: m \rightarrow k'$.
  \end{enumerate}
\end{lem}
\begin{proof}
  (2) $\Rightarrow$ (1) Suppose $Lt=Ms$ with $t\tsp{t}=\idmat{k}$ and $s\tsp{s}=\idmat{k'}$. Then
  \[
    L\tsp{L} = L(t\tsp{t})\tsp{L} = (Lt)\tsp{(Lt)} = (Ms)\tsp{(Ms)} = M(s\tsp{s})\tsp{M} = M\tsp{M}.
  \]
  (1) $\Rightarrow$ (2). Without loss of generality, suppose $k \geq k'$ and let $\pi \coloneq \begin{bmatrix} \idmat{k'} & \zmat \end{bmatrix}\in \R^{k' \times k}$ be the projection matrix onto the first $k'$ coordinates. A simple computation shows $\pi$ is a coisometry. Then, letting $\tilde{M} \coloneq M\pi \in \R^{n\times k}$, we have the following derivation
  \[ \tilde{M}\tsp{\tilde{M}} = M\pi\tsp{\pi}\tsp{M} = M\tsp{M} = L\tsp{L}. \]
  By \cite[Proposition~6]{Stein2026arxiv},\footnote{\cite{Stein2026arxiv} is a version of \cite{Stein2026} containing an appendix with detailed proofs.} there is a matrix $R\in\mathbb R^{k\times k}$ such that $R\tsp{R} = \idmat{k}$ and $LR=\tilde{M}$, that is, $LR = M\pi$. We conclude that $Lt = Ms$ where $t = R$ and $s =\pi$ are coisometries.
\end{proof}
A corollary of \cref{lem:charcoisoregrade} is a characterisation of when two morphisms in $\GaussSharp$ are in the same equivalence class in $\Coll(\GaussSharp)$.
\begin{cor}\label{cor:charcollcong}
  For any two morphisms $(A,\mean,L)$ and $(B,\meanb,M)$ in $\GaussSharp$, $(A,\mean,L) \collcong (B,\meanb,M)$ if and only if $A=B$, $\mean = \meanb$ and $L\tsp{L} = M\tsp{M}$.
\end{cor}
\begin{proof}
  ($\Leftarrow$) This follows directly from (1) $\Rightarrow$ (2) in \cref{lem:charcoisoregrade}.

  ($\Rightarrow$) The relation between morphisms in $\GaussSharp$ defined by $(A,\mean,L) \sim (B,\meanb,M)$ whenever $A = B$, $\mean = \meanb$ and $L\tsp{L} = M\tsp{M}$ is clearly transitive. Thus, it is enough to show $(A,\mean, L) \utg (B,\meanb,M)$ implies $(A,\mean,L) \sim (B,\meanb,M)$ because $\collcong$ is the transitive closure of $\utg$. Without loss of generality, suppose that $t \regrade (A,\mean,L)= (B,\meanb, M)$ for some coisometry $t$. By definition of "regrading" in the "para construction", this means that $(A,\mean,Lt) = (B,\meanb, M)$. Thus, $A = B$, $\mean=\meanb$, and we can use (2) $\Rightarrow$ (1) in \cref{lem:charcoisoregrade} to get that $L\tsp{L} = M\tsp{M}$.
\end{proof}
\begin{thm}\label{thm:gausscollapse}
  The mapping $\collclass{(A,\mean,L)} \mapsto (A,\mean,L\tsp{L})$ yields an isomorphism of "props" $\Coll(\GaussSharp) \cong \Gauss$.
\end{thm}
\begin{proof}
  We define an identity-on-objects functor $F: \Coll(\GaussSharp) \rightarrow \Gauss$ by $F\collclass{(A,\mean,L)} = (A,\mean,L\tsp{L})$. It is well-defined by \cref{cor:charcollcong}. We also need to check it is a symmetric strict monoidal functor, namely, that it preserves composition and monoidal product of maps. We have the following derivations.
  \begin{align*}
    F(\collclass{(A,\mean,L)} \compsimple \collclass{(B,\meanb,M)})
     & = F(\collclass{(A,\mean,L)\comp(B,\meanb,M)})                                                     &  & \text{$\collcong$ is a congruence}         \\
     & = F(\collclass{(BA, B\mean+\meanb, \begin{bmatrix}BL & M\end{bmatrix})})                          &  & \text{def.\@ of $\comp$ in $\GaussSharp$}  \\
     & = (BA, B\mean+\meanb, \begin{bmatrix}BL & M\end{bmatrix}\tsp{\begin{bmatrix}BL & M\end{bmatrix}}) &  & \text{def.\@ of $F$}                       \\
     & = (BA, B\mean+\meanb, B(L\tsp{L})\tsp{B}+M\tsp{M})                                                                                                \\
     & = (A,\mean,L\tsp{L}) \compsimple (B,\meanb,M\tsp{M})                                              &  & \text{def.\@ of $\compsimple$ in $\Gauss$} \\
     & = F\collclass{(A,\mean,L)} \compsimple F\collclass{(B,\meanb,M)}                                  &  & \text{def.\@ of $F$}
  \end{align*}
  \begin{align*}
    F(\collclass{(A,\mean,L)} \otimes \collclass{(B,\meanb,M)})
     & = F(\collclass{(A,\mean,L)\tensor (B,\meanb,M)})               &  & \text{$\collcong$ is a congruence}          \\
     & = F(\collclass{(A\oplus B, \mean\oplus\meanb, L\oplus M)})     &  & \text{def.\@ of $\tensor$ in $\GaussSharp$} \\
     & = (A\oplus B, \mean\oplus\meanb, (L\oplus M)\tsp{(L\oplus M)}) &  & \text{def.\@ of $F$}                        \\
     & = (A\oplus B, \mean\oplus\meanb, L\tsp{L}\oplus M\tsp{M})                                                       \\
     & = (A,\mean,L\tsp{L}) \otimes (B,\meanb,M\tsp{M})               &  & \text{def.\@ of $\otimes$ in $\Gauss$}      \\
     & = F\collclass{(A,\mean,L)} \otimes F\collclass{(B,\meanb,M)}   &  & \text{def.\@ of $F$}
  \end{align*}
  It remains to show it is fully faithful. For fullness, it is a standard fact in linear algebra (see e.g.\@ \cite[7.38]{LinearAlgebraAxler}) that any positive semidefinite matrix $\cov$ is of the form $L\tsp{L}$ for some matrix $L$, meaning that $(A,\mean,\cov) = F(A,\mean,L)$. Faithfulness follows by \cref{cor:charcollcong}.
\end{proof}

\subsection{Modular Presentation of \texorpdfstring{$\Gauss$}{Gauss}}\label{sec:gausspres}

\cref{thm:gausscollapse} identifies $\Gauss$ as the "regrading collapse" of the "para construction" $\Para{\CoIso}{\Aff}$. We have general results that provide string diagrammatic presentations of the "regrading collapse" (\cref{thm:collapsepresentation}) and the "para construction" (\cref{thm:gradedpresfrompres}) provided we have presentations of the input categories. The presentation of $\Aff$ was recalled in \cref{thm:Affpresentation}, and we now turn to $\CoIso$. This relies on the following lemma.

For $\theta\in [0,2\pi)$, recall that $R_\theta$ denotes the matrix $\begin{bmatrix}\cos\theta&-\sin\theta\\\sin\theta&\cos\theta\end{bmatrix}$ that rotates a two-dimensional plane about the origin by $\theta$ radians. A \textbf{Givens rotation} on $\R^m$ is a linear map which acts as $R_\theta$ on a plane spanned by two standard basis vectors and as the identity on the remaining basis vectors. Equivalently, it is a matrix of the form $\tsp{P}(R_\theta\oplus\idmat{m-2})P$, where $P$ is a permutation matrix.

\begin{lem}\label{lem:decomposecoiso}
  Any coisometry $A \in \R^{n\times m}$ can be decomposed into $\pi R_1\cdots R_k D$, where $\pi = \begin{bmatrix}\idmat{n} & \zmat \end{bmatrix}$, each $R_i$ is a Givens rotation, and $D$ is a diagonal matrix with entries in $\pm 1$.
\end{lem}
\begin{proof}
  Let $A \in \R^{n\times m}$ be a coisometry, i.e.\@ it satisfies $A\tsp{A} = \idmat{n}$. By the dual of \cite[7.49]{LinearAlgebraAxler}, its rows are orthonormal (they are vectors of norm $1$ and pairwise orthogonal). We can add rows to form an orthonormal basis of $\R^m$ (by \cite[6.36]{LinearAlgebraAxler}), and we let $B \in \R^{m\times m}$ be the resulting matrix. It satisfies $B\tsp{B} = \idmat{m} = \tsp{B}B$, so $B$ is an orthogonal matrix, and $A = \pi B$ where $\pi = \begin{bmatrix}\idmat{n} & \zmat \end{bmatrix} \in \R^{n\times m}$ is the projection onto the first $n$ coordinates. Finally, in the proof of \cite[Proposition~13]{Stein2026arxiv}, it is shown that any orthogonal matrix can be decomposed as $R_1\cdots R_k D$, where each $R_i$ is a Givens rotation, and $D$ is a diagonal matrix with entries in $\pm 1$. Therefore, we conclude that our coisometry $A$ can be decomposed as $\pi R_1\cdots R_k D$, as desired.
\end{proof}

\begin{defi}\label{defn:coisopresentation}
  The "monoidal theory" $\thCoIso$ has generators $\inlinefig{0.75}{deletegen}$, $\inlinefig{0.75}{scalarnegone}$, and $\inlinefig{0.75}{rotategen}$ for all $\theta \in [0,2\pi)$. To specify its axioms, interpret the generators of $\thCoIso$ as diagrams in $\thAff$, extend to terms in $\thCoIso$, and impose all the equations that hold inside $\vync{\thAff}$. Let $\mathcal{F}$ be the free prop on these generators and $\iota$ be the symmetric strict monoidal functor $\mathcal{F}\rightarrow \vync{\thAff}$ extending the assignment below, where $R_{\theta}$ is the Givens rotation for the angle $\theta \in [0,2\pi)$.
  \[
    \iota(\inlinefig{0.75}{deletegen}) = \inlinefig{0.75}{deletegen} \qquad
    \iota(\inlinefig{0.75}{scalarnegone}) = \inlinefig{0.75}{scalarnegone} \qquad
    \iota(\inlinefig{0.75}{rotategen}) = \inlinefig{0.75}{diagRtheta}
  \]
  The axioms of $\thCoIso$ are all equations $f=g$ between parallel terms $f,g:m \rightarrow n$ in $\mathcal{F}$ such that $\iota(f)=\iota(g)$ (as morphisms in $\vync{\thAff}$, not as terms).
\end{defi}
By postulating that the axioms of $\thCoIso$ are precisely the equations that hold under the assignment $\iota$, it automatically becomes a model and we get a functor $\iota: \vync{\thCoIso} \rightarrow \vync{\thAff}$. We show that the functor $\bar{\iota}:\vync{\thCoIso} \rightarrow \Aff$, obtained by post-composing $\iota$ with the isomorphism in \cref{thm:Affpresentation}, is an isomorphism of "props".

\begin{thm}\label{thm:coisopres}
  The "monoidal theory" $\thCoIso$ presents $\CoIso$. 
\end{thm}
\begin{proof}
  Using functorial semantics, we obtained an identity-on-objects symmetric strict monoidal functor $\iota: \vync{\thCoIso} \rightarrow \vync{\thAff}$. It is definitionally faithful because $\iota(f) = \iota(g)$ implies $f=g$ is an axiom of $\thCoIso$. Thus, the composite $\bar{\iota}:\vync{\thCoIso} \xrightarrow{\iota} \vync{\thAff} \cong \Aff$ is a faithful identity-on-objects symmetric strict monoidal functor. It remains to check that its image is exactly $\CoIso$.

  The image of each generator is a coisometry and so are the identity and symmetry. Since coisometries are closed under composition and monoidal product, the recursive definition of $\bar{\iota}$ forces its image to lie in $\CoIso$.

  The hard part requires showing that for every coisometry $A \in \R^{n \times m}$, there is a term $f:m \rightarrow n$ in $\thCoIso$ such that $\bar{\iota}(f)=(A,0)$. By \cref{lem:decomposecoiso}, we can decompose $A$ as $\pi R_1\cdots R_kD$. Thus, it suffices to construct a term in $\thCoIso$ representing each factor, and their composition will represent $A$.

  First, a diagonal matrix $D$ whose entries are $d_1, \dots, d_m$ is equal to a direct sum of $1\times 1$ matrices: $\bigoplus_{i=1}^m\begin{bmatrix}d_i\end{bmatrix}$. Hence, when each $d_i = \pm 1$, we can form a term representing $D$ as a monoidal product of several identity wires and copies of $\inlinefig{0.75}{scalarnegone}$ in the appropriate order.

  Next, let $R$ be a Givens rotation. By definition, there are a permutation matrix $P$ and an angle $\theta \in [0,2\pi)$ such that $R=\tsp{P}(R_\theta\oplus\idmat{m-2})P$. Permutation matrices are obtained in $\Aff$ as compositions and products of symmetries and identities, so if we take the term $p$ representing this composition and view it inside $\thCoIso$, it will satisfy $\bar{\iota}(p) = P$. Similarly, we find the term $p'$ such that $\bar{\iota}(p') = \tsp{P}$. Then the term corresponding to $R$ is $p\compsimple(\inlinefig{0.75}{rotategen} \otimes \id_{m-2})\compsimple p'$.

  Finally, the projection $\pi$ is represented by $\mathrm{id}_n\otimes \bigotimes_{i = 1}^{m-n} \inlinefig{0.75}{deletegen}$.
\end{proof}
\begin{rem}
  While the "monoidal theory" $\thCoIso$ presents $\CoIso$, it does so with `brute force', taking all true equations as axioms. Nevertheless, it is still amenable to diagrammatic reasoning because its generators are simple and it exploits the presentation of $\Aff$. Finding a succinct axiomatisation of $\CoIso$ remains an interesting direction for future work. In any case, the axioms of $\thCoIso$ are absent from our presentation of $\Gauss$, since we will pass to the "flattened theory", in which all grading wires are erased.
\end{rem}

We can now use \cref{thm:gradedpresfrompres} to get a presentation of $\GaussSharp$.
\begin{defi}[$\thGaussSharp$]\label{defn:gausssharppresentation}
  The $\CoIso$-"graded theory" $\thGaussSharp$ contains all the generators and axioms of $\thAff$, regarded as $0$-graded, and one additional generator $\noisegen:0\garrow{1}1$ and the following axioms:\footnote{We draw the generators of $\thCoIso$ vertically and coloured as we did for the presentation of $\BImP$.}
  \begin{equation}\label{eqn:axiomsgausssharp}
    \scalebox{0.8}{\input{slidingdelbis.tikz}} \  \scalebox{0.8}{\input{slidingscaleneg.tikz}} \  \scalebox{0.8}{\input{slidingrotate.tikz}} \text{$\theta \in [0,2\pi)$}\   \scalebox{0.8}{\input{slidingswapbis.tikz}}.
  \end{equation}
\end{defi}
\begin{thm}\label{thm:gausssharppresentation}
  The "graded theory" $\thGaussSharp$ presents $\GaussSharp$.
\end{thm}

Finally, by \cref{thm:collapsepresentation,thm:gausscollapse}, we get a presentation of $\Gauss$ by "flattening" $\thGaussSharp$.

\begin{thm}\label{thm:gaussnewpres}
  The "flattened theory" $\flatt{\thGaussSharp}$, as defined in \cref{defn:flatth}, presents $\Gauss$.
\end{thm}

\begin{exa}
  We draw diagrams for two different implementations in $\GaussSharp$ of the Gaussian map $1 \rightarrow 1$ sending $x \in \R$ to $x+z$, where $z\sim\normd(0,1)$. The first implementation (below left) has grade $1$: it draws one noisy input from the grading wire and directly adds it to the ordinary input wire. The second implementation (below right) has grade $2$: it draws two noisy inputs from the grading wires, adds them together, scales the result by $\sfrac{1}{\sqrt{2}}$, and adds that to the ordinary input.
  \[f_1 = \input{exmpgsharp1.tikz} \qquad \qquad f_2 = \input{exmpgsharp2.tikz}\]
  These clearly represent different morphisms in $\GaussSharp$ because they have different grades. However, $f_2$ is obtained from $f_1$ by regrading with the coisometry $\inlinefig{0.75}{coiso1sqrt2}$ as shown below: the first and second equality are by \cref{eqn:axiomsgausssharp}, the third is by a choice of term representing $R_{\frac{7\pi}{4}}$ (using $\cos\frac{7\pi}{4} = \sfrac{1}{\sqrt{2}}$ and $\sin\frac{7\pi}{4} = \sfrac{-\!1}{\sqrt{2}}$), and the fourth and fifth are by axioms of $\thAff$.
  \[\input{deriveexmpgsharp.tikz}\]
  Therefore, $f_1$ and $f_2$ implement the same Gaussian map.

  This can also be seen semantically. Using the isomorphism $\sync{\thGaussSharp}\cong\GaussSharp$, we find $\sem{f_1} = (1,0,L)$ and $\sem{f_2}= (1,0,M)$ with $L = \begin{bmatrix}1\end{bmatrix}$ and $M=\begin{bmatrix}
      \sfrac{1}{\sqrt{2}} & \sfrac{1}{\sqrt{2}}
    \end{bmatrix}$.
  Since $L\tsp{L} = M\tsp{M}= \idmat{1}$, under the isomorphism of \cref{thm:gausscollapse}, the equivalence classes of $\sem{f_1}$ and $\sem{f_2}$ are both sent to the Gaussian map $(1,0,\idmat{1})$.
\end{exa}

\paragraph{Comparison with \texorpdfstring{\cite{Stein2026}}{[SZPS26]}}
If we examine $\flatt{\thGaussSharp}$ more closely, we find that it differs from the theory $\thGauss$ only by two axioms which are derivable in the latter.

In the "flattened theory" $\flatt{\thGaussSharp}$, $\noisegen$ becomes a generator $0\rightarrow 1$ and we can check that its interpretation coincides with the interpretation of $\gaussgen$ in $\thGauss$. Indeed, going back to the proof of \cref{thm:gradedpresfrompres}, we find that the presentation of \cref{thm:gausssharppresentation} sends $\noisegen$ to $\id_1 \act \id_0 \in \GaussSharp_{1}(0,1)$. That is the identity affine map $\R^1 \oplus \R^0 \rightarrow \R^1$ which is viewed in $\GaussSharp$ as the triple $(\emat,0,\idmat{1})$. Therefore, when applying the isomorphism of \cref{thm:gausscollapse}, we obtain $(\emat,0,\idmat{1})$, which is the interpretation of $\gaussgen$.

All the axioms coming from $\thAff$ remain unchanged under flattening, and they also belong to $\thGauss$. For the new axioms in \cref{eqn:axiomsgausssharp}, flattening does the following:
\begin{itemize}
  \item The first becomes the axiom \cref{eqn:normaldelete} in $\thGauss$.
  \item The second becomes $\raisebox{0.3ex}{\scalebox{0.8}{$\input{normalnegone.tikz}$}}$ which is derivable from the axioms of $\thGauss$ as explained in \cite[Proposition~12]{Stein2026arxiv}.
  \item The third becomes the axiom \cref{eqn:normalrotate} of $\thGauss$.
  \item The fourth becomes $\raisebox{0.3ex}{\scalebox{0.8}{$\input{normalswap.tikz}$}}$ which is derivable from the axiom \cref{eqn:month:overthewire} of "SMCs".
\end{itemize}

Let us also compare the methodologies to arrive at the same presentation of $\Gauss$. The biggest difference is the central role of coisometries in our presentation, whereas orthogonal matrices played a similar role in \cite{Stein2026}.

The characterisation in \cref{eqn:charcoiso} generalises \cite[Proposition~13]{Stein2026arxiv}, which states that standard normal distributions are invariant under orthogonal transformations. However, unlike coisometries, orthogonal matrices do not change the dimension. Thus, taking $\CoIso$ as a grading category allows us to identify implementations of the same Gaussian map with different numbers of parameters. Similarly, \cref{lem:charcoisoregrade} generalises \cite[Proposition~6]{Stein2026arxiv} from orthogonal matrices to coisometries. Although presenting the prop of coisometries requires some bureaucratic work, the essential step is finding suitable generators. The generators for orthogonal matrices already appear implicitly in \cite{Stein2026} through the decomposition result. Moreover, taking \textit{all the equations as axioms} in $\thCoIso$ considerably simplifies the proof and introduces no additional complexity into the final presentation of $\Gauss$: "flattening" erases all grading wires, and hence the details of the chosen presentation of $\CoIso$ disappear from the "flattened theory".

Our proof otherwise closely follows that of \cite{Stein2026} except that their argument for completeness is factored through three conceptually disjoint (reusable) results in our work:
\begin{enumerate}
  \item the presentations of the "regrading collapse" and the "para construction" (\cref{thm:collapsepresentation,thm:gradedpresfrompres}),
  \item the semantic decomposition of $\Gauss$ as $\Coll(\Para{\CoIso}{\Aff})$ (\cref{thm:gausscollapse}), and
  \item the presentation of $\CoIso$ (\cref{thm:coisopres}).
\end{enumerate}

Finally, we emphasise that (2) provides a new categorical understanding of the idea that Gaussian maps are affine maps with additional noisy parameters, making a structural separation between the deterministic and stochastic part of these maps. This seems to be closely connected to the relationship between the CPM and discard constructions in~\cite{Carette2021}. We remark on this in the conclusion.

\section{Causal Interventions via Grading}\label{sec:interventions}
In this section, we introduce a "graded monoidal theory" $\thI$ with the aim of providing a graphical representation of probabilistic programs with causal interventions.
We begin with an overview of the string-diagrammatic treatment of causal models and interventions.

\subsection{String Diagrams and Causal Models}
Below, we reproduce a running example of~\cite{JacobsKZ21} on modelling the causal relationship between smoking and cancer.
\begin{figure}[!ht]
  \begin{minipage}{0.4\textwidth}
    \begin{tikzcd}[cramped,sep=small]
      && {\underline{S}} \\
      H && {\underline{T}} \\
      && {\underline{C}}
      \arrow[from=1-3, to=2-3]
      \arrow[from=2-1, to=1-3]
      \arrow[from=2-1, to=3-3]
      \arrow[from=2-3, to=3-3]
    \end{tikzcd}
  \end{minipage}
  \begin{minipage}{0.5\textwidth}
    \input{smokingscenario.tikz}
  \end{minipage}
  \caption{A Bayesian network and its string diagrammatic translation.}\label{fig:smokingscenario}
\end{figure}

On the left of \cref{fig:smokingscenario} there is a Bayesian network. The nodes labelled $H$, $S$, $T$, and $C$ represent random variables: $H$ stands for hidden factors, while $S$, $T$, and $C$ stand for different observations on a person, respectively, whether they are a \underline{s}moker, whether they have \underline{t}ar in the lungs, and whether they have \underline{c}ancer. The arrows indicate the assumed causal relationships: $H$ influences $S$ and $C$, $S$ influences $T$, and $T$ influences $C$. We note that $S$ may indirectly influence $C$, through $T$, because there is no direct arrow from $S$ to $C$.

On the right of \cref{fig:smokingscenario}, there is a string diagram picturing stochastic matrices (morphisms in $\BStoch$) connected according to the causal structure of the Bayesian network. The outputs of the diagram are the observations $S$, $T$, and $C$, which are underlined in the Bayesian network.

Several papers have explored the diagrammatic approach to study causal interventions~\cite{JacobsKZ21,Yimu2021,Lorenz2023}. In all these works, the diagrams describing a model before and after an intervention are related by a transformation defined externally to the diagrammatic language. Namely, they define a functor\footnote{While~\cite{Yimu2021} does not adopt the functorial point of view, the approach is essentially the same.} which applies to a diagram representing a model and yields a diagram representing that model after the intervention.
For example, one can intervene on $S$ in the model of \cref{fig:smokingscenario} to remove the dependence on $H$ (and hence quantify the degree of dependence of $C$ on $S$). In~\cite{JacobsKZ21}, the authors visualise this by \textit{cutting} $s$ and replacing it with a process that discards the input from $h$ and randomly fixes the output:
\begin{equation}\label{eqn:smokingafterintervention}
  \scalebox{0.8}{\input{smokingafterint.tikz}}
\end{equation}

The usual approach depicted above is \emph{meta-theoretic}. Our aim is to incorporate interventions \emph{within} the diagrammatic syntax with grading. With graded string diagrams, we are able to represent some breaking interventions (in the terminology of~\cite{Lorenz2023}) with a "regrading". We will proceed as follows: in \cref{sec:BInt}, we define the grading category that represents the possible interventions; in \cref{sec:EPara}, we define the "graded category" where our diagrams will be interpreted; and in \cref{sec:diagraminterventions}, we define the diagrammatic syntax and demonstrate our approach by revisiting the smoking scenario discussed above.

\subsection{A Grading Category for Interventions}\label{sec:BInt}

We first need to define the "SMC" $\BInt$ that will take the role of the grading category in the following sections. Intuitively, it allows carrying around the information of what intervention is being done. You can either cut a wire in a diagram and replace it with a fixed distribution on $\B$, or do nothing. Since a distribution $\dist$ on $\B$ is determined by an element $p \in \probs = [0,1]$ such that $\dist = p\ket{\b1} + (1-p)\ket{\b0}$, the set of possible interventions is $\intro*\probstar \coloneq \probs + \{\b*\}$, where $\b*$ means `do nothing'.

The category $\BInt$ is a subcategory of the Kleisli category $\Kl(-+\probstar)$ for the exception monad $-+\probstar$ on $\Set$. \AP Let us give an explicit definition of $\intro*\BInt$ and its "SMC" structure.
\begin{defi}[$\BInt$]
  \
  \begin{itemize}
    \item Objects of $\BInt$ are natural numbers $a\in \N$ identified with the finite sets $\fset{a} = \left\{ 0,\dots, a-1 \right\}$.
    \item A morphism $t \in \BInt(a, b)$ is a function $t:a \rightarrow b+\probstar$ that is \textbf{left-injective} or \textbf{injective on $b$}, namely, if $t(i) = t(j) \in b$, then $i = j$ (the images of distinct $i$ and $j$ can coincide if they belong to $\probstar$). In other words, $t$ is injective on the preimage of $b$. 
    \item The identity on $a$ is the left injection $a \hookrightarrow a+\probstar$.
    \item The composition of $t: a \rightarrow b$ and $s: b \rightarrow c$ is the function $t\compsimple s : a \rightarrow c+\probstar$ defined by $(t\compsimple s)(i) = s(t(i))$ when $t(i) \in b$ and $(t\compsimple s)(i) = t(i)$ otherwise.
  \end{itemize}
  Since the left adjoint $\Set \rightarrow \Kl(-+\probstar)$ preserves colimits, $\Kl(-+\probstar)$ has coproducts given by the disjoint union, yielding a symmetric monoidal structure on the Kleisli category. Furthermore, $\BInt$ inherits this SMC structure because finite sets and left-injective functions are closed under coproducts, and symmetries are left-injective.\footnote{Since $\Set$ is not strict monoidal, the associators and unitors are not directly inherited, they are taken to be identities so that $\BInt$ is a prop.}

  Explicitly, given two morphisms $t: a \rightarrow b$ and $t': a' \rightarrow b'$, their monoidal product acts like $t$ on $a$ and $t'$ on $a'$. It is the composition
  \[
    \begin{tikzcd}
      {a+a'} & {(b+\probstar)+(b'+\probstar)} & {(b+b')+(\probstar+\probstar)} &&& {b+b'+\probstar}
      \arrow["{t+t'}", from=1-1, to=1-2]
      \arrow["\cong"{description}, draw=none, from=1-2, to=1-3]
      \arrow["{\copair{\id_{b+b'},\copair{\id_{\probstar},\id_{\probstar}}}}", from=1-3, to=1-6]
    \end{tikzcd}.\]
\end{defi}

The "SMC" $\BInt$ is similar to the "SMC" $\mathbf{PF}$ of partial functions (the Kleisli category of $-+\{\b*\}$) that was given a string diagrammatic presentation in~\cite{Zanasi2016}. We give a "monoidal theory" $(\sigv,\axv)$ that presents $\BInt$ inspired by that of $\mathbf{PF}$. Its generators are
\[\sigv \coloneq \left\{ \begin{tikzpicture}
	\begin{pgfonlayer}{nodelayer}
		\node [style=black] (0) at (0.75, 0) {};
		\node [style=none] (1) at (-0.25, 0) {};
	\end{pgfonlayer}
	\begin{pgfonlayer}{edgelayer}
		\draw [style=single] (1.center) to (0);
	\end{pgfonlayer}
\end{tikzpicture}
, \begin{tikzpicture}
	\begin{pgfonlayer}{nodelayer}
		\node [style=black] (3) at (-0.75, 0) {};
		\node [style=none] (4) at (0.25, 0) {};
	\end{pgfonlayer}
	\begin{pgfonlayer}{edgelayer}
		\draw [style=single] (4.center) to (3);
	\end{pgfonlayer}
\end{tikzpicture}
 \right\} \cup \left\{ \bintp{p} \mid p \in \probs \right\},\]
and its axioms $\axv$ state that plugging $\begin{tikzpicture}
	\begin{pgfonlayer}{nodelayer}
		\node [style=black] (3) at (-0.75, 0) {};
		\node [style=none] (4) at (0.25, 0) {};
	\end{pgfonlayer}
	\begin{pgfonlayer}{edgelayer}
		\draw [style=single] (4.center) to (3);
	\end{pgfonlayer}
\end{tikzpicture}
$ into $\begin{tikzpicture}
	\begin{pgfonlayer}{nodelayer}
		\node [style=black] (0) at (0.75, 0) {};
		\node [style=none] (1) at (-0.25, 0) {};
	\end{pgfonlayer}
	\begin{pgfonlayer}{edgelayer}
		\draw [style=single] (1.center) to (0);
	\end{pgfonlayer}
\end{tikzpicture}
$ or into $\bintp{p}$ yields the empty diagram (the identity on $0$):
\begin{equation}\label{eqn:axBInt}
  \begin{tikzpicture}
	\begin{pgfonlayer}{nodelayer}
		\node [style=black] (0) at (-0.5, 0) {};
		\node [style=black] (1) at (0.5, 0) {};
	\end{pgfonlayer}
	\begin{pgfonlayer}{edgelayer}
		\draw [style=single] (0) to (1);
	\end{pgfonlayer}
\end{tikzpicture}
 = \inlinefig{0.8}{idzer} \qquad\qquad \begin{tikzpicture}
	\begin{pgfonlayer}{nodelayer}
		\node [style=black] (0) at (-0.65, 0) {};
		\node [style=flipop] (1) at (0.75, 0) {$p$};
	\end{pgfonlayer}
	\begin{pgfonlayer}{edgelayer}
		\draw [style=single] (0) to (1);
	\end{pgfonlayer}
\end{tikzpicture}
 = \inlinefig{0.8}{idzer} \quad (p \in \probs).
\end{equation}
\begin{prop}\label{prop:presBInt}
  The assignment $F$ defined below extends to an isomorphism of "SMCs" $\vync{\sigv,\axv} \cong \BInt$.
  \[F(\inlinefig{1}{genbintcounit}) = i \mapsto \b* \qquad F(\raisebox{0ex}{\scalebox{1}{\bintp{p}}}) = i \mapsto p \in \probs \qquad F(\inlinefig{1}{genbintunit}) = ()\]
\end{prop}
\begin{proof}[Proof idea]
  The proof of this statement follows the ideas in~\cite{Zanasi2016}, but we leave the technical details in \cref{appendix:interventions} since they are tangential to the main story---recall the grading category need not have a presentation. Informally, one can look at a diagram $a \rightarrow b$ in $\vync{\sigv,\axv}$ to determine the function it corresponds to. Indeed, for any $i \in a$, the $i$th input wire will be connected to exactly one of $\inlinefig{0.8}{genbintcounit}$, $\raisebox{0.3ex}{\scalebox{0.8}{\bintp{p}}}$ for some $p \in\probs$, or an output wire. The corresponding function sends $i$ to $\b* \in \probstar$ in the first case, to $p\in \probstar$ in the second, and to $j \in b$ if it is connected to the $j$th output wire. Parts of the diagram that are connected to no input and no output are removed by the axioms \cref{eqn:axBInt}, and they do not affect the corresponding function. Here are two examples:
  \[\scalebox{0.75}{\input{exmpdiagBInt.tikz}} \mapsto \begin{cases}
      0 \mapsto 1 & 2 \mapsto q \\1 \mapsto p & 3 \mapsto 0
    \end{cases}\qquad \scalebox{0.75}{\input{exmpdiagBInt2.tikz}} \mapsto \begin{cases}
      0 \mapsto \b* & 2 \mapsto 3 \\1 \mapsto \b*&3 \mapsto 0
    \end{cases}\]
  This process only looks at which input is connected to which generator or which output. Therefore, any legal deformation of a diagram (via laws of "SMCs") does not affect its corresponding function because they do not disconnect wires.
\end{proof}

We now introduce a technical lemma that will be used to define the "regrading" action of $\BInt$ on the "graded SMC" defined in \cref{sec:EPara}.
\AP It relies on viewing a function $t: a \rightarrow b+\probstar$ as an $\probstar$-predicate transformer. Namely, given $f: b \rightarrow \probstar$, we can define $\intro*\transformer{t}(f) : a \rightarrow \probstar$ by $\transformer{t}(f)(i) = f(t(i))$ if $t(i) \in b$ and $\transformer{t}(f)(i) = t(i)$ otherwise. This defines a function $\transformer{t}: \probstar^{\!\!\!b} \rightarrow \probstar^{\!\!\!a}$, and letting $a$ and $b$ vary, we obtain a symmetric monoidal functor.
\begin{lem}\label{lem:transformerfunctor}
  Mapping an object $a \in \BInt$ to $\probstar^{\!\!\!a} \in \Set$ and a morphism $t\in \BInt(b,a)$ to the function $\transformer{t} \in \Set(\probstar^{\!\!\!a}, \probstar^{\!\!\!b})$ yields a symmetric monoidal functor $(\op{\BInt}, +, 0) \rightarrow (\Set,\times,\mathbf{1})$.
\end{lem}
\begin{proof}
  We first need to show the following equations for all $t: a \rightarrow b+\probstar$, $s: b \rightarrow c+\probstar$, and $t': a' \rightarrow b'+\probstar$.\\[0.4em]
  \begin{minipage}{0.5\textwidth}
    \begin{equation}\label{eqn:transformercomp}
      \transformer{t\compsimple s} = \transformer{t} \circ \transformer{s}
    \end{equation}
  \end{minipage}
  \begin{minipage}{0.49\textwidth}
    \begin{equation}\label{eqn:transformerpar}
      \transformer{t+t'} = \transformer{t}\times \transformer{t'}
    \end{equation}
  \end{minipage}\\[0.4em]
  This only requires a careful analysis of cases.

  For \cref{eqn:transformercomp}, given $f : c \rightarrow \probstar$ and $i \in a$, there are three cases.
  \begin{enumerate}
    \item If $(t \compsimple s)(i) \in c$, then $\transformer{t\compsimple s}(f)$ sends it to $f(s(t(i)))$, and $t(i)$ must belong to $b$ and $s(t(i))$ to $c$, hence $\transformer{t}(\transformer{s}(f))$ sends $i$ to $f(s(t(i)))$.
    \item If $t(i) \in b$ but $s(t(i)) \notin c$, then $(t\compsimple s)(i) = s(t(i)) \in \probstar$. By definition, $\transformer{t\compsimple s}(f)(i) = s(t(i))$. On the right-hand side, since $s(t(i)) \notin c$, we have $\transformer{s}(f)(t(i)) = s(t(i)) \in \probstar$. Consequently, $\transformer{t}( \transformer{s}(f))(i) =  \transformer{s}(f)(t(i)) = s(t(i))$.
    \item If $t(i) \notin b$, then $t(i) \in \probstar$ and $(t\compsimple s)(i) = t(i)$. Thus $\transformer{t\compsimple s}(f)(i) = t(i)$. The definition of $\transformer{t}$ then gives $\transformer{t}(\transformer{s}(f))(i) = t(i)$.
  \end{enumerate}

  For \cref{eqn:transformerpar}, recall that we can identify $\probstar^{\!\!\!b+b'}$ with $\probstar^{\!\!\!b} \times \probstar^{\!\!\!b'}$ and $\probstar^{\!\!\!a+a'}$ with $\probstar^{\!\!\!a} \times \probstar^{\!\!\!a'}$. Let $h \in \probstar^{\!\!\!b+b'}$ correspond to the pair $(f: b \rightarrow
    \probstar, g: b' \rightarrow \probstar)$. Let $i \in a+a'$. If $i \in a$, then $(t+t')(i) = t(i)$ (embedded in $(b+b')+\probstar$). Then, we have two cases
  \begin{enumerate}
    \item If $t(i) \in b$, then $\transformer{t+t'}(h)(i) = h(t(i)) = f(t(i))$. On the other side, the left projection yields $\transformer{t}(f)(i) = f(t(i))$.
    \item If $t(i) \in \probstar$, then $\transformer{t+t'}(h)(i) = t(i)$ and $\transformer{t}(f)(i) = t(i)$.
  \end{enumerate}
  The case where $i \in a'$ is symmetric, involving $t'$ and $g$.  We see that the action on $h$ behaves as the component-wise action on $f$ and $g$, and we conclude that \cref{eqn:transformerpar} holds.

  It remains to prove that the symmetries are preserved. The symmetries in $\BInt$ are the symmetries in $(\FSet,+,\emptyset)$. In particular, they never have an image in $\probstar$, so their predicate transformers are simply precompositions: $\transformer{\swap{a,b}} = -\circ \swap{a,b}$. Identifying $\probstar^{\!\!\!a+b}$ with $\probstar^{\!\!\!a} \times \probstar^{\!\!\!b}$, precomposition by $\swap{a,b}$ is precisely the symmetry for the cartesian product $(f,g) \mapsto (g,f)$.
\end{proof}

\subsection{The Graded Category \texorpdfstring{$\EPBS$}{BIntStoch}}\label{sec:EPara}
We can now define an $\op{\BInt}$-"graded prop" $\EPBS$ (with a single colour) that will be the place to interpret the graded string diagrams of the next section. We progressively introduce the structure and show it satisfies the required properties in \cref{fig:axiomssyncat}. While we take the concrete route here, we note that $\EPBS$ can be constructed abstractly with a graded version of the external "para construction" defined in~\cite[Definition 4.2]{Smithe2022b} (also called category by proxy in~\cite[Definition 3.1]{Smithe2022})---just like $\BImP$ was constructed with the graded version of the "para construction". This justifies the name $\EPBS$ for \underline{e}xternal \underline{p}arametrisation of $\underline{\mathbf{B}}\mathbf{in}\underline{\mathbf{S}}\mathbf{toch}$.

Recall that an $a$-graded morphism $n \rightarrow m$ in $\BImP$ can be seen as a family of $2^a$ stochastic matrices indexed by the possible outcomes of $a$ nondeterministic binary choices. Very similarly, an $a$-graded morphism $n \rightarrow m$ in $\EPBS$ will be a family of matrices now indexed by $\probstar^{\!\!\!a}$ for $a$ possible interventions.

\AP The objects of $\intro*\EPBS$ are natural numbers $n \in \N$. The set of $a$-graded morphisms from $n$ to $m$, $\EPBS_a(n,m)$, is the set of functions of type $\probstar^{\!\!\!a} \rightarrow \BStoch(n,m)$. Intuitively, $f: n \garrow{a} m$ is a choice of $2^m \times 2^n$ stochastic matrix for each possible intervention (either replace with $\flip{p}$ or do nothing) on $a$ available places. 

For every morphism between grades $t: a \rightarrow b \in \BInt$ and graded morphism $f: n \garrow{a} m \in \EPBS$, the "regrading" $t \regrade f$ is defined as follows.
\begin{equation}\label{eqn:regradeBCaus}
  t \regrade f \coloneq \probstar^{\!\!\!b} \xrightarrow{\transformer{t}} \probstar^{\!\!\!a} \xrightarrow{f} \BStoch(n,m)
\end{equation}
Both \cref{eqn:th:trivialregrade} and \cref{eqn:th:seqregrade} hold because $\transformer{(\phantom{t})}$ is a functor (\cref{lem:transformerfunctor}). Indeed, "regrading" by the identity does nothing because $\transformer{\id} = \id$, and "regrading" by a composition $t\compsimple s$ is the same as "regrading" by $t$ then $s$ because of \cref{eqn:transformercomp}.

The composition of two graded morphisms $f: n \garrow{a} m$ and $g: m \garrow{b} \ell$ is given by
\begin{equation}\label{eqn:compBCaus}
  f \comp g \coloneq \probstar^{\!\!\!a+b} \cong \probstar^{\!\!\!a} \times \probstar^{\!\!\!b} \xrightarrow{f \times g} \BStoch(n,m) \times \BStoch(m,\ell) \xrightarrow{\compsimple} \BStoch(n,\ell).
\end{equation}
In words, we see the $a+b$ choices as a pair of $a$ choices and $b$ choices, and we multiply the corresponding matrices given by $f$ on the first $a$ choices and $g$ on the last $b$ choices. Thus, associativity \cref{eqn:th:assocseq} holds by associativity of matrix multiplication. For \cref{eqn:th:parregrade}, we see that "regrading" by $s+t$ becomes precomposition with $\transformer{s}\times \transformer{t}$ by \cref{eqn:transformerpar}, hence \[(f \times g) \circ (\transformer{s}\times \transformer{t}) = (f \circ \transformer{s}) \times (g \circ \transformer{t}) = (s \regrade f) \times (t \regrade g),\]
and applying $\compsimple$ yields the desired $(s+t) \regrade (f\comp g) = (s \regrade f) \comp (t \regrade g)$.

The monoidal product is addition on objects, and on morphisms $f: n \garrow{a} m$ and $f' : n' \garrow{a'} m'$, $f \tensor f'$ is defined as
\begin{equation}\label{eqn:parBCaus}
  \probstar^{\!\!\!a+a'} \cong \probstar^{\!\!\!a} \times \probstar^{\!\!\!a'} \xrightarrow{f \times f'} \scalebox{0.92}{$\BStoch(n,m)$} \times \scalebox{0.92}{$\BStoch(n',m')$} \xrightarrow{\otimes} \scalebox{0.92}{$\BStoch(n+n',m+m')$}.
\end{equation}
This is very similar to composition in the sense that $f \tensor g$ is the choice-wise product of the family of matrices $f$ and $g$. This makes associativity \cref{eqn:th:assocpar} as well as interchange \cref{eqn:th:interchangeupto} obvious by the corresponding properties for composition and product in $\BStoch$ (up to a reordering of choices). Moreover, \cref{eqn:th:parregradepar} holds by the same reasoning as for \cref{eqn:th:parregrade} above.

The identity morphism $\id_n : n \garrow{0} n$ is the function picking the identity $2^n\times 2^n$ matrix $\probstar^0 = \mathbf{1} \rightarrow \BStoch(n,n)$. It is clear that composition with $\id_n$ and monoidal product with $\id_0$ do nothing, i.e.\@ \cref{eqn:th:neutralid} and \cref{eqn:th:neutralidzer} hold.

Finally, the symmetry $\swap{n,m}: n + m \garrow{0} m+n$ is the function picking the symmetry in $\BStoch$: $\probstar^0 = \mathbf{1} \rightarrow \BStoch(n+m,m+n)$. Hence, it is involutive because the symmetry in $\BStoch$ is, i.e.\@ \cref{eqn:th:swapswap} holds. It also ensures that \cref{eqn:th:overthewire} holds.

We note that the "underlying SMC" of $\EPBS$ obtained by taking $0$-graded morphisms (recall \cref{prop:underlyingcat}) is precisely $\BStoch$.

In the following section, we define graded string diagrams for stochastic processes with interventions. They will be interpreted in $\EPBS$.

\subsection{Diagrammatic Reasoning with Interventions}\label{sec:diagraminterventions}
We define an $\BInt$-"graded theory" and a functor from its "syntactic category@@GRD" to $\EPBS$. 

First, recall that the grading category $\BInt$ was given a presentation in \cref{sec:BInt}. Here, we will flip it around to get a presentation of $\op{\BInt}$. We also orient it vertically and colour it to match the format in the previous sections: $\op{\BInt}$ is presented by the generators
\[\sigv = \{  \begin{tikzpicture}
	\begin{pgfonlayer}{nodelayer}
		\node [style=delgrade] (0) at (0, -0.5) {};
		\node [style=none] (1) at (0, 0.5) {};
	\end{pgfonlayer}
	\begin{pgfonlayer}{edgelayer}
		\draw [style=green] (0) to (1.center);
	\end{pgfonlayer}
\end{tikzpicture}
, \begin{tikzpicture}
	\begin{pgfonlayer}{nodelayer}
		\node [style=delgrade] (0) at (0, 0.5) {};
		\node [style=none] (1) at (0, -0.5) {};
	\end{pgfonlayer}
	\begin{pgfonlayer}{edgelayer}
		\draw [style=green] (0) to (1.center);
	\end{pgfonlayer}
\end{tikzpicture}
  \} \cup \{ \begin{tikzpicture}
	\begin{pgfonlayer}{nodelayer}
		\node [style=flipup] (0) at (0, -0.5) {$p$};
		\node [style=none] (1) at (0, 0.5) {};
	\end{pgfonlayer}
	\begin{pgfonlayer}{edgelayer}
		\draw [style=green] (0) to (1.center);
	\end{pgfonlayer}
\end{tikzpicture}
 \mid p \in \probs \}\]
and the axioms $\inlinefig{0.7}{vbintclosedcounit} = \inlinefig{0.8}{idzervert}$ and $\inlinefig{0.7}{vbintclosedp} = \inlinefig{0.8}{idzervert}$ (for $p \in \probs$), which are the vertical counterparts of \cref{eqn:axBInt}.
Second, we use the presentation of $\BStoch$ with $\thCC$, whose generators we recall:
\[\BStoch \text{ is a "prop" generated by } \sig = \{\begin{tikzpicture}
	\begin{pgfonlayer}{nodelayer}
		\node [style=none] (0) at (-0.5, 0) {};
		\node [style=black] (1) at (0.5, 0) {};
	\end{pgfonlayer}
	\begin{pgfonlayer}{edgelayer}
		\draw [style=single] (0.center) to (1);
	\end{pgfonlayer}
\end{tikzpicture}
, \input{copygen.tikz}, \input{andgen.tikz}, \begin{tikzpicture}[circuit logic US]
	\begin{pgfonlayer}{nodelayer}
		\node [style=none] (13) at (-1, 0) {};
		\node [style=not] (14) at (0, 0) {};
		\node [style=none] (15) at (1, 0) {};
	\end{pgfonlayer}
	\begin{pgfonlayer}{edgelayer}
		\draw [style=single] (13.center) to (14);
		\draw [style=single] (14) to (15.center);
	\end{pgfonlayer}
\end{tikzpicture}
, \begin{tikzpicture}
	\begin{pgfonlayer}{nodelayer}
		\node [style=flip] (0) at (-0.75, 0) {$p$};
		\node [style=none] (1) at (0.65, 0) {};
	\end{pgfonlayer}
	\begin{pgfonlayer}{edgelayer}
		\draw [style=single] (0) to (1.center);
	\end{pgfonlayer}
\end{tikzpicture}
\}.\]

\AP The "graded monoidal theory" $\intro*\thI$ combines these generators and the axioms in $\thCC$ with one additional generator and three additional axioms \cref{eqn:intergencounit,eqn:intergenp,eqn:intergenunit} listed below.\\
\setcounter{theqn}{0}
\begin{minipage}{0.28\textwidth}
  \begin{gather*}
    \input{intergen.tikz}: 1 \garrow{1} 1
  \end{gather*}
\end{minipage}
\begin{minipage}{0.35\textwidth}
  \begin{gather*}\label{eqn:intergencounit}
    \inlinefig{1}{intergencounit} \stepcounter{theqn}\tag{I\arabic{theqn}}\\
    \label{eqn:intergenp}
    \inlinefig{1}{intergenp} \stepcounter{theqn}\tag{I\arabic{theqn}}
  \end{gather*}
\end{minipage}
\begin{minipage}{0.35\textwidth}
  \begin{gather*}
    \label{eqn:intergenunit}
    \inlinefig{1}{intergenunit} \stepcounter{theqn}\tag{I\arabic{theqn}}
  \end{gather*}
\end{minipage}\\[0.2em]
With \cref{eqn:intergencounit,eqn:intergenp}, we see that "regrading" becomes the representation of an intervention internal to the diagrammatic theory. Indeed, the new generator indicates where there are points of possible interventions in a diagram. Also, "regrading" at an intervention point by $\inlinefig{0.7}{vgenbintcounit}$ removes the point (that is, the `do nothing' intervention), and "regrading" by $\inlinefig{0.7}{vgenbintp}$ replaces the point with $\flip{p}$ after discarding the input (that is, a \textit{cutting} or \textit{breaking} intervention). The last axiom \cref{eqn:intergenunit} intuitively reads as ``if you discard the outcome after an intervention, you can discard the input and the intervention''.

Let $\sync{\thI}$ denote the "graded prop" generated by this "theory@@GRD" as in \cref{defn:syntacticcat}. We will now interpret the diagrams in this "graded theory" by defining a "(graded symmetric monoidal) functor@graded functor" $\sem{-}: \sync{\thI} \rightarrow  \EPBS$. Since the $0$-graded morphisms in $\EPBS$ are simply stochastic matrices, we can let $\sem{-}$ act on the generators of $\thCC$ by sending them to their corresponding morphism in $\BStoch$. We automatically get that, after extending $\sem{-}$ to other $0$-"graded terms", the axioms in $\thI$ coming from $\thCC$ are satisfied.

Next, we define the value of $\sem{-}$ on the new generator as follows:
\[\sem{\!\inlinefig{0.7}{intergen}} : \probstar \rightarrow \BStoch(1,1) = \left\{
  \b* \mapsto \scalebox{0.75}{$\begin{bmatrix}
        1 & 0 \\
        0 & 1
      \end{bmatrix}$}\quad
  p \mapsto \scalebox{0.75}{$\begin{bmatrix}
        p   & p   \\
        1-p & 1-p
      \end{bmatrix}$}\right\}.
\]
This definition is the only one that ensures axioms \cref{eqn:intergencounit,eqn:intergenp} are satisfied. Indeed, "regrading" $\sem{\!\inlinefig{0.7}{intergen}}$ by $\inlinefig{0.7}{vgenbintcounit}$ or $\inlinefig{0.7}{vgenbintp}$ yields the matrices which are the image of $\b*$ and $p$ respectively, so they should coincide with $\sem{\!\inlinefig{0.7}{lonewire}}$ and $\sem{\!\inlinefig{0.7}{deletep}}$, respectively. The last axiom \cref{eqn:intergenunit} would be satisfied by any other definition.

Consequently, we can extend $\sem{-}$ to all graded terms to obtain an interpretation of these diagrams as stochastic processes with intervention points. Let us revisit the smoking scenario from~\cite{JacobsKZ21} with graded string diagrams.
\begin{exa}\label{exmp:gradedsmoking}
  Below on the left, we redraw the string diagram from \cref{fig:smokingscenario}, but we indicate two possible points of interventions on $S$ and $T$. On the right, we intervene on $S$ by "regrading" with $\inlinefig{0.7}{vgenbinthalf}$ and on $T$ by "regrading" with $\inlinefig{0.7}{vgenbintcounit}$.
  \[\scalebox{0.8}{\input{gradedsmoking.tikz}} \qquad \scalebox{0.8}{\input{gradedsmokingafterint.tikz}}\]
  Applying \cref{eqn:intergenp,eqn:intergencounit} to the diagram on the right yields \cref{eqn:smokingafterintervention}. The intervention points are chosen to match~\cite{JacobsKZ21}: the intervention on $S$ is the one performed in \cref{eqn:smokingafterintervention}, while the point on $T$ illustrates that an annotated point need not be used, as the `do nothing' "regrading" by $\inlinefig{0.7}{vgenbintcounit}$ simply removes it. Any wire of a diagram can be annotated in this way, and interventions at different points are performed independently of each other.
\end{exa}

Unlike the axiomatisation of $\BImP$ with string diagrams (\cref{sec:bimp}), the functor $\sem{-}:\sync{\thI} \rightarrow \EPBS$ is not full: there exist morphisms in $\EPBS$ that are not the semantics of any diagram. In particular, $\sem{-}$ is not an isomorphism of categories. We leave open whether $\sem{-}$ is faithful, equivalently, whether the axioms in $\thI$ are complete with respect to this semantics.

Intuitively, the matrices picked by a morphism in $\EPBS$ for each choice of intervention are entirely independent. For example, a morphism $f:0\garrow{1}1$ is a family $(f(\iota))_{\iota\in\probstar}$ of distributions on $\B$ with no coherence required between $f(\b*)$ (no intervention) and $f(p)$ for $p\in\probs$. One can have $f(\b*)=\ket{\b1}$ while $f(p)=\ket{\b0}$ for all $p\in\probs$. However, thinking of the grade as an intervention point, it does not make sense for a process to produce the true value $\b1 \in \B$ without any intervention, and to produce \textit{only} the false value $\b0 \in \B$ when intervened on. An intervention does not change the nature of the process, it only replaces some information flowing within it.\footnote{\cref{prop:nonuniversal} in \cref{appendix:interventions} shows that such arrows are not in the image of $\sem{-}$.}

While we could, in principle, define a semantic domain that corresponds to the subcategory of $\EPBS$ generated by the image of $\sem{-}$, we chose to use $\EPBS$ for convenience. Its morphisms are simple to describe and proving that it is a "graded category" closely followed the proof that the "para construction" yields a "graded category".

\begin{rem}\label{rem:nomodularity}
  The modularity result of \cref{sec:modularcompleteness} does not apply to $\EPBS$, for two reasons. First, the unit $0$ is not strictly terminal in the grading category $\op{\BInt}$: a morphism $0 \rightarrow a$ in $\op{\BInt}$ is a morphism $a \rightarrow 0$ in $\BInt$, i.e.\@ an arbitrary function $\fset{a} \rightarrow \probstar$, and there are many of those. Second, $\EPBS$ is not of the form $\Para{\G}{\catC}$ for an action of $\op{\BInt}$ on $\BStoch$: an $a$-graded morphism of $\EPBS$ is a family of stochastic matrices indexed by the uncountable set $\probstar^{\!\!\!a}$, which cannot be encoded as a single stochastic matrix with finitely many additional inputs. For these reasons, \cref{thm:gradedpresfrompres} cannot be invoked to obtain completeness of $\thI$ with respect to $\EPBS$, which we leave as an open question. The results of \cref{sec:gradedtheories}, on the other hand, do apply: $\sync{\thI}$ is a "graded prop", and $\sem{-}$ is a "graded functor" by \cref{lem:soundness}, so that equational reasoning in $\thI$ is sound with respect to $\EPBS$.
\end{rem}

\subsection{Identifiability of Single Interventions}\label{sec:identifiability}
The interest of interventional distributions lies in the fact that, under suitable assumptions on the causal structure, they can be computed from observational ones: in the smoking scenario, this is what allows the scientist to estimate the effect of smoking on cancer without running a trial. In the terminology of causal inference~\cite{Pearl2009}, the effect of the intervention is then \emph{identifiable}. The main result of~\cite{JacobsKZ21} in this sense is \cite[Theorem~8.1]{JacobsKZ21}: it isolates a shape of causal structures for which the distribution resulting from an intervention on a single variable is computable from the distribution without any intervention, by means of a generalised form of string diagram disintegration (\textit{cf.}\@ \cite{ChoJacobs2019}) called \emph{comb disintegration}. We now reformulate this result with graded string diagrams. This is conceptually simple: in our syntax, the observational model and its interventions are "regradings" of a single graded diagram $\omega$, so that the theorem becomes a statement about the single morphism $\sem{\omega}$ of $\EPBS$---namely, that its component at the trivial intervention (do nothing) determines all the others. Compared to the original result, the graded formulation simplifies some of the hypotheses and may be stated without assuming a specific value of the intervention. We refer to \cref{rem:combjkz} below for a more detailed comparison with~\cite{JacobsKZ21}.

\begin{thm}\label{thm:gradedidentifiability}
  Let $\omega : 0 \garrow{1} A \otimes X \otimes B \otimes C$ be a diagram in $\sync{\thI}$ of the form
  \begin{equation}\label{eqn:combshape}
    \omega \quad \coloneq \quad \vcenter{\hbox{\scalebox{0.8}{\input{combshape.tikz}}}}
  \end{equation}
  where $f_1 : 0 \rightarrow A \otimes X \otimes D$, $g : A \otimes X \rightarrow B$, and $f_2 : B \otimes D \rightarrow C$ are diagrams in $\Circ$, that is, they contain no intervention point. Let $\omega_{\b*}$ and $\omega_p$ (for $p \in \probs = [0,1]$) be the $0$-graded diagrams obtained by "regrading" $\omega$ with $\vcenter{\hbox{{\scalebox{0.7}{\begin{tikzpicture}
	\begin{pgfonlayer}{nodelayer}
		\node [style=delgrade] (0) at (0, -0.5) {};
		\node [style=none] (1) at (0, 0.5) {};
	\end{pgfonlayer}
	\begin{pgfonlayer}{edgelayer}
		\draw [style=green] (0) to (1.center);
	\end{pgfonlayer}
\end{tikzpicture}
}}}}$ and with $\vcenter{\hbox{{\scalebox{0.7}{\begin{tikzpicture}
	\begin{pgfonlayer}{nodelayer}
		\node [style=flipup] (0) at (0, -0.5) {$p$};
		\node [style=none] (1) at (0, 0.5) {};
	\end{pgfonlayer}
	\begin{pgfonlayer}{edgelayer}
		\draw [style=green] (0) to (1.center);
	\end{pgfonlayer}
\end{tikzpicture}
}}}}$ respectively:
  \begin{align}
    \omega_{\b*} \;\coloneq\; \vcenter{\hbox{\scalebox{0.75}{\input{combshaperegcounit.tikz}}}}\qquad \qquad
    \omega_p \;\coloneq\; \vcenter{\hbox{\scalebox{0.75}{\input{combshapereg.tikz}}}}  \label{eqn:combshapes}
  \end{align}
  so that $\sem{\omega_{\b*}} = \sem{\omega}(\b*)$ is the observational distribution of the model and $\sem{\omega_p} = \sem{\omega}(p)$ is the distribution obtained by intervening on $X$ with $\flip{p}$.

  If $\sem{\omega_{\b*}}$ has \emph{full support} (all the entries of this column vector are strictly positive), then every interventional distribution $\sem{\omega_p}$ can be computed from the observational distribution $\sem{\omega_{\b*}}$ alone, with no knowledge of $f_1$, $g$, $f_2$, and $D$. In other words, the graded morphism $\sem{\omega} \in \EPBS_1(0,A \otimes X \otimes B \otimes C)$ is determined by its component $\sem{\omega}(\b*)$.
\end{thm}

Before giving the proof, let us read \cref{eqn:combshape} as a causal model. The process $f_1$ produces the observed variables $A$ and $X$ together with the hidden ones $D$, which are not an output of $\omega$ but feed $f_2$. The variable $B$ depends on $A$ and $X$ through $g$ and $C$ depends on $B$ and $D$ through $f_2$. Both $B$ and $C$ are observable. Hence $X$ influences $C$ only through $B$, while $D$ is a hidden \emph{confounder}: it may be correlated with $X$, being produced jointly with it, and it influences $C$ directly (in the smoking scenario, $H$ is a common cause of $S$ and $C$). This is why the conditional distribution of $C$ given $X$ in $\sem{\omega_{\b*}}$ does not measure the effect of $X$ on $C$, whereas it does in $\sem{\omega_p}$: the "regrading" cuts any influence of $A$ and $D$ on $X$ (see \cref{eqn:combshapeint} below). The theorem states that $\sem{\omega_p}$ can nonetheless be computed from $\sem{\omega_{\b*}}$ alone, by the formula \cref{eqn:combformula} below. This is a feature of the shape \cref{eqn:combshape}: as discussed before \cref{rem:nomodularity}, a generic graded morphism of $\EPBS$ is an arbitrary family of states, with no relation whatsoever between its components.

\begin{proof}[Proof of \cref{thm:gradedidentifiability}]
  We first fix some notation. Recall that $\sem{\omega_{\b*}} \in \BStoch(0,A \otimes X \otimes B \otimes C)$. We write $P(\vvu x \vvv \vvw)$ for the weight of $\sem{\omega_{\b*}}$ at $\ket{\vvu x \vvv \vvw}$, where $\vvu \in \B^A$, $x \in \B$, $\vvv \in \B^B$, and $\vvw \in \B^C$ are the values carried by the output wires. We write $P(\vvu)$, $P(x \mid \vvu)$, $P(\vvv \mid \vvu x)$, and $P(\vvw \mid \vvu x \vvv)$ for the marginal and conditional probabilities computed from $P$ in the usual way, e.g.\@ $P(\vvu) = \sum_{x,\vvv,\vvw} P(\vvu x \vvv \vvw)$ and $P(x \mid \vvu) = \sum_{\vvv, \vvw} P(\vvu x \vvv \vvw)/P(\vvu)$. The assumption of full support ensures that no division by zero occurs. These probabilities are the entries of the following stochastic matrices in $\BStoch$:
  \begin{align*}
    \dist_{\!A}: 0 \rightarrow A = \textstyle\sum_{\vvu} P(\vvu)\ket{\vvu}\qquad & \qquad
    \dist_{X|A} : A \rightarrow X = \ket{\vvu} \mapsto  \textstyle\sum_x P(x\mid \vvu)\ket{x}                                                                         \\
    \dist_{B|AX} : A \otimes X \rightarrow B                                     & = \ket{\vvu x} \mapsto \textstyle\sum_{\vvv} P(\vvv\mid \vvu x)\ket{\vvv}          \\
    \dist_{C|AXB} : B \otimes X \otimes A \rightarrow  C                         & =\ket{\vvv x \vvu} \mapsto \textstyle\sum_{\vvw} P(\vvw\mid \vvu x \vvv)\ket{\vvw}
  \end{align*}
  (The inputs of $\dist_{C|AXB}$ are ordered as in \cref{eqn:combrecomp} below.) In the terminology of~\cite{ChoJacobs2019}, these matrices are the \emph{disintegrations} of $\sem{\omega_{\b*}}$ along its successive outputs, which are unique because of full support, and the chain rule $P(\vvu x \vvv \vvw) = P(\vvu)P(x\mid\vvu)P(\vvv\mid \vvu x)P(\vvw \mid \vvu x \vvv)$ expresses $\sem{\omega_{\b*}}$ as their composite, with suitable copying of $A$, $X$, and $B$. Finally, for $p \in \probs$, we write $\dist_p \coloneq \sem{\flip{p}} = p\ket{\b1} + (1-p)\ket{\b0}$.

  We shall prove that, for every $p \in \probs$,
  \begin{equation}\label{eqn:combrecomp}
    \sem{\omega_p} \;=\; \sem{\vcenter{\hbox{\scalebox{0.75}{\input{combshapereg.tikz}}}}} \;=\; \vcenter{\hbox{\scalebox{0.75}{\input{combrecomp.tikz}}}}
  \end{equation}
  where the first equality just recalls the definition of $\omega_p$, and the right-hand side is a composite of stochastic matrices ($\dist_p$ and the disintegrations of $\sem{\omega_{\b*}}$). The corresponding formula is
  \begin{equation}\label{eqn:combformula}
    \sem{\omega_p} = \sum_{\vvu, x, \vvv, \vvw} P(\vvu)\dist_p(x)P(\vvv \mid \vvu x)(\sum_{x' \in \B} P(x' \mid \vvu)P(\vvw \mid \vvu x' \vvv))\ket{\vvu x \vvv \vvw}.
  \end{equation}
  Since the right-hand side of \cref{eqn:combformula} only involves $p$ and the entries of $\sem{\omega_{\b*}}$, this proves the theorem. In essence, the proof is analogous to the one given in~\cite{JacobsKZ21}, preceded by one additional step accounting for the grading.

  \emph{Step 1: regrading.} The right-hand side of \cref{eqn:combformula} is built from the entries of $\sem{\omega_{\b*}}$, while its left-hand side is $\sem{\omega_p}$. To compare the two, we express both in terms of the components $f_1$, $g$, and $f_2$ of $\omega$, following \eqref{eqn:combshapes} and applying axioms \cref{eqn:intergencounit,eqn:intergenp}:
  \begin{align}
    \omega_{\b*} \;\overset{\text{\cref{eqn:intergencounit}}}{=}\; \vcenter{\hbox{\scalebox{0.75}{\input{combshapeerased.tikz}}}} \qquad \qquad
    \omega_p \;\overset{\text{\cref{eqn:intergenp}}}{=}\; \vcenter{\hbox{\scalebox{0.75}{\input{combshapeint.tikz}}}} \label{eqn:combshapeint}
  \end{align}
  These diagrams are $0$-graded, so they can be seen as diagrams in $\Circ$: they are composites of $f_1$, $g$, $f_2$, and generators of $\thCC$, with no intervention point left. We now compute the stochastic matrices they denote, $\sem{\omega_{\b*}}$ and $\sem{\omega_p}$.

  \emph{Step 2: comb disintegration.} Let $F_1 \coloneq \sem{f_1}$, $G \coloneq \sem{g}$, and $F_2 \coloneq \sem{f_2}$, with entries $F_1(\vvu x \vvd)$, $G(\vvv \mid \vvu x)$, and $F_2(\vvw \mid \vvv \vvd)$, where $\vvd \in \B^D$ ranges over the values carried by the hidden wire $D$. Reading the right-hand sides of \cref{eqn:combshapeint} as composites of stochastic matrices (composition is matrix multiplication, so the entries of a composite are obtained by summing over the values carried by the internal wires), we obtain, respectively,
  \begin{align}
    P(\vvu x \vvv \vvw)              & = \textstyle\sum_{\vvd} F_1(\vvu x \vvd) G(\vvv \mid \vvu x) F_2(\vvw \mid \vvv \vvd), \text{ and }\label{eqn:obsentries}       \\
    \sem{\omega_p}(\vvu x \vvv \vvw) & = \textstyle\sum_{x', \vvd} F_1(\vvu x' \vvd) \dist_p(x) G(\vvv \mid \vvu x)\, F_2(\vvw \mid \vvv \vvd), \label{eqn:intentries}
  \end{align}
  where, in \cref{eqn:intentries}, the value $x'$ output by $f_1$ on the wire $X$ is summed over because it is discarded, while $x$ is the fresh value produced by $\dist_p$.
  We now compute the disintegrations of $\sem{\omega_{\b*}}$ from \cref{eqn:obsentries}, that is, we express them through $F_1$, $G$, and $F_2$. Summing \cref{eqn:obsentries} over $\vvw$, then $\vvv$, then $x$, and using that $F_2$ and $G$ are stochastic, gives $P(\vvu x \vvv) = F_1(\vvu x) G(\vvv \mid \vvu x)$, $P(\vvu x) = F_1(\vvu x)$, and $P(\vvu) = \sum_x F_1(\vvu x)$, where $F_1(\vvu x) \coloneq \sum_{\vvd} F_1(\vvu x \vvd)$. Therefore,
  \begin{gather*}
    P(x \mid \vvu) = \frac{F_1(\vvu x)}{P(\vvu)}, \qquad P(\vvv \mid \vvu x) = G(\vvv \mid \vvu x), \qquad
    P(\vvw \mid \vvu x \vvv) = \textstyle\sum_{\vvd} \frac{F_1(\vvu x \vvd)}{F_1(\vvu x)}F_2(\vvw \mid \vvv \vvd),
  \end{gather*}
  where all denominators are positive by full support of $\sem{\omega_{\b*}}$ (for instance, $F_1(\vvu x) = P(\vvu x) \geq P(\vvu x \vvv \vvw) > 0$). Substituting in the right-hand side of \cref{eqn:combformula}, we get
  \begin{align*}
     & P(\vvu)\,\dist_p(x)\,P(\vvv \mid \vvu x)\sum_{x'} P(x' \mid \vvu)\,P(\vvw \mid \vvu x' \vvv)                                                                                    \\
     & \qquad= P(\vvu)\,\dist_p(x)\,G(\vvv \mid \vvu x)\sum_{x'}\left( \frac{F_1(\vvu x')}{P(\vvu)} \sum_{\vvd}\frac{F_1(\vvu x' \vvd)}{F_1(\vvu x')}\,F_2(\vvw \mid \vvv \vvd)\right) \\
     & \qquad= \dist_p(x)\,G(\vvv \mid \vvu x)\sum_{x',\vvd} F_1(\vvu x' \vvd)\,F_2(\vvw \mid \vvv \vvd),
  \end{align*}
  which is \cref{eqn:intentries}. This proves \cref{eqn:combformula} whose diagrammatic rendering is \cref{eqn:combrecomp}.
\end{proof}

\begin{rem}\label{rem:combjkz}
  We comment on the relationship with~\cite[Theorem~8.1]{JacobsKZ21} and its proof.
  \begin{enumerate}
    \item[(i)] The graded formulation brings two benefits over~\cite[Theorem~8.1]{JacobsKZ21}. The first concerns the hypotheses. In~\cite{JacobsKZ21}, the variable $X$ is the output of a designated generator $x$ of the syntax, and intervening on $X$ is a functor $\mathrm{cut}_X$ on diagrams, which replaces \emph{every} occurrence of $x$ by the uniform distribution on $X$, preceded by the discarding of the inputs of $x$. Accordingly, \cite[Theorem~8.1]{JacobsKZ21} assumes that the model has the shape \cref{eqn:combshape} with the intervention point replaced by $x$, and moreover that $x$ does not occur in $f_1$, $g$, and $f_2$, for otherwise $\mathrm{cut}_X$ would modify the model elsewhere than at the displayed position. In our setting, the intervention point marks the exact location where the intervention takes place, so that the second condition is not needed: the generator producing $X$ is simply part of $f_1$, and no condition on its other occurrences is required. The two interventional models coincide when \cite{JacobsKZ21} applies. Indeed, by naturality of discarding in $\Circ$ (axioms (D1)--(D3) of \cref{fig:binstocheqns}), discarding the output of $x$ amounts to discarding its inputs, so that the right-hand diagram of \cref{eqn:combshapeint} for $p = \sfrac{1}{2}$ is precisely $\mathrm{cut}_X$ applied to $\omega_{\b*}$. The second benefit concerns the intervention itself. In~\cite{JacobsKZ21}, it is always the cut by the uniform distribution. Interventions by regrading allow fixing $X$ to an arbitrary $\flip{p}$ with $p \in \probs$, and \cref{thm:gradedidentifiability} is proved once for all interventions. In this sense, the identifiability of the effect of $X$ is independent of the value forced on $X$.
    \item[(ii)] Step 2 of the proof is the comb disintegration of~\cite[Theorem 6.3]{JacobsKZ21}, which we have unfolded into a computation with the entries of $\sem{\omega_{\b*}}$ to keep the paper self-contained. In the terminology of that paper, the matrices $\dist_{\!A}$, $\dist_{X|A}$, and $\dist_{C|AXB}$ assemble, together with suitable copies of $A$ and $X$, into a \emph{2-comb} $B \rightarrow (A \otimes X) \otimes C$ with a `hole' of type $A \otimes X \rightarrow B$, and $\dist_{B|AX}$ is the stochastic matrix plugged into the hole. Plugging $\dist_{B|AX}$ directly gives back $\sem{\omega_{\b*}}$, while plugging it after the cut $\inlinefig{0.7}{deletep}$ on $X$ gives $\sem{\omega_p}$, exactly as in~\cite[Equation~(14)]{JacobsKZ21}. Uniqueness of comb disintegration is what makes the right-hand side of \cref{eqn:combrecomp} independent of the decomposition $(f_1,g,f_2)$ of $\omega$, which is not recoverable from $\sem{\omega_{\b*}}$: the same observational distribution is compatible with many choices of hidden wire $D$ and of $f_1$ and $f_2$.
    \item[(iii)] The shape \cref{eqn:combshape} covers the single-variable versions of the classical \emph{front-door} and \emph{back-door} criteria~\cite{Pearl2009}, see~\cite[Section 8]{JacobsKZ21} for a discussion. When $A$ is absent ($A = 0$), $B$ mediates the effect of $X$ on $C$ while $D$ is a hidden common cause of $X$ and $C$, and \cref{eqn:combformula} is the front-door adjustment formula. This is the smoking scenario of \cref{exmp:gradedsmoking}, with $X = S$, $B = T$, $C = C$, $D = H$, $g = t$, $f_2 = c$, and $f_1$ consisting of $h$ followed by a copy and $s$ (after the point on $T$ has been removed by "regrading" it with $\inlinefig{0.7}{vgenbintcounit}$, using \cref{eqn:intergencounit}). When $D$ is absent ($D = 0$), all confounding between $X$ and its effects goes through the observed $A$, and the marginal of \cref{eqn:combformula} on $X$ and $C$ is the back-door adjustment formula $\dist_p(x)\sum_{\vvu} P(\vvu) P(\vvw \mid \vvu x)$, with $A$ as the set of adjustment variables.
    \item[(iv)] The same argument applies to a diagram carrying $j$ intervention points at the same position of the shape \cref{eqn:combshape}, on a bundle $X$ of $j$ wires. For every $t \in \BInt(j,0)$, that is, for every choice of either doing nothing or intervening with some $\flip{p}$ at each of the points, $\sem{t \regrade \omega}$ is given by a formula analogous to \cref{eqn:combformula}, in which only the coordinates of $x$ that are intervened upon are replaced by fresh values. Hence $\sem{\omega} \in \EPBS_j(0,A \otimes X \otimes B \otimes C)$ is again determined by its component at $(\b*,\dots,\b*)$.
  \end{enumerate}
\end{rem}

\section{Conclusion}
In this work, we laid the foundations of graded monoidal algebra by introducing suitable notions of theory, syntactic "prop" generated by a theory, model, and string diagram. We proved a modularity result, which allows reducing the task of axiomatising a "graded prop" to that of axiomatising the base category and the grading category, provided the unit of the grading category is strictly terminal. As a first application, we axiomatised a "graded prop" whose morphisms represent imprecise probabilistic processes. We showed how a simple probabilistic programming language with nondeterminism may be interpreted in this category, and how this model may be extended with exact conditioning. As a second application, we showed how the presentation of the category of Gaussian maps in~\cite{Stein2026} can arise as the "flattening" of a "graded monoidal theory" obtained via our modularity result. As a third application, we used graded string diagrams to represent causal interventions in probabilistic processes.

Our theoretical investigations were first motivated by the study of imprecise probability through graded distributive Markov categories, as introduced by~\cite{LiellCock2025}. The diagrammatic language we introduced does not capture the full structure exploited in~\cite{LiellCock2025}. In particular, graded string diagrams do not naturally accommodate coproducts. This could be achieved by combining our approach with graphical languages for distributive categories~\cite{Comfort2020,Bonchi2023}, leading to a theory of graded distributive string diagrams.

Our novel approach to presenting $\Gauss$ naturally emerged from the resemblance between the proof of \cref{lem:decompositiongradedterm} and the normal form for terms in $\thGauss$ used in~\cite{Stein2026}. Similar normal forms underlie several other completeness results for monoidal presentations. For example, the presentation of $\Aff$ in~\cite{Bonchi2019c} relies on the normal form displayed in \cite[(17)]{Bonchi2019c} which pulls the generator $\inlinefig{0.75}{one}$ all the way to the left (bottom in loc.~cit.). Another example is the discard construction of~\cite{Carette2021} which uses a similar-looking purification result \cite[Lemma~3.5]{Carette2021}. The latter work also considers a quotient analogous to our "regrading collapse", but formulated in terms of isometries instead of the coisometries in our instantiation. We believe both these examples could be formulated within our framework following \cref{sec:gauss}. We leave the systematic investigation of further applications of this pattern to future work.

Our last application shares its motivation with the preliminary ideas of~\cite{Stein2025} for a compositional semantics of causal interventions, which were developed independently of ours. The "graded theory" $\thI$ is a natural starting point for giving semantics to a language of causal models with interventions, as it treats interventions internally rather than as program transformations, which are harder to handle compositionally.

Realising this vision would first require tools for automated equational reasoning with graded string diagrams, such as hypergraph rewriting systems akin to those developed for classical string diagrams in~\cite{RewriteTh1,RewriteTh2,RewriteTh3}. Looking further ahead, one could allow the outcome of a process to influence subsequent interventions by adding a cograding wire emerging from the top---or, more generally, a double-categorical structure---to enable information flow from the process category to the grading category.

We plan to consider other instances of our framework. First, the axiomatisation of the category $\mathbf{Gauss}$ in~\cite{Stein2026} could be combined with grading wires that model nondeterminism, provided we describe the appropriate action of $\op{\FInj}$ on $\mathbf{Gauss}$. Second, it is natural to ask whether our modular approach to axiomatisation may be applied to categories of algebras for graded (co)monads, as studied in programming language semantics~\cite{Katsumata2014,Gaboardi2016,Dorsch2019,Gaboardi2021,Ford2022,Forster2024}. Finally, we would like to study grading of string diagrammatic calculi for automata and other state-based systems, as introduced in~\cite{Piedeleu2023}, drawing inspiration from the use of graded monads in trace semantics~\cite{Milius2015}.

On a more abstract level, we are also interested in characterising the relation between our "graded monoidal theories" and graded algebraic theories (see e.g.\@ \cite{Kura2020,Sanada2023}) as the authors of~\cite{Bonchi2018} have done in the classical setting. This would allow obtaining graded monoidal presentations starting from graded algebraic presentations~\cite{Katsumata2022}.

\bibliographystyle{alphaurl}
\bibliography{refs}
\newpage
\appendix
\section{Appendix to \texorpdfstring{\cref{sec:background}}{Section 2}}\label{appendix:background}
\counterwithin{figure}{section}
\setcounter{figure}{0}

The additional coherence laws in a symmetric monoidal "actegory" are the following.
\begin{figure}[!ht]
  \begin{equation*}
    \begin{tikzcd}
      {A \act (B \act (X \otimes Y))} & {A \act (X \otimes (B \act Y))} & {X \otimes (A \act (B \act Y))} \\
      {(A\otimes B) \act (X \otimes Y)} && {X \otimes ((A \otimes B) \act Y)}
      \arrow["{A\act \kappa}", from=1-1, to=1-2]
      \arrow["\mu"', from=1-1, to=2-1]
      \arrow["\kappa", from=1-2, to=1-3]
      \arrow["{X\otimes \mu}", from=1-3, to=2-3]
      \arrow["\kappa"', from=2-1, to=2-3]
    \end{tikzcd}
  \end{equation*}
  \begin{equation*}
    \begin{tikzcd}
      {I \act (X \otimes Y)} && {X \otimes (I \act Y)} \\
      & {X\otimes Y}
      \arrow["\kappa", from=1-1, to=1-3]
      \arrow["\eta"', from=1-1, to=2-2]
      \arrow["{X\otimes \eta}", from=1-3, to=2-2]
    \end{tikzcd}
  \end{equation*}
  \begin{equation*}
    \begin{tikzcd}
      {A \act (X \otimes (Y \otimes Z))} & {X \otimes (A \act (Y \otimes Z))} & {X\otimes (Y \otimes (A \act Z))} \\
      {A \act ((X\otimes Y) \otimes Z)} && {(X \otimes Y)\otimes (A \act Z)}
      \arrow["\kappa", from=1-1, to=1-2]
      \arrow["{A \act \alpha}"', from=1-1, to=2-1]
      \arrow["{X \otimes \kappa}", from=1-2, to=1-3]
      \arrow["\alpha", from=1-3, to=2-3]
      \arrow["\kappa"', from=2-1, to=2-3]
    \end{tikzcd}
  \end{equation*}
  \caption{Coherence laws in symmetric monoidal "actegories".}\label{fig:coherencesymmact}
\end{figure}

\section{Appendix to \texorpdfstring{\cref{sec:gradedtheories}}{Section 3}}\label{appendix:gradedtheories}

We restate the definition of a "graded SMC" in full detail.

\begin{defi}[Graded SMC]
  A \textbf{$\G$-"graded symmetric strict monoidal category"} $\catC$ consists of the following data.
  \begin{itemize}
    \item A class of \textbf{objects} $\mathrm{Ob}(\catC)$. We write $X \in \catC$ instead of $X \in \mathrm{Ob}(\catC)$.
    \item For all objects $X,Y \in \catC$ and grade $\grade \in \G$, a set $\catC_{\grade}(X,Y)$ of \textbf{morphisms of grade $\grade$}. We often write $f: X \garrow{\grade} Y$ to mean $f \in \catC_{\grade}(X,Y)$.
    \item For all objects $X,Y \in \catC$ and morphism between grades $t: \grade \rightarrow \gradeb$, a function $t \regrade - : \catC_{\gradeb}(X,Y) \rightarrow \catC_{\grade}(X,Y)$ called \textbf{"regrading" by $t$}. We require that "regrading" by an identity morphism does nothing, and a composition of "regradings" is equal to "regrading" by a composition:
          \begin{equation}\label{eqn:gsmc:regradeisfunctor}\arraycolsep=1pt
            \forall s: \grade'' \rightarrow  \grade', t: \grade' \rightarrow \grade, f : X \garrow{\grade} Y, \qquad \begin{array}{rcl} \id_{\grade} \regrade f & = & f                           \\
             s \regrade (t \regrade f)     & = & (s\compsimple t) \regrade f
            \end{array}\ .
          \end{equation}
    \item For all objects $X,Y,Z \in \catC$ and grades $\grade,\gradeb \in \G$, there is a \textbf{composition} function ${\comp} : \catC_{\grade}(X,Y) \times \catC_{\gradeb}(Y,Z) \rightarrow \catC_{\grade\gradeb}(X,Z)$ that commutes with "regrading" in the following sense:
          \begin{equation}\label{eqn:gsmc:compisnatural}
            \forall\ \begin{matrix}
              t: \grade' \rightarrow  \grade & f : X \garrow{\grade} Y \\t': \gradeb' \rightarrow  \gradeb &
              g: Y \garrow{\gradeb} Z
            \end{matrix}, \quad (t \regrade f) \comp (t' \regrade g) = (tt') \regrade (f \comp g).
          \end{equation}
          We also require that composition is associative:\footnote{We implicitly use the fact that $\G$ is strict, so the grades of both morphisms in \cref{eqn:gsmc:compisassoc} coincide: $\grade(\grade'\grade'') = (\grade\grade')\grade''$.}
          \begin{equation}\label{eqn:gsmc:compisassoc}
            \forall f: X \garrow{\grade} Y, g: Y \garrow{\grade'} Z, h: Z \garrow{\grade''} W, \qquad f \comp (g \comp h) = (f \comp g) \comp h.
          \end{equation}
    \item For all objects $X \in \catC$, there is an \textbf{identity morphism} $\id_X \in \catC_{\unitG}(X,X)$ that is unital for composition:\footnote{Again, we implicitly need strictness of $\G$ to use $\unitG\grade = \grade = \grade \unitG$. Henceforth, we will not mention other similar uses of strictness.}
          \begin{equation}\label{eqn:gsmc:idcomp}
            \forall f: X \garrow{\grade} Y, \qquad \id_X \comp f = f = f \comp \id_Y.
          \end{equation}
    \item There is a binary operation $\tensor : \mathrm{Ob}(\catC) \times \mathrm{Ob}(\catC) \rightarrow \mathrm{Ob}(\catC)$ called the \textbf{monoidal product} that is associative $X\tensor (X' \tensor X'') = (X \tensor X') \tensor X''$ and has a unit $\unitC \in \catC$ ($\unitC \tensor X = X = X \tensor \unitC$).
    \item The monoidal product is also defined on morphisms as a $(\G\times \G)$-indexed family of functions $\tensor : \catC_{\grade}(X,Y) \times \catC_{\gradeb}(X',Y') \rightarrow \catC_{\grade\gradeb}(X\tensor X', Y \tensor Y')$ satisfying five properties. First, $\tensor$ is associative:    \begin{equation}\label{eqn:gsmc:assoctens}
            \forall \{f_i: X_i \garrow{\grade_i} Y_i\}_{i=1,2,3}, \qquad f_1 \tensor (f_2 \tensor f_3) = (f_1 \tensor f_2) \tensor f_3.
          \end{equation}
          Second, monoidal product with $\id_{\unitC} \in \catC_{\unitG}(\unitC,\unitC)$ does nothing:
          \begin{equation}\label{eqn:gsmc:idtens}
            \forall f : X \garrow{\grade} Y, \qquad \id_{\unitC} \tensor f = f = f\tensor \id_{\unitC}.
          \end{equation}
          Third, the monoidal product of two identity morphisms yields an identity morphism:
          \begin{equation}\label{eqn:gsmc:tensorofids}
            \forall X, X' \in \catC, \qquad \id_X \tensor \id_{X'} = \id_{X \tensor X'}.
          \end{equation}
          Fourth, monoidal product commutes with "regrading" in the following sense:
          \begin{equation}\label{eqn:gsmc:tensorisnatural}
            \forall\ \begin{matrix}
              t: \grade' \rightarrow  \grade & f : X \garrow{\grade} Y \\t': \gradeb' \rightarrow  \gradeb &
              f': X' \garrow{\gradeb} Y'
            \end{matrix}, \quad (t \regrade f) \tensor (t' \regrade f') = tt' \regrade (f \tensor f').
          \end{equation}
          Finally fifth, the \textbf{interchange law} holds "up to regrading" by a symmetry. Note that the interchange law cannot hold as usual because the monoidal product in $\G$ is not commutative, and this implies the grades of the morphisms in the following equation do not coincide:
          \[(f \tensor f') \comp (g \tensor g') \neq (f \comp g) \tensor (f' \comp g').\]
          The grade of the L.H.S.~can be obtained by applying a symmetry in $\G$ to the grade of the R.H.S., and we require the law to hold "up to regrading" by that symmetry. Formally,     \begin{gather}
            \forall f : X \garrow{\grade} Y,  g: Y \garrow{\gradeb} Z, f' : X' \garrow{\grade'} Y',  g': Y' \garrow{\gradeb'} Z',\notag\\
            (f \tensor f') \comp (g \tensor g') = \id_{\grade}\swapg{\gradeb\grade'}\id_{\gradeb'} \regrade ((f \comp g) \tensor (f' \comp g')).\label{eqn:gsmc:interchange}
          \end{gather}
    \item For all objects $X,X' \in \catC$, there is a morphism $\swap{X,X'}: X \tensor X' \garrow{\unitG} X' \tensor X$ of grade $\unitG$ called the \textbf{symmetry}. We require that the composite of two symmetries be the identity:
          \begin{equation}\label{eqn:gsmc:symmiso}
            \forall X,X' \in \catC, \qquad \swap{X,X'} \comp \swap{X'\!\!,X} = \id_{X \tensor X'},
          \end{equation}
          that symmetries are \textbf{coherent},
          \begin{equation}\label{eqn:gsmc:symmcoherent}
            \forall X,Y,Z \in \catC, \qquad \swap{X,Y\otimes Z} = (\swap{X,Y} \tensor \id_Z) \comp (\id_Y \tensor \swap{X,Z}),
          \end{equation}
          and the \textbf{sliding equation} holds, once again, up to a (necessary) "regrading" by a symmetry.    \begin{equation}\label{eqn:gsmc:sliding}
            \forall f: X \garrow{\grade} Y, g: X' \garrow{\grade'} Y', (f \tensor g) \comp \swap{Y\!\!,Y'} = \swapg{\grade\grade'} \regrade (\swap{X\!\!,X'} \comp (g \tensor f)).
          \end{equation}
  \end{itemize}
\end{defi}

When defining the closure $\closed{\axh}$ in \cref{defn:closure}, we require that it is a congruence. Explicitly, this amounts to asking that $\closed{\axh}$ is closed under the inference rules of \cref{fig:infrulessyncat}.
\begin{figure}[!ht]
  \begin{gather*}
    \begin{bprooftree}
      \AxiomC{$\phantom{\inlinefig{0.45}{gradedletterf}}$}
      \RightLabel{\hypertarget{rulerefl}{\scriptsize\textsc{Refl}}}
      \UnaryInfC{$\inlinefig{0.45}{gradedletterf} {\scriptstyle =} \inlinefig{0.45}{gradedletterf}$}
    \end{bprooftree} \quad \begin{bprooftree}
      \AxiomC{$\inlinefig{0.45}{gradedletterf} {\scriptstyle =} \inlinefig{0.45}{gradedletterg}$}
      \RightLabel{\hypertarget{rulesymm}{\scriptsize\textsc{Symm}}}
      \UnaryInfC{$\inlinefig{0.45}{gradedletterg} {\scriptstyle =} \inlinefig{0.45}{gradedletterf}$}
    \end{bprooftree} \quad \begin{bprooftree}
      \AxiomC{$\inlinefig{0.45}{gradedletterf} {\scriptstyle =} \inlinefig{0.45}{gradedletterg}$}
      \AxiomC{$\inlinefig{0.45}{gradedletterg} {\scriptstyle =} \inlinefig{0.45}{gradedletterh}$}
      \RightLabel{\hypertarget{ruletrans}{\scriptsize\textsc{Trans}}}
      \BinaryInfC{$\inlinefig{0.45}{gradedletterf} {\scriptstyle =} \inlinefig{0.45}{gradedletterh}$}
    \end{bprooftree}\\
    \begin{bprooftree}
      \AxiomC{$\inlinefig{0.55}{gradedletterf1} {\scriptstyle =} \inlinefig{0.55}{gradedletterf2}$}
      \AxiomC{$\inlinefig{0.55}{gradedletterg1} {\scriptstyle =} \inlinefig{0.55}{gradedletterg2}$}
      \RightLabel{\hypertarget{ruleseq}{\scriptsize\textsc{Seq}}}
      \BinaryInfC{$\inlinefig{0.55}{gradedletterf1g1} {\scriptstyle =} \inlinefig{0.55}{gradedletterf2g2}$}
    \end{bprooftree}\quad
    \begin{bprooftree}
      \AxiomC{$\inlinefig{0.55}{gradedletterf1} {\scriptstyle =} \inlinefig{0.55}{gradedletterf2}$}
      \AxiomC{$\inlinefig{0.55}{gradedletterg1} {\scriptstyle =} \inlinefig{0.55}{gradedletterg2}$}
      \RightLabel{\hypertarget{rulepar}{\scriptsize\textsc{Par}}}
      \BinaryInfC{$\inlinefig{0.55}{tensorf1g1} {\scriptstyle =} \inlinefig{0.55}{tensorf2g2}$}
    \end{bprooftree}\\
    \begin{bprooftree}
      \AxiomC{$\inlinefig{0.55}{gradedletterf} {\scriptstyle =} \inlinefig{0.55}{gradedletterg}$}
      \RightLabel{\hypertarget{rulereg}{\scriptsize\textsc{Reg}}}
      \UnaryInfC{$\inlinefig{0.55}{regradeft} {\scriptstyle =} \inlinefig{0.55}{regradegt}$}
    \end{bprooftree}
  \end{gather*}
  \caption{Rules to make $\closed{\axh}$ a congruence relative to $\comp$, $\tensor$, and $\regrade$.}\label{fig:infrulessyncat}
\end{figure}

We also mentioned that the axioms in \cref{fig:axiomssyncat} ensure that all the properties of "graded SMCs" \cref{eqn:gsmc:regradeisfunctor,eqn:gsmc:compisnatural,eqn:gsmc:compisassoc,eqn:gsmc:idcomp,eqn:gsmc:assoctens,eqn:gsmc:idtens,eqn:gsmc:tensorofids,eqn:gsmc:tensorisnatural,eqn:gsmc:interchange,eqn:gsmc:symmiso,eqn:gsmc:symmcoherent,eqn:gsmc:sliding} hold. In fact, each axiom is the diagrammatic rendering of the corresponding property as listed below:
\[\scalebox{0.9}{$\begin{array}{cccccccccccc}
        \cref{eqn:gsmc:regradeisfunctor}                      & \cref{eqn:gsmc:compisnatural} & \cref{eqn:gsmc:compisassoc} & \cref{eqn:gsmc:idcomp}  & \cref{eqn:gsmc:assoctens} & \cref{eqn:gsmc:idtens}     & \cref{eqn:gsmc:tensorofids} & \cref{eqn:gsmc:tensorisnatural} & \cref{eqn:gsmc:interchange}   & \cref{eqn:gsmc:symmiso} & \cref{eqn:gsmc:symmcoherent} & \cref{eqn:gsmc:sliding}   \\
        \Updownarrow                                          & \Updownarrow                  & \Updownarrow                & \Updownarrow            & \Updownarrow              & \Updownarrow               & \Updownarrow                & \Updownarrow                    & \Updownarrow                  & \Updownarrow            & \Updownarrow                 & \Updownarrow              \\
        \cref{eqn:th:trivialregrade},\cref{eqn:th:seqregrade} & \cref{eqn:th:parregrade}      & \cref{eqn:th:assocseq}      & \cref{eqn:th:neutralid} & \cref{eqn:th:assocpar}    & \cref{eqn:th:neutralidzer} & \cref{eqn:th:coherentid}    & \cref{eqn:th:parregradepar}     & \cref{eqn:th:interchangeupto} & \cref{eqn:th:swapswap}  & \cref{eqn:th:coherentswaps}  & \cref{eqn:th:overthewire}
      \end{array}$}
\]

\section{Appendix to \texorpdfstring{\cref{sec:application}}{Section 5}}\label{appendix:bimp}
In \cref{fig:binstocheqns}, we list all the axioms in the "theory@@MON" $\thCC$ as listed in~\cite{Piedeleu2025b}.
Note that any axiom which uses parameters $p,q,\dots$ implicitly quantifies over $\probs$, unless stated otherwise. We use the following syntactic sugar for the or and if gates.
\definecolor{comb}{HTML}{ddf1fb}
\definecolor{seq}{HTML}{6ae4a3}
\definecolor{neutral}{HTML}{f8f8f2}
\definecolor{unit}{HTML}{b9b9b9}
\newcommand{\thickness}{0.5mm}
\newcommand{\stringwidth}{0.4mm}
\newcommand{\mediumwidth}{7.5mm}
\newcommand{\smallwidth}{5mm}
\newcommand{\tinywidth}{5.75mm}
\newcommand{\boxsize}{\Large}
\newcommand{\corners}{}
\newcommand{\innersep}{1mm}
\newcommand{\dotsize}{3pt}
\tikzset{ihbase/.style={inner sep=0,circle,draw,fill=lightgray,minimum size={\dotsize}, node contents={}}}
\tikzset{ihblack/.style={ihbase,fill=black}}
\tikzset{ihwhite/.style={ihbase,fill=white}}
\tikzset{mat/.style={draw,fill=white,rectangle,node font=\scriptsize}}
\tikzset{ha/.style={mat,rounded rectangle,rounded rectangle left arc=none}}
\tikzset{haop/.style={mat,rounded rectangle,rounded rectangle right arc=none}}
\tikzset{blackha/.style={mat,rounded rectangle,rounded rectangle left arc=none,font=\color{white},fill=black}}
\tikzset{blackhaop/.style={mat,rounded rectangle,rounded rectangle right arc=none,font=\color{white},fill=black}}
\tikzset{anti/.style={inner sep=0,isosceles triangle,fill=black,draw=black, minimum width=0.75em, node contents={}}}
\tikzset{antiop/.style={anti,shape border rotate=180}}
\tikzset{antisq/.style={inner sep=0,rectangle,fill=black, minimum height=1em, minimum width=0.6em, node contents={}}}
\tikzset{count/.style={above,inner ysep=0.15em,font=\scriptsize}}
\tikzset{axiom/.style={above,font=\small}}
\tikzset{dir/.style={-Latex}}
\tikzset{st/.style={decoration={markings,
          mark={at position 0.5 with {\draw (0, 2pt) to (0, -2pt);}}},
      postaction=decorate}}
\tikzstyle{flip}=[draw={rgb,255: red,86; green,86; blue,86}, fill={rgb,255: red,86; green,86; blue,86},
rounded rectangle, rounded rectangle right arc=none, text=white, node font={\scriptsize}, inner sep =0.2em]
\tikzstyle{cflip}=[circle, draw={rgb,255: red,128; green,128; blue,128}, fill={rgb,255: red,128; green,128; blue,128}, inner sep=0pt, minimum size=4pt]
\tikzstyle{and}=[fill=white, draw=black, and gate]
\tikzstyle{or}=[fill=white, draw=black, or gate, anchor=center]
\tikzstyle{not}=[fill=white, draw=black, not gate]
\tikzstyle{xor}=[fill=white, draw=black, xor gate]
\tikzstyle{if}=[trapezium, trapezium angle=60, draw,inner xsep=0pt,outer sep=0pt,minimum height=10pt, rotate=270, text width=10pt, fill=white]

\tikzstyle{plain}=[inner sep=0pt]
\tikzstyle{black}=[circle, draw=black, fill=black, inner sep=0pt, minimum size={\dotsize}]
\tikzstyle{black-faded}=[circle, draw=light-gray, fill=light-gray, inner sep=0pt, minimum size=4pt]
\tikzstyle{white}=[circle, draw=black, fill=white, inner sep=0pt, minimum size={\dotsize}]
\tikzstyle{white-faded}=[circle, draw=light-gray, fill=white, inner sep=0pt, minimum size=4.5pt]
\tikzstyle{reg}=[draw, fill=white, rounded rectangle, rounded rectangle left arc=none, minimum height=1.2em, minimum width=1.4em, node font={\scriptsize}]
\tikzstyle{effect}=[draw, fill=white, rounded rectangle, rounded rectangle left arc=none, minimum height=1.2em, minimum width=1.4em]
\tikzstyle{coreg}=[draw, fill=white, rounded rectangle, rounded rectangle right arc=none, minimum height=1.2em, minimum width=1.4em, node font={\scriptsize}]
\tikzstyle{box}=[shape=rectangle, text height=1.5ex, text depth=0.25ex, yshift=0.2mm, fill=white, draw=black, minimum height=3mm, minimum width=5mm, font={\small}]
\tikzstyle{basic box}=[shape=rectangle, text height=1.5ex, text depth=0.25ex, yshift=0.2mm, fill={probcircuitcolor}, draw={probcircuitcolor}, minimum height=3mm, minimum width=5mm, text={probcircuitcolorop}, inner sep=4pt]
\tikzstyle{small box}=[draw, fill={probcircuitcolor}, rectangle, minimum height=1.2em, minimum width=1.4em, node font={\scriptsize}, text={probcircuitcolorop}]
\definecolor{probcircuitcolor}{rgb}{0.33,0.33,0.33}
\definecolor{probcircuitcolorop}{rgb}{1.0,1.0,1.0}
\definecolor{boolcircuitcolor}{rgb}{1.0,1.0,1.0}
\definecolor{boolcircuitcolorop}{rgb}{0.0,0.0,0.0}
\tikzstyle{prob circuit}=[shape=rectangle, text height=1.5ex, text depth=0.25ex, yshift=0.2mm, fill={probcircuitcolor}, draw={probcircuitcolor}, minimum height=3mm, minimum width=5mm, text={probcircuitcolorop}, inner sep=4pt]
\tikzstyle{small prob circuit}=[draw={probcircuitcolor}, fill={probcircuitcolor}, rectangle, minimum height=1.2em, minimum width=1.4em, node font={\scriptsize}, text={probcircuitcolorop}]
\tikzstyle{tall prob circuit}=[shape=rectangle, text height=1.5ex, text depth=0.25ex, yshift=0.2mm, fill={probcircuitcolor}, draw={probcircuitcolor}, minimum height=8mm, minimum width=5mm, text={probcircuitcolorop}]
\tikzstyle{bool circuit}=[shape=rectangle, text height=1.5ex, text depth=0.25ex, yshift=0.2mm, fill={boolcircuitcolor}, draw=black, minimum height=3mm, minimum width=5mm, font={\small}, text={boolcircuitcolorop}]
\tikzstyle{small bool circuit}=[draw, fill={boolcircuitcolor}, rectangle, minimum height=1.2em, minimum width=1.4em, node font={\scriptsize},text={boolcircuitcolorop}]
\tikzstyle{tall bool circuit}=[rectangle, draw, fill={boolcircuitcolor}, minimum height=8mm, minimum width=5mm ,text={boolcircuitcolorop}]

\newcommand{\wht}[1]{node[ihwhite,#1]}
\newcommand{\blk}[1]{node[ihblack,#1]}
\newcommand{\gry}[1]{node[ihbase,#1]}
\newcommand{\blkn}[1]{node[ihblack,#1]}
\newcommand{\blkm}[1]{node[ihblack,#1]}
\newcommand{\ha}{node[ha]}
\newcommand{\haop}{node[haop]}
\newcommand{\mat}{node[mat]}
\newcommand{\anti}{node[anti]}
\newcommand{\antiop}{node[antiop]}
\newcommand{\antisq}[1]{node[antisq,#1]}
\newcommand{\myeq}[1]{\mathrel{\overset{\makebox[0pt]{\mbox{\normalfont\tiny\sffamily #1}}}{=}}}

\newcommand{\genericcomult}[2]{
  \raisebox{-0.8em}{\scalebox{0.6}{\begin{tikzpicture}
        \node at (1, 0) [ihbase,solid,name=copy,#1];
        \draw[#2] (copy) .. controls (1.25, 0.5) .. (2, 0.5);
        \draw[#2] (0, 0) -- (copy);
        \draw[#2] (copy) .. controls (1.25, -0.5) .. (2, -0.5);
      \end{tikzpicture}}}
}
\newcommand{\genericcomultn}[2]{
  \tikz {
    \draw (1, 0) node[ihbase,name=copy,#1] .. controls (1.25, 0.5) .. (1.5, 0.5)
    -- node[count] {$#2$} (2.25, 0.5);
    \draw (0, 0) -- node[count] {$#2$} (copy) .. controls (1.25, -0.5) .. (1.5, -0.5)
    -- node[count] {$#2$} (2.25, -0.5);
  }
}
\newcommand{\genericcounit}[2]{
  \tikz \draw[#2] (0, 0) -- (1, 0) node[ihbase,#1, solid];
}
\newcommand{\genericcounitn}[2]{
  \tikz \draw (0, 0) -- node[count] {$#2$} (1, 0) node[ihbase,#1];
}
\newcommand{\genericmult}[2]{
  \tikz {
    \node at (1,0) (copy) [ihbase,#1,solid];
    \draw[#2] (0,  0.5) .. controls (0.75,  0.5) .. (copy);
    \draw[#2] (0, -0.5) .. controls (0.75, -0.5) .. (copy);
    \draw[#2] (copy) -- (2, 0);
  }
}
\newcommand{\genericmultn}[2]{
  \tikz {
    \draw (0,  0.5) -- node[count] {$#2$} (0.75,  0.5)
    .. controls (1,  0.5) .. (1.25, 0) node[ihbase,name=copy,#1];
    \draw (0, -0.5) -- node[count] {$#2$} (0.75, -0.5)
    .. controls (1, -0.5) .. (copy) -- node[count] {$#2$} (2.25, 0);
  }
}
\newcommand{\genericunit}[2]{
  \tikz \draw[#2] (0, 0) node[ihbase,#1, solid] -- (1, 0);
}
\newcommand{\genericunitn}[2]{
  \tikz \draw (0, 0) node[ihbase,#1, solid] -- node[count] {$#2$} (1, 0);
}

\newcommand{\Bcomult}{\genericcomult{ihblack}{}}
\newcommand{\Bcomultn}[1]{\genericcomultn{ihblack}{#1}}
\newcommand{\Bcounit}{\genericcounit{ihblack}{}}
\newcommand{\Bcounitn}[1]{\genericcounitn{ihblack}{#1}}
\newcommand{\Bmult}{\genericmult{ihblack}{}}
\newcommand{\Bmultn}[1]{\genericmultn{ihblack}{#1}}
\newcommand{\Bunit}{\genericunit{ihblack}{}}
\newcommand{\Bunitn}[1]{\genericunitn{ihblack}{#1}}

\newcommand{\Andgate}{
  \raisebox{-0.8em}{\scalebox{0.5}{
      \begin{tikzpicture}[circuit logic US]
        \begin{pgfonlayer}{nodelayer}
          \node [style=none] (91) at (-1.5, 0.5) {};
          \node [style=and] (93) at (0, 0) {};
          \node [style=none] (96) at (-1.5, -0.5) {};
          \node [style=none] (103) at (1.25, 0) {};
        \end{pgfonlayer}
        \begin{pgfonlayer}{edgelayer}
          \draw [in=0, out=150] (93) to (91.center);
          \draw [in=-150, out=0, looseness=0.75] (96.center) to (93);
          \draw (93) to (103.center);
        \end{pgfonlayer}
      \end{tikzpicture}}}
}
\[\input{causcirc/or.tikz} \coloneq \input{causcirc/or-def.tikz} \qquad \input{causcirc/if.tikz} \coloneq \input{causcirc/if-def.tikz}\]
\begin{figure}[p!]
  {\scriptsize
    \centering
    \begin{gather*}
      \input{causcirc/copy-associative.tikz} \;\myeq{A1} \input{causcirc/copy-associative-1.tikz} \qquad \input{causcirc/copy-unital-left.tikz} \myeq{A2l} \input{lonewire.tikz}\myeq{A2r} \input{causcirc/copy-unital-right.tikz}\qquad  \input{causcirc/copy-commutative.tikz} \myeq{A3} \input{causcirc/copy.tikz}
    \end{gather*}
    \noindent\rule{\linewidth}{0.4pt}
    \begin{gather*}
      \input{causcirc/and-associative-left.tikz}\myeq{B1}\input{causcirc/and-associative-right.tikz}\qquad \input{causcirc/and-unit-left.tikz}\myeq{B2l}\input{lonewire.tikz} \myeq{B2r}\input{causcirc/and-unit-right.tikz}\qquad \input{causcirc/swap-and.tikz}\myeq{B3}\input{causcirc/and.tikz}
    \end{gather*}
    \noindent\rule{\linewidth}{0.4pt}
    \begin{gather*}
      \input{causcirc/not-not.tikz} \myeq{B4} \input{lonewire.tikz}\qquad \qquad \qquad \qquad \input{causcirc/copy-and.tikz}\myeq{B5} \input{lonewire.tikz}\\
      \input{causcirc/copy-idxnot-and.tikz} \myeq{B6} \input{causcirc/del-false.tikz} \qquad \input{causcirc/and-or.tikz}\myeq{B7} \input{causcirc/orx2-and.tikz}
    \end{gather*}
    \noindent\rule{\linewidth}{0.4pt}
    \vspace{-3mm}
    \begin{gather*}
      \input{causcirc/0-copy.tikz}\myeq{C1} \input{causcirc/0x0.tikz}\qquad \qquad\qquad \qquad\qquad  \input{causcirc/1-copy.tikz}\myeq{C2} \input{causcirc/1x1.tikz} \\ \input{causcirc/and-copy.tikz}\myeq{C3} \input{causcirc/copyxcopy-andxand.tikz}\qquad \input{causcirc/not-copy.tikz}\myeq{C4} \input{causcirc/copy-not.tikz}
    \end{gather*}
    \vspace{-3mm}
    \noindent\rule{\linewidth}{0.4pt}
    \vspace{-3mm}
    \begin{gather*}
      \input{causcirc/and-del.tikz} \myeq{D1} \input{causcirc/delxdel.tikz} \qquad \input{causcirc/not-del.tikz} \myeq{D2}\input{causcirc/del.tikz} \qquad \input{causcirc/flip-del.tikz} \myeq{D3}\input{causcirc/empty-diag.tikz}
    \end{gather*}
    \noindent\rule{\linewidth}{0.4pt}
    \begin{gather*}
      \input{causcirc/flip-not.tikz}\myeq{E1} \;\input{causcirc/flip-1-p.tikz}
    \end{gather*}
    \begin{gather*}
      \scalebox{0.8}{\input{causcirc/ax-general-inverse.tikz}}\\       \text{where }  \tilde{r}=rp+(1-r)q \;\;\text{ and }\begin{array}{ll}
        \tilde{p} = \frac{rp}{\tilde{r}}       & \text{ if } \tilde{r}\neq 0, \text{ and anything otherwise}  \\
        \tilde{q} = \frac{r(1-p)}{1-\tilde{r}} & \text{ if }  \tilde{r}\neq 1, \text{ and anything otherwise}
      \end{array}
    \end{gather*}
    \begin{gather*}
      \input{causcirc/ax-bary.tikz}\\ \text{where } \tilde{p}=pq \;\;\text{ and } \tilde{q}=\frac{p(1-q)}{1-pq}\text{ for }pq\neq 1
    \end{gather*}
    \begin{gather*}
      \input{causcirc/ax-swap.tikz}
    \end{gather*}
  }
  \caption{Axioms for causal circuits.}\label{fig:binstocheqns}
\end{figure}

\section{Computing Inside \texorpdfstring{$\BImP$}{BImP}}\label{appendix:computing}
Both notations that we use for computations with morphisms in $\BStoch$, stochastic matrices and kets $\ket{-}$, are easily adapted to do computations inside $\BImP$.

First, a morphism $f: n \garrow{a} m$ in $\BImP$ was identified as a family of submatrices $f\ket{\vvx}: n \rightarrow m$ indexed by $\vvx \in \B^a$. Our interpretation of $f$ as an imprecise probabilistic process suggests that $\vvx$ can also be seen as a nondeterministic input to $f$ which determines the submatrix that will act on the rest of the input.

To further emphasise the distinction between the nondeterministic and probabilistic parts of the behaviour of $f$, we will write the input determining the submatrix with square kets $\nket{\vvx}$, and the rest of the input in the usual kets $\ket{\vvu}$. When tensored, we combine the notations to $\nket{\vvx}\ket{\vvu}$.

We can thus give the following equivalent definitions of "regrading" and sequential and parallel composition inside $\BImP$.
\begin{align*}
  (t\regrade g) \nket{\ttx_0\cdots \ttx_{b-1}}\ket{\vvu}  & = g \nket{\ttx_{t(0)}\cdots\ttx_{t(a-1)}}\ket{\vvu}      \\
  (f\comp g)\nket{\vvx}\nket{\vvy}\ket{\vvu}              & = g(\nket{\vvy}\otimes f\nket{\vvx}\ket{\vvu})           \\
  (f\tensor f')\nket{\vvx}\nket{\vvy}\ket{\vvu}\ket{\vvv} & = f\nket{\vvx}\ket{\vvu} \otimes f\nket{\vvy}\ket{\vvv}
\end{align*}

Second, representing a graded morphism $f: n \garrow{a} m$ in $\BImP$ as a single matrix divided into $2^a$ submatrices (with separators) is very convenient because we can compute sequential and parallel compositions methodically by slightly modifying the usual algorithms for multiplication and direct sum.

As we have said, the order of the submatrices cannot be ignored. Thankfully, the order of the columns does not change if we see a morphism $f: n \garrow{a} m$ as $2^a$ stochastic matrices or as a single stochastic matrix in $\BStoch(a+n, m)$. Indeed, the $j$th submatrix of $f$ is $f\nket{\vvx}\ket{-}$ with $\vvx = \mathtt{bin}(2^a-j)$, and the $i$th column of the $j$th submatrix of $f$ is $f\nket{\vvx}\ket{\vvu}$ with $\vvu= \mathtt{bin}(2^n-i)$. This coincides with the $((j-1)\cdot 2^n + i)$th column of $f$, i.e.\@ $f\nket{\vvx}\ket{\vvu}$.

This ordering of submatrices and the vertical line notation behave rather well with composition inside $\BImP$. For example, let $f: 0 \garrow{1} 1$ and $g: 1 \garrow{2} 1$ be defined below. 
\[f= \scalebox{0.8}{$\left[\!\! \begin{array}{c|c}
          1 & 0.5 \\
          0 & 0.5
        \end{array} \!\!\right]$} \qquad  g = \scalebox{0.8}{$\left[\!\! \begin{array}{cc|cc|cc|cc}
          1 & 1 & 1 & 0 & 0 & 1 & 0 & 0 \\
          0 & 0 & 0 & 1 & 1 & 0 & 1 & 1
        \end{array}\!\! \right]$}\]
To compute $f\comp g$, we proceed with a simple variant of the usual matrix multiplication algorithm. First focusing on the first column (or submatrix) of $f$, we multiply with each submatrix of $g$, and we add a vertical line between all the results, then we add another separator before doing the same thing for the second column of $f$. It yields the matrix
\[f\comp g = \scalebox{0.8}{$\left[\!\! \begin{array}{c|c|c|c|c|c|c|c}
          1 & 1 & 0 & 0 & 1 & 0.5 & 0.5 & 0 \\
          0 & 0 & 1 & 1 & 0 & 0.5 & 0.5 & 1
        \end{array} \!\!\right]$},\]
which has $2^3$ submatrices of dimension $2^1 \times 2^0$, and $f\comp g$ indeed belongs to $\BImP_3(0,1)$.

In general, given $f: n \garrow{a} m$ and $g: m \garrow{b} \ell$, for any $\vvx \in \B^a$ and $\vvy \in \B^b$, the submatrix $(f\comp g)\nket{\vvx}\nket{\vvy}\ket{-}$ is equal to the multiplication $(g\nket{\vvy}\ket{-})(f\nket{\vvx}\ket{-})$ (you can verify this by using the definition of composition in $\BImP$). Therefore, the procedure above is always sound.

Similarly, for parallel composition, given $f: n \garrow{a} m$ and $f' : n' \garrow{b} m'$, for any $\vvx \in \B^a$ and $\vvy \in \B^b$, the submatrix $(f\otimes f')\nket{\vvx}\nket{\vvy}\ket{-}$ is equal to the tensor $f\nket{\vvx}\ket{-} \otimes f'\nket{\vvy}\ket{-}$. Hence, the algorithm to write the matrix for $f \otimes f'$ mirrors that for sequential composition, except that multiplication is replaced with the Kronecker product. You first focus on the first submatrix of $f$, it gets tensored with each submatrix of $f'$ with vertical lines separating each result, then you move to the next submatrix of $f$, and so on. Below is an example.
\[f = \scalebox{0.8}{$\left[\!\!\begin{array}{cc|cc}
          1 & 0 & 0 & 1 \\0&1&1&0
        \end{array}\!\! \right]$} \quad f' = \scalebox{0.8}{$\left[\!\! \begin{array}{c|c}
          .5 & 0 \\.5&1
        \end{array} \!\!\right]$} \quad f\otimes f' = \scalebox{0.7}{$\left[ \!\!\begin{array}{cc|cc|cc|cc}
          .5 & 0  & 0 & 0 & 0  & .5 & 0 & 0 \\
          .5 & 0  & 1 & 0 & 0  & .5 & 0 & 1 \\
          0  & .5 & 0 & 0 & .5 & 0  & 0 & 0 \\
          0  & .5 & 0 & 1 & .5 & 0  & 1 & 0
        \end{array}\!\! \right]$}\]

\section{Appendix to \texorpdfstring{\cref{sec:gauss}}{Section 7}}
In \cref{fig:affeqns}, we list all the axioms of $\thAff$ as listed in~\cite{Stein2026}.
\begin{figure}
  \begin{gather*}
    \input{aff/assocadd.tikz} \qquad \input{aff/assoccom.tikz} \qquad \input{aff/assocunit.tikz} \qquad \input{aff/onecopy.tikz}
  \end{gather*}
  \begin{gather*}
    \input{aff/assoccopy.tikz} \qquad \input{aff/commcopy.tikz} \qquad \input{aff/unitcopy.tikz} \qquad \input{aff/onedelete.tikz}
  \end{gather*}
  \begin{gather*}
    \input{aff/addscalar.tikz} \qquad \input{aff/scalezero.tikz} \qquad \input{aff/copyscalar.tikz} \qquad \input{aff/deletescalar.tikz}
  \end{gather*}
  \begin{gather*}
    \input{aff/multscalar.tikz} \qquad \input{aff/plusscalar.tikz} \qquad \input{aff/zeroscalar.tikz} \qquad \input{aff/onescalar.tikz}
  \end{gather*}
  \begin{gather*}
    \input{aff/addcopy.tikz} \qquad \input{aff/zerocopy.tikz} \qquad \input{aff/zerodelete.tikz} \qquad \input{aff/adddelete.tikz}
  \end{gather*}
  \caption{Axioms for affine maps.}\label{fig:affeqns}
\end{figure}

\section{Appendix to \texorpdfstring{\cref{sec:interventions}}{Section 8}}\label{appendix:interventions}
\begin{proof}[Proof of \cref{prop:presBInt}]
  The assignment of \cref{prop:presBInt} extends to an identity-on-objects functor $F: \vync{\sigv,\axv} \rightarrow \BInt$: the axioms \cref{eqn:axBInt} hold in $\BInt$ because both sides of each axiom are interpreted as the unique function $\fset{0} \rightarrow \fset{0} + \probstar$. It remains to verify that $F$ is fully faithful.

  We first identify convenient representatives of morphisms in $\vync{\sigv,\axv}$. Every generator in $\sigv$ has at most one input and one output, so every connected component of a diagram in $\vynct{\sigv}$ is a path that starts at an input of the diagram or at $\inlinefig{0.75}{genbintunit}$, and ends at an output of the diagram, at $\inlinefig{0.75}{genbintcounit}$, or at $\bintp{p}$ for some $p \in \probs$. A component that starts at $\inlinefig{0.75}{genbintunit}$ and ends at $\inlinefig{0.75}{genbintcounit}$ or at $\bintp{p}$ is connected to no input and no output, and it is equal, via the laws of "SMCs", to the left-hand side of one of the axioms \cref{eqn:axBInt}. Hence, such components are equal, by \cref{eqn:axBInt}, to the identity on the monoidal unit $0$, and they can be ignored. We conclude that every morphism of $\vync{\sigv,\axv}$ except $\id_0$ is represented by a diagram in which every connected component is connected to an input or to an output. We call such diagrams \textbf{reduced}.

  Next, for every $a,b \in \N$, we partition the sets $\vync{\sigv,\axv}(a,b)$ and $\BInt(a,b)$ into sets $\vync{\sigv,\axv}^k(a,b)$ and $\BInt^k(a,b)$ indexed by $0 \leq k \leq a$.
  \begin{itemize}
    \item A morphism $t: a \rightarrow b$ belongs to $\vync{\sigv,\axv}^k(a,b)$ if and only if it has a reduced representative containing the generator $\bintp{p}$ for exactly $k$ distinct values of $p$. In a reduced diagram, each occurrence of $\bintp{p}$ is connected to an input, so $k$ is between $0$ and $a$, and $F(t)$ takes exactly $k$ distinct values in $\probs$; since $F$ is well-defined on morphisms of $\vync{\sigv,\axv}$, this shows that $k$ does not depend on the choice of reduced representative. Therefore, $\vync{\sigv,\axv}(a,b)$ is the pairwise disjoint union of all $\vync{\sigv,\axv}^k(a,b)$ with $0 \leq k \leq a$.
    \item A left-injective function $t : a \rightarrow b + \probstar$ belongs to $\BInt^k(a,b)$ if and only if it takes exactly $k$ distinct values in $\probs \subseteq \probstar$. It is clear that $k$ is uniquely determined and between $0$ and $a$ since $t$ takes at most $a$ distinct values. Therefore, $\BInt(a,b)$ is the pairwise disjoint union of all $\BInt^k(a,b)$ with $0\leq k \leq a$.
  \end{itemize}

  Next, we show that the functor $F$ induces a bijection between $\vync{\sigv,\axv}^k(a,b)$ and $\BInt^k(a,b)$ for all $a,b \in \N$ and $0\leq k \leq a$. This implies that $F$ is fully faithful.

  Given a function $t \in \BInt^k(a,b)$, we can decompose it into 1) the ordered sequence $p_1 < \cdots < p_k$ of values it takes in $\probs$ and 2) a sequence of injections
  \begin{equation}\label{eqn:decomposeleftinjective}
    \begin{tikzcd}
      a & {z_*} & {z_1} & \cdots & {z_k} & b
      \arrow["{{t_*}}"', hook', from=1-2, to=1-1]
      \arrow["{{t_1}}"', hook', from=1-3, to=1-2]
      \arrow["{{t_2}}"', hook', from=1-4, to=1-3]
      \arrow["{{t_k}}"', hook', from=1-5, to=1-4]
      \arrow["t|_{z_k}", hook, from=1-5, to=1-6]
    \end{tikzcd}, \text{ where}
  \end{equation}
  \begin{itemize}
    \item $z_*$ is the preimage of $b+\probs$ under $t$, meaning it contains only those elements whose images are not $\b*$, and $t_*$ is the inclusion,
    \item $z_1$ is the preimage of $b+\probs\setminus \{p_1\}$, meaning it contains only those elements whose images are not $\b*$ nor $p_1$, and $t_1$ is the inclusion,
    \item $z_2$ is the preimage of $b+\probs\setminus \{p_1,p_2\}$, meaning it contains only those elements whose images are not $\b*$ nor $p_1$ nor $p_2$, and $t_2$ is the inclusion,
    \item similarly for each $z_\ell$ and $t_\ell$,
    \item $z_k$ is therefore the preimage of $b$ under $t$, and $t|_{z_k}$ is $t$ restricted to its preimage on $b$, and it is injective by definition of left-injective.
  \end{itemize}

  Conversely, given an ordered sequence $(p_1 < \cdots < p_k)$ in $\probs$ along with a sequence of injections as in \cref{eqn:decomposeleftinjective}, we can reconstruct the function $t: a \rightarrow b+\probstar$ as follows. Given $i \in a$, we consider three possible cases: if $i$ is not in the image of $t_*$, then $t(i) \coloneq \b*$; else if there exists $1 \leq \ell \leq k$ such that $i$ is in the image of $t_*\circ t_1 \circ \cdots \circ t_{\ell-1}$ but not in the image of $t_*\circ t_1 \circ \cdots \circ t_{\ell}$ (with $1\leq \ell \leq k$) then $t(i) \coloneq p_\ell$; otherwise, $i$ is in the image of $t_*\circ t_1 \circ \cdots \circ t_k$, so it has a preimage $i^* \in z_k$, and $t(i) \coloneq t|_{z_k}(i^*)$.

  By appropriate renamings, we can force all the sets in \cref{eqn:decomposeleftinjective} to be of the form $\{0,\dots,n-1\}$ and all the right-to-left injections to be order-preserving. Indeed, each $z_\ell$ is a subset of the set to its left and the function $t_\ell$ is an inclusion, so you can take $n$ to be the cardinality of $z_\ell$ and $t_\ell$ to send $i$ to the $i$th smallest element in $z_\ell$. This makes the mappings described above inverses to each other, simply because there are no other possible choices for the $z_i$s in \cref{eqn:decomposeleftinjective}. Therefore, letting $\OrdInj$ denote the wide subcategory of $\FInj$ with order-preserving injections only, we obtain a bijective correspondence between $\BInt^k(a,b)$ and
  \[\probs^{(k,<)}\times \sum_{0 \leq z \leq a,b} \left(\OrdInj_{k+1}(z,a) \times \FInj(z,b)\right), \text{ where}\]
  \begin{itemize}
    \item $\probs^{(k,<)}$ is the set of strictly increasing sequences of length $k$ in $\probs$,
    \item $\sum_{0 \leq z \leq a,b}$ denotes a dependent sum, corresponding to the choice of $z_k$ in \cref{eqn:decomposeleftinjective},
    \item $\OrdInj_{k+1}(z,a)$ is the set of paths of length $k+1$ from $z$ to $a$ in $\OrdInj$, corresponding to the choice of one injection for each of $t_*, t_1, \dots, t_k$ in \cref{eqn:decomposeleftinjective}, and
    \item $\FInj(z,b)$ is the set of injections from $z$ to $b$.
  \end{itemize}

  Now, we exploit known results about diagrammatic presentations of the two categories above. While $\FInj$ is presented by the \textit{symmetric} "monoidal theory" with generator $\inlinefig{0.75}{genbintunit}$ and no axioms, $\OrdInj$ is presented by the same "theory@@MON" considered \textit{non-symmetric}. Concretely, the diagrams that correspond to morphisms in $\OrdInj$ do not contain any swaps ($\inlinefig{0.5}{swap}$). We note that the presentation of the opposite category is obtained by flipping the generators and axioms, and how we draw the generator is not relevant. Thus, if we consider a presentation of $\op{\OrdInj}$ with generator $\inlinefig{0.75}{genbintcounit}$ and one with generator $\bintp{p_\ell}$ for each $1\leq \ell \leq k$, we get a correspondence between $\BInt^k(a,b)$ and diagrams of shape
  \begin{equation}\label{eqn:diag:decomposeleftinj}
    \input{decomposeleftinj.tikz},
  \end{equation}
  where $d_*$ contains no swaps and only the generator $\inlinefig{0.75}{genbintcounit}$, each $d_\ell$ contains no swaps and only the generator $\bintp{p_\ell}$, and $d_t$ may contain swaps and the generator $\inlinefig{0.75}{genbintunit}$.

  Finally, we observe that any morphism in $\vync{\sigv,\axv}^k(a,b)$ is equal to a diagram of the shape in \cref{eqn:diag:decomposeleftinj}. Indeed, starting from a reduced representative (so that every generator lies on a path starting at an input or an output of the diagram), and using the interchange law \cref{eqn:month:interchange} and naturality of symmetries \cref{eqn:month:overthewire}, we can ensure that the generators $\bintp{p_1}$ up to $\bintp{p_k}$ appear in that order, after all occurrences of $\inlinefig{0.75}{genbintcounit}$, and before any occurrence of $\inlinefig{0.5}{swap}$ or $\inlinefig{0.75}{genbintunit}$. We conclude that $\BInt^k(a,b)$ is in bijective correspondence with $\vync{\sigv,\axv}^k(a,b)$. One can further check that this correspondence is induced by $F$, i.e.\@ the process we described above is inverse to $F$.
\end{proof}

\begin{prop}\label{prop:nonuniversal}
  There are morphisms in $\EPBS$ which are not equal to $\sem{t}$ for any diagram $t \in \sync{\thI}$.
\end{prop}
\begin{proof}
  We only need to look at the specific cases of morphisms $0 \garrow{1} 1$ which correspond to a choice of distributions on $\B$ for each element of $\probstar$.

  A diagram $t: 0 \garrow{1} 1$ in $\sync{\thI}$ must contain exactly one dangling wire on the bottom. It is easy to see that this wire must be connected to either $\inlinefig{0.7}{vgenbintunit}$ or $\inlinefig{0.7}{intergen}$. Then, using the laws of "graded SMCs" (\cref{fig:axiomssyncat}) and axioms \cref{eqn:intergencounit,eqn:intergenp}, we can rewrite $t$ in either of these shapes:
  \[\input{twithunit.tikz} \qquad \text{or} \qquad \input{twithgen.tikz},\]
  where $n \in \N$, and $t_1$ and $t_2$ are diagrams in $\sync{\thCC}$.

  In the first case, $t_1$ represents a distribution $\dist$, and since the intervention is discarded, $\sem{t}$ would send every element in $\probstar$ to $\dist$. In the second case, $t_1$ would represent a distribution on $\B^n\otimes \B^1$ that we write
  \[\sem{t_1} = \sum_{\vec{u} \in \B^n} q^{\b0}_{\vec{u}}\ket{\vec{u}\b0} + q^{\b1}_{\vec{u}}\ket{\vec{u}\b1}.\]
  We can write the distribution after intervening as follows:
  \begin{align*}
    \sem{\scalebox{0.8}{\input{twithgenbeforeint.tikz}}}
     & = \sem{\scalebox{0.8}{\input{twithgenafterint.tikz}}}                                                                                                                                         \\
     & = p\sum_{\vec{u} \in \B^n}(q^{\b0}_{\vec{u}} + q^{\b1}_{\vec{u}})\sem{t_2}\ket{\vec{u}\b0} + (1-p)\sum_{\vec{u} \in \B^n}(q^{\b0}_{\vec{u}} + q^{\b1}_{\vec{u}})\sem{t_2}\ket{\vec{u}\b1},
  \end{align*}
  where the second equality is justified by the semantics of diagrams in $\thCC$.

  Now, we see that $\sem{t}$ cannot be the function $\left\{ \b* \mapsto \ket{\b1}, p \mapsto \ket{\b0} \right\}$. Indeed, if $\sem{t}$ sends $p$ to $\ket{\b0}$, then $\sem{t_2}\ket{\vec{u}\b0}$ and $\sem{t_2}\ket{\vec{u}\b1}$ must be $\ket{\b0}$ for every $\vec{u}$ such that $q^{\b0}_{\vec{u}} + q^{\b1}_{\vec{u}} > 0$. This in turn means that
  \[\sem{t_1\compsimple t_2} = \sum_{\vec{u} \in \B^n} q^{\b0}_{\vec{u}}\sem{t_2}\ket{\vec{u}\b0} + q^{\b1}_{\vec{u}}\sem{t_2}\ket{\vec{u}\b1} = \ket{\b0}.\]
  We identified one function in $\EPBS_1(0,1)$ which is not the semantics of a diagram in $\sync{\thI}$, which is enough to conclude that $\sem{-}$ is not full.
\end{proof}

\end{document}